\documentclass[a4paper,10pt]{article}
\everymath{\displaystyle}
\usepackage{amssymb,amsmath,amsthm}
\usepackage{graphicx}

\usepackage{amsfonts}
\usepackage{hyperref}
\usepackage{latexsym}
\usepackage{amsthm}
\usepackage{amsmath}
\usepackage{amssymb}
\usepackage{amsopn}
\usepackage{geometry}
\usepackage{amscd}
\usepackage{mathrsfs}
\usepackage{graphicx}
\usepackage[english]{babel}
\usepackage[english]{babel}

\newcommand{\Om}{\Omega}
\newcommand{\la}{\langle}
\newcommand{\ra}{\rangle}

\newenvironment{pf}{\noindent{\sc Proof}.\enspace}{\rule{2mm}{2mm}\medskip}
\newenvironment{pfn}{\noindent{\sc Proof} \enspace}{\rule{2mm}{2mm}\medskip}
\newtheorem{theorem}{Theorem}[section]
\newtheorem{proposition}{Proposition}[section]
\newtheorem{lemma}{Lemma}[section]
\newtheorem{corollary}{Corollary}[section]

\newtheorem{remark}{Remark}[section]
\newtheorem{remarks}{Remark}[section]
\newtheorem{definition}{Definition}[section]
\newcommand{\be}{\begin{equation}}
\newcommand{\ee}{\end{equation}}
\newcommand{\teta}{\theta}
\newcommand{\om}{\omega}

\newcommand{\e}{\varepsilon}

\newcommand{\ov}{\overline}
\newcommand{\wtilde}{\widetilde}
\newcommand{\R}{\mathbb R}
\newcommand{\C}{\mathbb C}

\newcommand{\Z}{\mathbb Z}

\newcommand{\N}{\mathbb N}
\newcommand{\T}{\mathbb T}

\renewcommand{\a }{\alpha }
\renewcommand{\b }{\beta }
\newcommand{\s }{\sigma }
\newcommand{\ii }{{\rm i} }
\renewcommand{\d }{\delta }
\newcommand{\D }{\Delta}
\newcommand{\g }{\gamma}

\renewcommand{\l }{\lambda }

\newcommand{\vphi}{\varphi }

\renewcommand{\t }{\tau }
\renewcommand{\o }{\omega }
\renewcommand{\O }{\Omega }
\newcommand{\norma}{|\!\!|}

\newcommand{\lin}{{\cal L}}
\newcommand{\matr}{{\cal M}}
\newcommand{\mapp}{{\cal H}_+}
\newcommand{\norso}[1]{|\!\!| #1  |\!\!|_{s_0}}
\newcommand{\nors}[1]{|\!\!| #1  |\!\!|_{s}}
\newcommand{\norsone}[1]{|\!\!| #1  |\!\!|_{s_1}}
\newcommand{\norS}[1]{|\!\!| #1  |\!\!|_{S}}
\newcommand{\norsob}[1]{|\!\!| #1  |\!\!|_{s_0 + b}}
\newcommand{\norsb}[1]{|\!\!| #1  |\!\!|_{s+ b}}

\newcommand{\dom}{{\cal D}}
\newcommand{\remain}{{\cal R}}

\newcommand{\Hb}{{\bf H}}

\begin{document}

\title{{\bf  Quasi-periodic solutions 
 with Sobolev regularity of NLS on $ \T^d $ with a multiplicative potential }}

\date{}

\author{Massimiliano Berti, Philippe Bolle}

\maketitle

\noindent
{\bf Abstract:}
We prove the existence of  
quasi-periodic solutions  for Schr\"odinger equations with a {\it multiplicative} potential
on $ \T^d $, $ d \geq 1 $, 
merely {\it differentiable} nonlinearities,
 and  tangential frequencies constrained along a {\it pre-assigned}  direction.
The solutions have only Sobolev regularity both in time and space.  
If the nonlinearity and the potential are  $ C^\infty $ then the solutions are  $ C^\infty $.
The proofs are  based on an improved 
Nash-Moser  iterative scheme, which assumes the weakest tame estimates
for the inverse linearized operators (``Green functions") along scales of Sobolev spaces. 
The key  off-diagonal decay estimates of the Green functions are proved via a new multiscale inductive analysis.
The main novelty concerns the measure and ``complexity" estimates.
\\[2mm]
\\[1mm]
2000AMS subject classification: 35Q55, 37K55, 37K50.

\section{Introduction}\setcounter{equation}{0}

The first existence results of 
quasi-periodic solutions of Hamiltonian PDEs  have been proved by Kuksin \cite{K1} and Wayne
\cite{Wa1} 
for one dimensional, analytic, nonlinear perturbations of linear wave and Schr\"odinger equations. 
The main difficulty, namely   the presence  of arbitrarily ``small divisors" in the expansion series of the solutions,
is handled via KAM theory.
These pioneering results 
were limited to 
Dirichlet boundary conditions because the eigenvalues of the Laplacian  had to be simple.
In this case one can impose the so-called  ``second order  Melnikov" non-resonance conditions
to solve the linear 
homological equations which arise at each KAM step, see also P\"oschel \cite{Po2}. 
Such equations are  linear PDEs with constant coefficients and can be solved using Fourier series. 
Already for periodic boundary conditions, where two consecutive eigenvalues are possibly equal, the 
second order  Melnikov non-resonance conditions are  violated.

\smallskip

Later on, another more direct bifurcation approach has been proposed 
by Craig and Wayne \cite{CW}, who introduced the 
Lyapunov-Schmidt decomposition method for PDEs and solved the small divisors problem, for periodic
solutions,  with an analytic Newton iterative scheme. 
The advantage of this approach is to require only the ``first order Melnikov" non-resonance conditions,
which are essentially
the minimal assumptions.  
On the other hand, the main difficulty of this strategy lies 
in the inversion of the linearized operators obtained at each step of the
iteration, and in achieving suitable estimates for their inverses in high (analytic) norms. Indeed
these  operators come from linear PDEs with {\it non-constant} coefficients and 
are small perturbations of a diagonal operator having 
arbitrarily small eigenvalues. 

In order to get estimates  in analytic norms for the inverses, called Green functions by the 
analogy with Anderson localization theory,  
Craig and Wayne developed a coupling technique inspired by the methods of Fr\"ohlich-Spencer \cite{FS}.
The key properties are:
\\[1mm] 
($i$) ``separations" between  singular sites, namely the Fourier indexes of the small divisors,
\\[1mm]
($ii$) ``localization" of the eigenfunctions of 
$ - \partial_{xx} + V(x) $ with respect to the exponentials.
\\[1mm]
Property ($ii$) implies that the matrix which represents, in the eigenfunction basis, 
the multiplication operator for an analytic function has an exponentially fast decay  off the diagonal. 
Then the ``separation properties"  ($i$)  
imply a  very ``weak interaction" between the singular sites.
Property ($ii$) holds in dimension $ 1$,  i.e. $ x \in \T^1 $, 
but, for $ x \in \T^d $, $ d \geq 2 $, 
some counterexamples are known, see \cite{FKT}.

\smallskip

The ``separation properties" ($i$) are quite different for 
periodic or quasi-periodic solutions.
In the first case the singular sites are ``separated at infinity", namely the distance between 
distinct singular sites increases when the Fourier indexes tend to infinity. 
This property is exploited in \cite{CW}. 
On the contrary, it 
never holds for quasi-periodic solutions, even 
for finite dimensional systems. For example, in the ODE 
case where 
the small divisors are $ \omega \cdot k  $, $ k \in \Z^\nu $, 
if the frequency vector $ \om  \in \R^\nu $  is  diophantine, then  the singular sites  $ k  $ where
$ | \omega \cdot k | \leq \rho $ are ``uniformly distributed" in a neighborhood of the hyperplane
$ \omega \cdot k = 0 $, 
with nearby indices at distance $ O( \rho^{- \a}) $ for some $ \a >  0 $. 

This difficulty has been overcome by Bourgain \cite{Bo1}, 
who extended the approach of  Craig-Wayne in \cite{CW} via a  multiscale inductive argument, 
proving the existence of quasi-periodic solutions of $ 1 $-dimensional wave and Schr\"odinger equations 
with polynomial nonlinearities. 
In order to get  estimates of the Green functions,  Bourgain 
imposed  lower bounds for the determinants of most ``singular sub-matrices" along the diagonal. 
This implies, by a  repeated use of  the ``resolvent identity"  (see \cite{FS}, \cite{B3}), 
a sub-exponentially  fast decay of the Green functions. As a 
consequence, at the end of the iteration, the quasi-periodic solutions are Gevrey regular.

\smallskip

At present, KAM theory for $1$-dimensional semilinear PDEs has been sufficiently understood, see e.g. 
\cite{K2},  \cite{KP}, \cite{C}, 
but much work remains for PDEs in higher space dimensions, due to the more complex 
properties of the eigenfunctions and eigenvalues of
$$
( - \Delta + V(x)) \, \psi_j(x) = \mu_j \, \psi_j (x) \, . 
$$
The main difficulties for PDEs in higher dimensions are:
\begin{enumerate}
\item the multiplicity of the eigenvalues $ \mu_j $ 
tends to infinity as $ \mu_j \to + \infty $, 
\item the eigenfunctions $ \psi_j (x) $ are (in general) ``not localized" with respect to the exponentials.
\end{enumerate}

Problem 2 has been often bypassed considering pseudo-differential PDEs substituting the 
multiplicative potential $ V(x) $ by a ``convolution potential"  
$$
V * ( e^{\ii j \cdot x}) = m_j  e^{\ii j \cdot x} \, , \ m_j \in \R \, ,   \ j \in \Z^d \, , 
$$
which, by definition, is diagonal on the exponentials. 
The scalars $ m_j $ are called the ``Fourier multipliers".

\smallskip

Concerning problem 1, since the  approach of Craig-Wayne and Bourgain requires 
only the  first order Melnikov non-resonance conditions, 
it works well, in principle,
in case of multiple eigenvalues, in  particular for  PDEs in higher spatial dimensions. 

Actually the first existence results of periodic solutions  for NLW and NLS on $ \T^d $, $ d \geq 2 $, 
have been established   
by Bourgain in \cite{B4}-\cite{B3}. 
Here the singular sites 
form huge clusters (not only points as in $ d = 1 $) but are still ``separated at infinity".
The nonlinearities are polynomial  and the solutions have Gevrey regularity in space and time. 

\smallskip

Recently these results were extended in  \cite{BB07}-\cite{BP} 
to prove the existence of periodic solutions, with only Sobolev regularity, for NLS and NLW
in any dimension and with merely 
differentiable nonlinearities. 
Actually in \cite{BBP}, \cite{BP} the  PDEs are defined not only on tori, but on any
compact Zoll manifold, Lie group and homogeneous space. 
These results are proved via  an abstract Nash-Moser implicit function theorem (a simple Newton method is not sufficient).
Clearly, a  difficulty when working with functions having only Sobolev regularity is that the Green functions will 
exhibit only a 
polynomial decay off the diagonal, and not exponential (or sub-exponential).
A key concept that one must exploit are the interpolation/tame estimates. 
For PDEs on Lie groups only weak properties of   ``localization"  ($ii$) of the eigenfunctions hold, see \cite{BP}. 
Nevertheless these properties imply a block diagonal decay, for the matrix which represents the 
multiplication operator in the eigenfunctions basis, sufficient to achieve the tame estimates.   

We also mention that existence of periodic solutions for NLS on $ \T^d $ has been proved, for analytic nonlinearities,
by Gentile-Procesi \cite{GP} via the Lindstedt series techniques,  and,  in the differentiable case, by Delort \cite{D}
using paradifferential calculus. 

\smallskip

Regarding quasi-periodic solutions,  Bourgain \cite{B3} was  the first to prove  their existence 
for PDEs in higher dimension, actually 
for nonlinear Schr\"odinger equations with Fourier multipliers and  polynomial nonlinearities 
on $ \T^d $ with $ d = 2 $. The Fourier multipliers, in number equal to  the tangential frequencies 
of the quasi-periodic solution, play the role of external parameters.
The main difficulty arises in the multiscale argument to estimate the decay of the Green functions.
Due to the degeneracy of the eigenvalues of the Laplacian the singular sub-matrices 
that one has to control are huge. 
If $ d = 2 $, careful estimates on the number of integer
vectors on a sphere, allowed anyway Bourgain to show that
the required non-resonance conditions are fulfilled for ``most"
Fourier multipliers.

More recently Bourgain \cite{B5} improved the techniques in \cite{B3} 
proving the existence of quasi-periodic solutions for  nonlinear 
wave and Schr\"odinger equations with Fourier multipliers on any $ \T^d $, $ d > 2 $, still for 
polynomial nonlinearities.
The improvement  in \cite{B5} comes from the use
of sophisticated techniques developed in the context of Anderson localization theory 
in Bourgain-Goldstein-Schlag \cite{BGS}, Bourgain \cite{B6}, 
see also Bourgain-Wang \cite{BW1}. 
These techniques (sub-harmonic functions, Cartan theorem, semi-algebraic sets) mainly concern fine properties
of rational and analytic functions, especially measure estimates of sublevels. 
 Actually the nonlinearities in \cite{B5} are  taken to be polynomials 
in order to use semialgebraic  techniques. Very recently, Wang  \cite{W2} has generalized
the results in \cite{B5} for NLS  with no  Fourier multipliers and with supercritical nonlinearities. 
The main step is a Lyapunov-Schmidt reduction in order  to 
 introduce parameters and then be able to apply the results of \cite{B5}.

\smallskip

We also remark that, 
in the last years,   the KAM approach has been extended by  Eliasson-Kuksin \cite{EK} for 
nonlinear Schr\"odinger equations on $ \T^d $
with a convolution potential 
 and analytic nonlinearities. 
The potential plays the role of ``external parameters". The quasi-periodic solutions are $ C^\infty $ in space.  
Clearly an advantage of the KAM approach 
is to provide also a stability result:
the linearized equations on the perturbed invariant tori  are reducible to constant coefficients,
see also 
 \cite{EK1}. 
 
For the cubic NLS in $ d = 2 $  the existence of quasi-periodic solutions 
has been recently proved by Geng-Xu-You \cite{GXY}
via a Birkhoff normal form and a modification of  the KAM approach in \cite{EK}, see also 
 Procesi-Procesi \cite{PP}, valid in any dimension.

\smallskip

In the present paper we prove  -see Theorem \ref{thm:main}-
the existence of quasi-periodic solutions
for nonlinear Schr\"odinger equations on $ \T^d $, $ d \geq 1 $, with:
\begin{enumerate}
\item 
merely {\it differentiable} nonlinearities, 
see \eqref{nonli},
\item a {\it multiplicative} (merely differentiable) potential $ V(x) $, see \eqref{eq:posi},
\item 
a {\it pre-assigned} (Diophantine) direction of the tangential frequencies, 
see \eqref{baromega}-\eqref{diophan0} . 
\end{enumerate}
The quasi-periodic solutions in  Theorem \ref{thm:main}
have the same Sobolev regularity both in time and space, see remark \ref{regus}.
Moreover, we prove that, if the potential and the nonlinearity are of class $ C^\infty $, then the 
quasi-periodic solutions are $ C^\infty $-functions of $(t,x)$. 

\smallskip
Let us make some comments on the results. 

\smallskip

1. Theorem \ref{thm:main} 
confirms the natural conjecture about the persistence of quasi-periodic solutions  for Hamiltonian PDEs 
into a setting of finitely many derivatives (as in the classical KAM theory \cite{MP}, \cite{Po}, \cite{Z}), 
stated for example by Bourgain \cite{B97}, page 97.
The nonlinearities in Theorem \ref{thm:main}, as well as the potential, are sufficiently many times differentiable,
depending on the dimension and the number of  the frequencies. 
Of course we can not expect the existence of quasi-periodic solutions of the Schr\"odinger equation under too weak regularity
assumptions on the nonlinearities. Actually, 
for finite dimensional Hamiltonian systems, it has been rigorously proved that,  
if the vector field is not sufficiently smooth, then all the invariant tori 
could be destroyed and only  discontinuous  Aubry-Mather invariant sets survive, see e.g. \cite{He1}.
We have not tried to estimate 
the minimal 
smoothness exponents, see however remark \ref{smoothness}.
This could be interesting for comparing  Theorem \ref{thm:main} 
with the well posedness results of the Cauchy problem.

\smallskip

2.  Theorem  \ref{thm:main} is the first existence result of quasi-periodic solutions 
with a multiplicative potential $ V(x) $ on $ \T^d $, $ d \geq 2 $.
We never exploit properties of ``localizations" of the eigenfunctions  of $ - \Delta + V(x) $ 
with respect to the exponentials, that  
actually might {\it not} be true, 
see \cite{FKT}.
Along the multiscale analysis we use the exponential basis
which diagonalizes $ - \Delta + m $ where $ m $ is the average of $ V(x) $, see \eqref{average}, 
and not the eigenfunctions of $ - \Delta + V(x) $.
 In \cite{B3} Bourgain considered analytic multiplicative periodic potentials of the special form
 $ V_1(x_1) + \ldots + V_d(x_d) $ to 
 ensure localization properties of the eigenfunctions,  
 leaving open the natural problem for a general  multiplicative potential $ V(x) $.

We also underline that  Theorem \ref{thm:main} holds for any {\it fixed} potential $ V(x) $: 
 we do  not extract parameters from $ V $,  the role of external parameters being played by 
 the frequency $ \om = \l \bar \omega $. 

\smallskip

3. For finite dimensional systems, 
the existence of quasi-periodic solutions with tangential frequencies constrained along a fixed direction
has been proved by Eliasson \cite{E} (with KAM theory) and Bourgain \cite{B1} (with a multiscale approach).
The main difficulty clearly relies in satisfying the  Melnikov non-resonance conditions, required at each step of the iterative process, 
using only {\it one} parameter. Bourgain raised in  \cite{B1} the question if a similar result
holds true also for infinite dimensional Hamiltonian systems.
This has been recently proved  in \cite{BBi} for $ 1 $-dimensional PDEs, 
verifying the second order Melnikov non-resonance conditions of KAM theory.
Theorem \ref{thm:main}   (and its method of proof) answers positively to Bourgain's conjecture
also for PDEs in higher space dimension. The non-resonance
conditions that we have to fulfill are of first order Melnikov type, see the end of section \ref{ideas}.

\smallskip

The proof of Theorem \ref{thm:main} is based on a 
Nash-Moser iterative scheme and a multiscale analysis of the  linearized operators as in \cite{B5}.
However, our approach presents many differences with respect to  Bourgain's one  \cite{B5}, about:
\begin{enumerate}
\item the  iterative scheme, 
\item the multiscale proof of  the Green's functions polynomial decay estimates. 
\end{enumerate}

Referring to section \ref{ideas} for a detailed exposition of our approach, we
outline here the main differences.

\smallskip

1. Since we deal with merely differentiable nonlinearities 
we need all the power of the Nash-Moser theory  in scales of Sobolev functions spaces.
A Newton method valid in analytic Banach  scales is not sufficient. 
This means  that the superexponential smallness of the error terms  
due to  finite dimensional truncations, see \eqref{rnhigh}, 
can  not be obtained, in Sobolev scales, decreasing the analyticity strips, but
using the structure of the iteration and the interpolation   estimates  of the   Green functions, 
see lemmas \ref{lemcon}, \ref{htn+1}, \ref{S5n+1}. This is a key idea 
when  dealing with matrices with  a merely polynomial off-diagonal decay. 

Actually, the Nash-Moser scheme developed in section \ref{sec:NM}
also improves the one  in \cite{BB07}-\cite{BBP}, 
requiring  the minimal tame properties  \eqref{normabassa}
 for the inverse linearized operators, see comments after (\ref{Ln1}).

Another comment is in order: we do not follow the 
 ``analytic smoothing technique" suggested by Moser in \cite{MP}
 of approximating  the differentiable Hamiltonian PDE by analytic ones.
This  technique is very efficient for finite dimensional Hamiltonian systems, see \cite{Po}, \cite{Z}, 
but  it seems quite delicate for PDEs (especially 
in dimensions $ d \geq 2 $) because of the presence of large clusters of small divisors.
So we prefer a more direct Nash-Moser iterative procedure more similar,  in spirit, to \cite{LZ}.
 
\smallskip

2.  The main difference between our multiscale  approach, which is developed
 to prove the Green functions  estimates \eqref{normabassa}, 
 and the one in \cite{B5}, \cite{BGS}, \cite{B6}, \cite{BW1}, 
concerns the way we prove inductively the existence of ``large sets" of  
$ N_n $-good parameters, see Definition \ref{def:freqgood}.
Quoting Bourgain \cite{B04} 
``...{\it the results in \cite{B5} make essential use of the general perturbative technology
(based on subharmonicity and 
semi-algebraic set theory) [...]. This technique enables us to deal with large sets of `singular sites' 
[...], something difficult to achieve  with 
conventional eigenvalue methods.}". 
Actually,  exploiting  that $ - \Delta + V(x) $ is positive definite,
 we are able  to prove the necessary measure and ``complexity" estimates by
using only elementary eigenvalue variation arguments, 
 see section \ref{sec:measure}.

\smallskip

Another deep difference is required for dealing with a multiplicative potential $ V(x) $:
we define  ``very regular" sites (see Definition \ref{regulars}) depending on the potential $ V $.

\smallskip

We hope that this novel approach will be useful
also for extending the results of \cite{B6}, \cite{B5}, \cite{BGS}, \cite{BW1}.

\smallskip

We tried to present the steps of proof in an abstract setting (as much as possible) 
in order to develop 
a systematic procedure, alternative to KAM theory, for the search of quasi-periodic solutions of PDEs. 
The proof of Theorem \ref{thm:main} is completely self-contained.
All the techniques employed are elementary and 
based on  abstract arguments valid for many PDEs. 
Only the ``separation properties" of the bad sites (section \ref{sec:sepa}) will change, of course, for different PDEs.

\smallskip

Since the aim of the present paper is to focus on the small divisors problem for 
quasi-periodic solutions with Sobolev regularity of NLS with a multiplicative potential on $ \T^d $  
and differentiable nonlinearities, 
we have considered, among many possible variations,
quasi-periodically  forced nonlinear perturbations of linear Schr\"odinger equations.
In this way, we avoid the Lyapunov-Schmidt decomposition. 
Clearly the small divisors difficulty for quasi-periodically forced NLS is the same as for autonomous NLS.

\smallskip

We now state precisely our results. 

\subsection{Main result}\setcounter{equation}{0}

We consider $ d $-dimensional nonlinear  Schr\"odinger equations with a potential $V$,  like
\be\label{eq:main} 
\ii u_{t} - \Delta u + V(x) u  = \e f(\om t, x, |u|^2 )u  + \e g( \om t, x)    \, , \quad x \in {\T}^d  \, , 
\ee 
where   
$ V \in C^{q} ({\T}^d;\R) $ for some $ q $ large  enough, 
$ \e > 0 $ is a small parameter, 
the frequency vector $ \om \in \R^\nu $ is non resonant (see \eqref{diophan0}),   
the nonlinearity is  quasi-periodic in time 
and only finitely many times differentiable, more precisely 
\be\label{nonli}
f \in C^q ( {\T}^\nu \times {\T}^d \times \R;\R)  \, , \quad
g \in C^q ( {\T}^\nu \times {\T}^d; {\C}) 
\ee
for some  $ q \in \N $ large enough. 
Moreover we suppose
\be\label{eq:posi}
- \Delta + V(x) \geq  \b_0 I  \, ,  \  \b_0 >  0 \, . 
\ee

\begin{remark}\label{rem:posi}
Condition \eqref{eq:posi} is  used for the measure estimates of section \ref{sec:measure}. Actually 
for autonomous NLS it can be always verified after 
a gauge-transformation  $ u \mapsto e^{- \ii \s t} u $  
for $ \s $ large enough.
\end{remark}

We assume that the frequency vector $ \om $ is a small dilatation  of a fixed
Diophantine  vector $ \bar \o \in \R^{\nu}  $, namely
\be\label{baromega}
\om = \l \bar \om \, , \quad \l \in \Lambda := [1/2, 3/2]  \, ,  \quad |\bar \om | \leq 1 \, , 
\ee
where,  for some $  \g_0 \in (0,1) $, $ \t_0  > \nu - 1 $,  
\be\label{diophan0} 
|\bar \om \cdot l | \geq \frac{2 \g_0}{|l|^{\tau_0}} \, , \quad \forall l \in \Z^{\nu} \setminus \{ 0 \} \, , 
\ee
and $ |l| := \max \{|l_1|, \ldots , |l_\nu| \} $. 
For definiteness we fix $ \t_0 := \nu $.

\smallskip
 
If $ g (\om t,x) \not\equiv 0 $ then $ u = 0 $ is not a solution of (\ref{eq:main}) for $ \e \neq  0 $.

\begin{itemize}
\item Question: do there exist quasi-periodic solutions of (\ref{eq:main}) for  sets of 
parameters $ (\e , \l ) $ of positive measure?
\end{itemize}

This means looking for $(2\pi)^{\nu +d}$-periodic  solutions $ u(\vphi,x) $ of
\be\label{eq:freq} 
\ii \om \cdot \partial_\vphi u  - \Delta u + V(x) u = \e f( \vphi , x, |u|^2 )u  + \e g(\vphi,  x) \, . 
\ee 
These solutions will be, for some $  (\nu + d) \slash 2 < s \leq q $,
in the Sobolev space
\begin{eqnarray}\label{def:Hs}
H^s  :=  H^s ( \T^\nu \times \T^d;  {\C}) & := & \Big\{  u(\vphi,x)  =  \sum_{ (l,j) \in \Z^\nu \times \Z^d}  u_{l,j} 
e^{\ii (l \cdot \vphi + j \cdot x)} \\
& & \quad   {\rm such \  \ that \  } \ \| u \|_s^2 := K_0 \sum_{i \in \Z^{\nu + d}}  |u_{i}|^2 \langle i \rangle^{2s} < + \infty \nonumber
\Big\} 
\end{eqnarray}
where 
$$
i := (l,j) \, ,  \quad  \langle i \rangle := \max(|l|,|j|,1)
 \, ,  \  \ |j| := \max \{|j_1|, \ldots , |j_d| \}.
$$
For the sequel we fix 
$ s_0 > (d+\nu )\slash 2 $ so that there is the continuous embedding    
\be\label{embedding}
H^s (\T^{\nu+d} ) \hookrightarrow  L^\infty (\T^{\nu+d} ) \, , \quad \forall s \geq s_0 \, , 
\ee
and $ H^s $ is a Banach algebra with respect to the multiplication of functions.
The constant $ K_0 >  0  $ in the definition  \eqref{def:Hs} of the Sobolev norm $ \| \ \|_s $ is independent of $ s $.
The value of $ K_0 $ is 
fixed (large enough) so that  $ |u|_{L^\infty}  \leq \| u \|_{s_0}$ and the interpolation inequality 
\be\label{lions}
\| u_1 u_2 \|_s \leq \frac12 \| u_1 \|_{s_0} \| u_2 \|_s + \frac{C(s)}{2}  \| u_1 \|_s \| u_2 \|_{s_0}   \, , \quad 
\forall s \geq s_0 \, , \ u_1, u_2 \in H^s \, , 
\ee
holds with $C(s) \geq 1$ and   $ C(s) = 1, \forall s \in [s_0,s_1] $; the 
constant $ s_1 $ is defined in \eqref{Sgr} and 
depends only on $d,\nu, \tau_0 := \nu $. 
With respect to the standard Moser-Nirenberg interpolation estimate in Sobolev spaces, see e.g. 
\cite{LM}, the additional property in (\ref{lions})  is that one of the constants  is independent
of $ s $. The proof  of (\ref{lions})  is given for example  in Appendix of \cite{BBP}, see also \cite{LM}. 

\smallskip

The main result of this paper is:

\begin{theorem}  \label{thm:main} Assume (\ref{diophan0}). 
There are $ s := s(d, \nu) $,  $ q := q(d, \nu) \in \N $, such that:
$ \forall \, V \in C^q $ satisfying \eqref{eq:posi}, $ \forall f , g  \in C^q $, there exist $ \varepsilon_0 > 0 $,  a  map 
$$ 
u \in C^1([0, \varepsilon_0] \times \Lambda ;H^s) 
\quad {\rm with }  \quad u (0, \l) = 0  \, , 
$$
and a Cantor like set $ {\cal C}_\infty \subset  [0, \varepsilon_0] \times \Lambda $ of 
asymptotically full Lebesgue measure, i.e.   
\be\label{Cmeas}
| {\cal C}_\infty | / \e_0 \to 1 \quad {\rm as } \quad \e_0 \to 0  \, ,  
\ee
such that, $ \forall (\e,\l) \in {\cal C}_\infty $, $ u (\e, \l)  $ is a solution of (\ref{eq:freq}) with $ \om = \l \bar \om $.
\\[1mm]
Moreover, if $ V, f, g $ are of class $ C^\infty $ then $ u(\e,\l) \in C^\infty (\T^d \times \T^\nu ; \C )$. 
\end{theorem}

We have not tried to optimize the estimates for $ q := q(d, \nu) $ and 
$ s := s(d, \nu) $. 

\begin{remark}\label{smoothness}
In \cite{BB07} we proved the existence of periodic solutions in 
$ H^s_t (\T ;  H^1_x (\T^d)) $ with $ s > 1 /2 $,  for  one dimensional NLW  equations
with  nonlinearities of class $ C^6 $, see the bounds (1.9), (4.28) in \cite{BB07}.
\end{remark}

\subsection{Ideas of the proof} \label{ideas}

{\bf Vector NLS.} We prove Theorem \ref{thm:main}
finding solutions of the ``vector" NLS equation
\be\label{vnls}\left\{\begin{array}{ll}
\ \ \ii \o \cdot \partial_\vphi u^+ - \Delta u^+ + V(x) u^+ =  \e f(\vphi , x, u^- u^+ ) u^+  + \e g( \vphi , x)  \\ 
- \ii \o  \cdot  \partial_\vphi u^- - \Delta u^-   + V(x) u^- =  \e f(\vphi , x, u^- u^+ ) u^-  + \e {\bar g}( \vphi , x)  
\end{array} \right.
\ee
where  
\be\label{bHs}
{\bf u} := (u^+,u^-) \in {\bf H}^s:= H^s \times H^s 
\ee
(the second equation is obtained  by formal  complex conjugation of the first one).
In the system (\ref{vnls}) the variables $ u^+ $, $ u^- $ are independent. 
However,  note that (\ref{vnls}) reduces to the scalar NLS equation (\ref{eq:main}) in the  set 
\be\label{subvs} 
{\cal U} :=  \Big\{ {\bf u} := (u^+,u^-) \   : \  \overline{ u^+} = u^- \Big\} 
\ee
in which $ u^- $ is  the complex conjugate of $ u^+ $ (and viceversa). 
\\[1mm]
{\bf Linearized equations.}
We  look for solutions of the vector NLS equation (\ref{vnls}) in $ {\bf H}^s \cap {\cal U} $ by a Nash-Moser iterative scheme.
The main step concerns the invertibility  of (any finite dimensional restriction of) the linearized operators at 
any $  {\bf u} \in {\bf H}^s \cap {\cal U} $, namely 
$$ 
{\cal L} ({\bf u}) := L_\om - \e T_1 = D_\om +  T 
$$
described in (\ref{Linve})-(\ref{Lu}), with suitable estimates of the inverse in high Sobolev norm.

An advantage of the vector NLS formulation, with respect to the scalar NLS equation \eqref{eq:freq}, 
is that the  operators $ {\cal L} ({\bf u}) $ are $\C$-linear and  selfadjoint.
This is convenient for proving the measure estimates via eigenvalue variation arguments.  
Moreover the matrix $ T $  is T\"oplitz, see \eqref{Tmatrix12}, 
and its entries on the lines parallel to the diagonal decay to zero at a polynomial rate. 
\\[1mm]
{\bf Matrices with off-diagonal decay.}
In section \ref{sec:off} we develop an abstract setting for dealing  with 
matrices with polynomial off-diagonal decay. 
In Definition \ref{defnormatr}   we introduce the $ s $-norm of a matrix 
and we prove the algebra and interpolation properties  \eqref{algebra}, \eqref{interpm}.
The  $ s $-norms are inspired to mimic the behavior 
of matrices representing  the multiplication operator by a function of $ H^s $.
This 
intrinsic setting is very convenient  (in particular for the multiscale Proposition \ref{propinv}) 
to estimate the decay of inverse matrices 
via Neumann series, because
product, and then powers, of matrices with finite $s$-norm 
 will  exhibit the same  off-diagonal decay. 
\\[1mm]
\noindent
{\bf Improved Nash-Moser iteration.} 
We construct inductively better and better  approximate solutions $ {\bf u}_n  $ of the NLS equation  (\ref{vnls}) 
by a Nash-Moser iterative scheme, see the ``truncated" equations $(P_n) $ in Theorem \ref{cor1}.
The $ {\bf u}_n \in H_n $, see \eqref{Hn},
are trigonometric polynomials  with a super-exponential number $ N_n $ of harmonics, see \eqref{defNn}.

At each step we impose that, 
for ``most" parameters $(\e,\l) \in [0,\e_0) \times [1/2,3/2] $, the eigenvalues of the
restricted linearized operators $ {\cal L}_n := P_n {\cal L}({\bf u}_n)_{|H_n} $ 
are in modulus bounded from below by  $ O(N_n^{-\tau}) $, see Lemma \ref{measure0}.
The proof exploits that $ - \Delta + V $ is positive definite, see  \eqref{eq:posi} and remark \ref{rem:posi}. 
Then  the  $ L^2 $-norm 
of the inverse  satisfies $ \| {{\cal L}_n^{-1}} \|_0= O(N_n^\tau) $. 
By Lemma \ref{norcomp} this implies  that the $s$-norm (see Definition \ref{defnormatr}) 
satisfies
$$
\nors{ {\cal L}_n^{-1}} 
\leq  N_n^{s+d+\nu} \| { {\cal L}_n^{-1}} \|_0 = 
O( N_n^{s+d+\nu+\tau}) \, , \  \forall s >  0 \, .
$$ 
Such an estimate is {\it not} sufficient for the convergence of the Nash-Moser scheme.
We need  sharper estimates for the Green functions (sublinear decay), of the form
\be\label{Ln1}
\nors{ {\cal L}_n^{-1}} = O( N_n^{\tau'+\d s}) \, , \   \ \d \in (0,1) \, ,  \ \t' >  0 \, , \  \forall s >  0 \, ,
\ee
which imply an off-diagonal decay of the inverse matrix coefficients like
$$
| ({\cal L}_n^{-1})_{i'}^{i}| \leq C \frac{N_n^{\tau'+\d s}}{ \langle i-i' \rangle^s} \, ,  \quad |i|, |i'| \leq N_n \, ,
$$
see Definition \ref{defM}. Actually the conditions \eqref{Ln1} are optimal
for the convergence of the Nash-Moser iterative scheme, 
as a famous  counter-example of Lojasiewicz-Zehnder \cite{LZ} shows:
if $ \d = 1 $ the scheme does not converge. 
By Lemma \ref{sobonorm} the bound (\ref{Ln1}) implies the interpolation estimate in Sobolev norms 
$$
\|{\cal L}_n^{-1} h \|_s \leq  C(s)( N_n^{\tau'+\d s} \|h \|_{s_1} + N_n^{\tau'+\d s_1} \| h \|_s) \, , \quad \forall s \geq s_1 \, ,
$$
which is sufficient for the  Nash-Moser  convergence.  
Note that the exponent $ \t' + \d s $ in  (\ref{Ln1}) grows with $ s $, unlike
the usual Nash-Moser theory, see e.g. \cite{Z}, where the ``tame" exponents are $ s $-independent. 
Actually it is easier to prove these weaker tame estimates, see, in particular, Step II of Lemma \ref{defY}.

In order to prove \eqref{Ln1}  we have to exploit (mild) ``separation properties" of the small divisors:
several  eigenvalues of ${\cal L}_n $ are  actually much bigger (in modulus) than   $ N_n^{-\tau} $.

\smallskip

\noindent
{\bf Estimates of Green functions.}
The core of the paper is to establish the  
Green functions  estimates (\ref{Ln1}) at each step of the iteration, see Lemma \ref{invLn+1}. 
These follow by an inductive application of the multiscale  Proposition \ref{propinv}, 
once verified the ``separation property" (H3),   see Lemma \ref{H1H3}.

The ``separation properties" of the $ N_n $-bad and singular sites are obtained by  Proposition \ref{prop:separation} 
for all the parameters $(\e,\l) $ 
which are $ N_n $-good, see Definition \ref{def:freqgood} and  assumption ($i$). 
We first use the covariance property \eqref{shifted}
and the ``complexity" information  \eqref{BNcomponents}   on the set $ B_{N}(j_0; \e, \l ) $ in \eqref{tetabad} 
(the set of the ``bad" $\teta$) 
to bound the number of ``bad" time-Fourier components,
see Lemma \ref{Ntime} (this idea goes back to \cite{E1}).
Next  we use also the information that the sites are ``singular" to bound the length of a ``chain" of 
$ N_n $-bad and singular sites (with ideas similar to \cite{B5}), see Lemma \ref{thm:separation}.  

In order to conclude the inductive proof  we have to verify  that  ``most"   parameters $(\e,\l) $ are $ N_n $-good.
For this, we do not  invoke 
sub-harmonic functions theory, Cartan theorem as in
 \cite{B5}, \cite{BGS}, \cite{B6}.

\smallskip

\noindent
{\bf Measure and ``complexity" estimates.}
Using Proposition \ref{PNmeas} we prove first 
that most parameters $(\e,\l) $ are $ N_n $-good in a weak sense.
The proof of Proposition \ref{PNmeas} is based on simple  eigenvalue variation arguments 
and Fubini theorem. 
The main novelty is to use that $ -\Delta + V(x) $ is positive definite, 
see  \eqref{eq:posi} and remark \ref{rem:posi},
and to perform the measure estimates in the new set of variables  \eqref{changevaria}. 
In this way we prove that for ``most" parameters $(\e,\l) $ the set 
$ B_{N}^0(j_0; \e, \l ) $ in \eqref{tetabadweak} (of ``strongly" bad $ \teta $) has a small measure. 
This fact and the Lipschitz dependence of the eigenvalues
with respect to parameters imply also the complexity bound \eqref{BNcomponent2}, see
Lemma \ref{complessita}.   Finally, 
using again the multiscale Proposition \ref{propinv} and the 
separation Proposition \ref{prop:separation} we  conclude inductively that most of these parameters $(\e,\l) $
 are actually $ N_n $-good (in the strong sense), see Lemma \ref{S3n+1}.

\smallskip

\noindent
{\bf Definition of regular sites.} In order to deal with a multiplicative potential the key idea is to define
``very regular" sites, i.e. in Definition \ref{regulars} the constant $ \Theta $ will be taken large
with respect to the potential $ V $,  so that the diagonal terms \eqref{diao}
dominate also the off diagonal part $ V_0(x) $ of the potential, see Lemma \ref{defmatrMN}. Taking
a large value for the constant $ \Theta $ does not affect the qualitative properties of the chains of
singular sites proved in Lemma \ref{thm:separation}. 
Then we achieve in section \ref{sec:sepa}
the separation properties for the clusters of small divisors, and the multiscale Proposition 
\ref{propinv} applies.  We refer  also to Lemmas \ref{Inizioind} and \ref{alta} for a similar 
construction at the initial step of the iteration.

\smallskip

\noindent
{\bf Melnikov non-resonance conditions.}\label{nonreso}  
An advantage of the Nash-Moser iterative scheme is to require  weaker non-resonance conditions
than for the  KAM approach. For clarity we collect all the non-resonance conditions 
that we  make along the paper below:
\\[0.3mm]
\indent
- $ \omega = \l \bar \om $ is diophantine, see (\ref{diophan0}), (\ref{diophan}). It is used only in Lemma \ref{Ntime}
to get separation properties of the bad sites in the time Fourier components. 
\\[0.3mm]
\indent
- $ \omega =\l \bar \om $ satisfies the 
 non-resonance condition (\ref{diofs}) of first order Melnikov type. 
Physically, this assumption means that the forcing frequencies $ \o $ do not enter in resonance with 
the first $ N_0 $ normal mode frequencies of the  linearized Schr\"odinger equation at the origin. 
This is 
used for the initialization of the Nash-Moser scheme, see subsection \ref{sec:ini}.
\\[0.3mm]
\indent
- $ (\l \bar \o, \e) $ satisfy the ``first order Melnikov" non-resonance conditions  at each step of the
Nash-Moser iteration:  
the eigenvalues of $ A_{N_{n}}(\l \bar \om,\e) $ have to be $ \geq 2 N_{n}^{-\t} $, see also
Lemma \ref{measure0}. 
\\[0.5mm]
\indent
- We also verify that most frequencies are $ N $-good (see Definition \ref{def:freqgood}) imposing conditions on 
the eigenvalues of the matrices $ A_{N,j_0}(\l \bar \om, \e, \theta) $ 
as in Lemma \ref{cor2}. These requirements can then be seen 
as other ``first order Melnikov" non-resonance conditions. 

\smallskip

\noindent
{\bf Sobolev regularities.} Along the proof we make use of  three different Sobolev regularity thresholds
$$ 
s_0 < s_1 < S \, .
$$
The scale $ s_0 > (d + \nu )/2 $ is simply required to establish the algebra and interpolation estimates,
see e.g. \eqref{lions}. The Sobolev index $ s_1 $ is large enough to have a sufficiently strong decay when
proving the multiscale  Proposition \ref{propinv}, see \eqref{s1}.
Finally the Sobolev regularity $ S $ is large enough (see \eqref{Sgr}) for proving the convergence of the 
Nash-Moser iterative scheme in section \ref{sec:NM}. 

\smallskip  
  
\noindent
{\bf Acknowledgments:} The authors thank Luca Biasco for useful comments on the paper.

\section{The linearized equation}\setcounter{equation}{0}

We  look for solutions of the vector NLS equation (\ref{vnls}) in $ {\bf H}^s \cap {\cal U} $ (see \eqref{subvs}) 
by a Nash-Moser iterative scheme.
The main step concerns the invertibility  of (any finite dimensional restriction of) the family of linearized operators 
\be\label{Linve}
{\cal L} ({\bf u}) := {\cal L}(\o, \e,  {\bf u}) := L_\om - \e T_1 
\ee
acting on $ {\bf H}^s $,  where   
$ {\bf u} =  (u^+,u^-) \in C^1([0,\e_0] \times \Lambda, {\bf H}^s \cap {\cal U})  $, 
\be\label{LomegaV}
L_{\om} := 
  \left(\begin{array}{cc} 
\ii \o \cdot  \partial_\vphi -   \Delta + V(x) & 0\\ 0 &
- \ii  \om \cdot \partial_\vphi -  \Delta + V(x)  \end{array}\right) 
\ee
and 
\be\label{T1}
T_1 := \left(\begin{array}{cc} p(\vphi,x) & q(\vphi,x)  \\ 
\bar q(\vphi,x) & p(\vphi,x)  \end{array}\right)
\ee
with
\be\label{pq}
p(\vphi,x) :=  f(\vphi, x, |u^+|^2) + f'(\vphi,x, |u^+|^2) |u^+|^2  \, , 
\ q(\vphi,x) :=   f'(\vphi, x, |u^+|^2) (u^+)^2 \, . 
\ee
Above  $ f' $ denotes  the derivative of $f(\vphi,x, s)$ with respect to $ s $. 
The functions $ p, q $ depend also on $ \e, \l $ through $ {\bf u} $.
Note that  $ u^+ u^- = |u^+|^2 \in \R $ since $ {\bf u} \in {\cal U} $, see (\ref{subvs}).

Decomposing the multiplicative potential
\be\label{average}
V(x) = m + V_0 (x)
\ee
where $ m $ is the average of $ V(x) $ and $ V_0 (x) $ has zero mean value,
we also write 
\be\label{LomDom}
L_\om = D_\om +  T_2 
\ee
where $ D_\om $ is the constant coefficient differential operator 
\be\label{Lomega}
D_\om :=  \left(\begin{array}{cc} 
\ii \o \cdot  \partial_\vphi -   \Delta + m & 0\\ 0 &
- \ii  \om \cdot \partial_\vphi -  \Delta + m \end{array}\right)  
\quad {\rm and} \quad  
T_2 := \left(\begin{array}{cc} V_0(x) & 0  \\ 
0 & V_0(x)  \end{array}\right) \, .
\ee
Hence the operator $ {\cal L}({\bf u})$ in \eqref{Linve} can also be written as
\be\label{Lu}
{\cal L}({\bf u}) =   D_\om + T \, , \qquad T := T_2 - \e T_1  \, . 
\ee

\begin{lemma}
$ {\cal L } ({\bf u})$ is {\it symmetric} in $ {\bf H}^0 $, i.e. 
$  ({\cal L } ({\bf u})h,k)_{L^2} = (h,{\cal L } ({\bf u})k)_{L^2} $ for all
$ h, k $ in the domain of ${\cal L } ({\bf u})$. 
\end{lemma}

\begin{pf}
The operator $ L_\om $ is symmetric with respect to the $ L^2 $-scalar product  in $ {\bf H}^0 $, because
  each $ \pm \ii \omega \cdot \partial_\vphi - \Delta + V(x) $  is symmetric in
$  H^0 (\T^\nu \times \T^d; \C ) $. Moreover $ T_2 $, $ T_1 $ are selfadjoint in  $ {\bf H}^0 $
because  $ V(x) $ and $ p(\vphi,x)  $ are real valued, being   $ |u^+|^2 \in \R $ and $ f $ real by (\ref{nonli}), see \cite{BP}. 
\end{pf}


The Fourier basis diagonalizes the differential operator $ D_\om $.  In what follows we sometimes  identify an operator
with the associated (infinite dimensional) matrix in the Fourier basis. 
The  operator  ${\cal L}(\om,\e, {\bf u}) $ 
is represented by the infinite dimensional Hermitian matrix 
\be\label{Aomega}
A(\om) := A (\om,\e,  {\bf u}) := D_\om  +  T \,  ,    
\ee
 where
\be\label{diagD}
D_\om := {\rm diag}_{i \in \Z^b} \left(\begin{array}{cc} 
- \om \cdot l + \|j\|^2 +m & 0\\ 0 &
\om \cdot l + \| j  \|^2 + m \end{array}\right) \, , 
\ee
\be\label{def:euc}
i := (l,j) \in \Z^{b}  := \Z^\nu \times \Z^d  \,, \qquad 
\| j \|^2 := j_1^2+ \ldots + j_d^2 \, , 
\ee
and
\be\label{Tmatrix}
T := (T_i^{i'})_{i \in \Z^b , i'\in \Z^b} \, , \ \
T_i^{i'} := - \e (T_1)_i^{i'} + (T_2)_i^{i'}
\, , 
\ee
\be\label{Tmatrix12}
(T_1)_i^{i'} = \left(
\begin{array}{cc}
p_{i-i'}& q_{i-i'}\\
(\overline{q})_{i-i'}& p_{i-i'} 
\end{array} \right) \, , \quad 
(T_2)_i^{i'} = 
\left(
\begin{array}{cc}
(V_0)_{j-j'}& 0 \\
0& (V_0)_{j-j'} 
\end{array} \right) \, ,
\ee
where $ p_i, q_i, (V_0)_j $ denote the Fourier coefficients of  $ p(\vphi, x), q(\vphi, x), V_0(x) $. 

Note that $ (T_i^{i'})^\dag = T_{i'}^{i} $ (the symbol $ \dag $ denotes the conjugate transpose ) 
because
$ (\overline{q})_{i-i'} = \overline{{q}_{i'-i}}$ and  $  \overline{p}_i = p_{-i} $, since $ p $ is  real-valued.
The matrix $ T $ is  {\it T\"oplitz}, namely $ T_i^{i'} $ depends only on the difference of the indices $ i - i' $.
Moreover, since the functions $ p, q  $ in (\ref{pq}), as well as the potential $ V $, 
are in $ H^s $,   then $ T_i^{i'} \to 0  $ as $ |i- i'| \to \infty $  at a polynomial rate.  
In the next section we introduce precise norms to measure such off-diagonal decay.

Moreover we shall introduce a further index $ a \in \{0,1 \} $ to distinguish the two
eigenvalues $ \pm \om \cdot l + \|j\|^2 +m  $ (see \eqref{diao}) and the four elements of each 
of these $ 2 \times 2 $ matrices, see Definition \ref{matrixrepre} and \eqref{Tiaa'}.

\smallskip

We introduce the one-parameter family of infinite dimensional matrices 
\be\label{Atranslated}
A(\o,\teta) := A(\om) + \theta Y :=   D_\om +  T +  \theta \, Y
\ee
where 
\be\label{def:Y}
Y := {\rm diag}_{i \in \Z^b}\left(\begin{array}{cc} 
- 1 & 0 \\ 
0 & 1 \end{array}\right) \, . 
\ee
The reason for adding  $ \teta Y $ is that,  translating  the time Fourier indices
$$ 
(l,j) \mapsto (l+ l_0, j)
$$
in $ A(\om) $, gives $ A (\o, \teta)$ with $ \teta = \om \cdot l_0 $, see 
(\ref{shifted}): the matrix $ T $  remains unchanged under translation because it is T\"oplitz.

\begin{remark}\label{adv2}
The covariance property \eqref{shifted} will be exploited in 
section \ref{sec:sepa} to prove ``separation properties" of the ``singular sites". 
\end{remark}

We shall study properties of the linearized systems $A(\om,\e, {\bf u}) v=h$ in
sections $3-6$. To apply the results of these sections to the Nash-Moser
scheme of section $7$, we have to keep in mind that $ {\bf u} $ itself depends 
on the parameters $(\om,\e)$ (in a $C^1$ way, with some bound on
$ \| {\bf u} \|_{s_1} +  \| \partial_{(\om,\e)} {\bf u} \|_{s_1} $). Therefore the frame of 
sections $3-6$ will be the following: we study parametrized families of
(infinite dimensional) matrices 
\be \label{matrpar}
A(\e,\l,\theta)=D(\l) +  T(\e,\l) + \theta Y,
\ee 
where $D(\l)$ is defined by (\ref{diagD}) with $\om=\l \bar{\om}$, and  
$T$ is a T\"oplitz matrix such that 
$\norsone{T} + \norsone{\partial_{(\lambda, \e)} T} \leq C $ ($C$ depending on $ V $).

\smallskip

The main goal of the following sections is to prove polynomial off-diagonal decay for the inverse
of  the $ 2(2N+1)^b $-dimensional sub-matrices of $ A(\e, \l, \theta) $ centered at $ (l_0, j_0) $ denoted by
\be\label{ANl0}
A_{N,l_0, j_0}(\e, \l, \theta) := A_{|l - l_0| \leq N, |j - j_0| \leq N}(\e, \l, \theta) 
\ee  
where 
\be\label{supeuc}
|l| :=  \max \{|l_1|, \ldots, |l_\nu|\} \, ,  \ \  |j| := \max \{|j_1|, \ldots, |j_d|\} \, 
, \  \ |j| \leq \|j\| \leq \sqrt{d} |j| \ , .
\ee 
If $ l_0 = 0 $ we use the simpler notation
\be\label{ANj0}
A_{N, j_0}(\e, \l, \theta) :=  A_{N,0,j_0} (\e, \l, \theta) \, .
\ee
If also $ j_0 = 0 $,  we simply write
$$
A_{N}(\e, \l, \theta) :=  A_{N,0} (\e, \l, \theta) \, ,
$$
and, for $ \teta = 0 $, we denote
$$
A_{N,j_0}(\e, \l) :=  A_{N,j_0} (\e, \l, 0) \, .  
$$
We have the following crucial {\it covariance} property 
\be\label{shifted}
A_{N, l_1,  j_1 } (\e, \l, \theta) = A_{N, j_1} (\e, \l , \teta +  \l  \bar{\om} \cdot l_1 ) \, , \ee
which will be exploited in Lemma \ref{Ntime}. 

A major role is played by
the eigenvalues of $ D(\lambda) + \teta Y $,  
$$
d_i^{\pm}  := d_i^{\pm} (\l , \theta) := \pm \l \bar{\om} \cdot l + \|j\|^2  + m \pm \teta  \, . 
$$
In order to distinguish between the $ \pm $ sites we introduce an index
$$
a \in \{0,1\} 
$$ 
and we  denote
\be\label{diao}
d_{i,a}(\l , \theta) = 
\begin{cases}
\ \ \l \bar{\om} \cdot l + \|j\|^2  + m + \teta \qquad {\rm if} \ a = 0 \cr
  - \l \bar{\om} \cdot l + \|j\|^2  + m - \teta  \qquad {\rm if} \ a = 1 \, . 
\end{cases}
\ee

\section{Matrices with off-diagonal decay}\setcounter{equation}{0}\label{sec:off}

Let us consider the  basis of the vector-space $ {\bf H}^s := H^s \times H^s  $ 
made up by  
\be\label{basisi}
e_{i,0}  := (e^{\ii (l \cdot \vphi + j \cdot x)},0),  \ 
e_{i,1} := (0,e^{\ii (l \cdot \vphi + j \cdot x)}),  \ \ i := (l,j) \in \Z^b := \Z^\nu \times \Z^d \, .
\ee
Then we write any $ {\bf u} = (u^+,u^-) \in H^s \times H^s $ as
$$
{\bf u} = \sum_{k \in \Z^b \times \{0,1\}}  u_{k} e_k  \, , \quad k := (i,a) \in \Z^b \times \{0,1\} \, , 
$$
where $ {\bf u}_{l,j,0} := u_{l,j}^+ $, resp. $ {\bf u}_{l,j,1} := u_{l,j}^- $, 
denote the Fourier indices of $ u^+ $, resp. $ u^- $, see \eqref{def:Hs}. 

For $ B \subset \Z^b \times \{ 0,1\} $, we introduce the subspace
$$
{\bf H}^s_B := \Big\{ {\bf u} \in H^s \times H^s \, : \, u_{k} = 0 \ {\rm if} 
\ k \notin B \Big\} \, .
$$
When $ B $ is finite, the space $ {\bf H}^s_B $ does not 
depend on $ s $ and will be denoted $ \Hb_B $. We define 
$$ 
\Pi_B : \Hb^s \to \Hb_B
$$ 
the $L^2$-orthogonal projector onto $ \Hb_B $. 

In what follows $ B, C, D, E $ are finite  subsets of $ \Z^b\times \{ 0,1\} $.

We identify the space
$ \lin^B_C $ of the linear maps $L : \Hb_B \to \Hb_C $  with the space of matrices
$$
\matr^B_C := \Big\{ M = (M^{k'}_k)_{k' \in B, k \in C} \, , \ M^{k'}_{k} \in \C  \Big\} 
$$ 
according to the following usual definition. 

\begin{definition}\label{matrixrepre}
The matrix $ M \in \matr^B_C  $ represents the linear operator $ L \in  \lin^B_C $, if 
$$ 
\forall k'= (i',a') \in B, \, k=(i,a) \in C \, ,  \quad
\Pi_k L e_{k'}=M_k^{k'} e_k  \, , 
$$
where $ e_{i,0}$, $ e_{i,1} $ are defined in \eqref{basisi} and $ M_k^{k'} \in \C $.
\end{definition}
For example, with the above notation, the matrix elements of the matrix $ (T_1)_i^{i'} $ in \eqref{Tmatrix12} 
are
\be\label{Tiaa'}
(T_1)_{i,0}^{i',0} = p_{i-i'} \, , \ (T_1)_{i,0}^{i',1} = q_{i-i'} \, , \ (T_1)_{i,1}^{i',0} = (\overline{q})_{i-i'} 
= \overline{q_{i'-i}} \, , \ (T_1)_{i,1}^{i',1}  = p_{i-i'} \, .
\ee
{\sc Notations.} 
For any subset $B$ of $\Z^b \times \{0,1\}$, we denote by 
\be\label{proj}
\ov{B} := {\rm proj}_{\Z^b} B 
\ee
the projection of $B$ in $\Z^b$.

Given $ B \subset B' $, $ C \subset C' $ $ \subset \Z^b \times \{0,1\} $ and $ M \in {\cal M}_{C'}^{B'} $ we can introduce
the restricted matrices  
\be\label{defMBC}
M^B_C := \Pi_C M_{|\Hb_B} \, ,  \quad
M_C := \Pi_C M \, , \quad M^B := M_{|\Hb_B} \, .
\ee
If $ D \subset {\rm proj}_{\Z^b} B'$, $E \subset {\rm proj}_{\Z^b} C'  $, then we define 
\be\label{BCtilde}
M_E^{D} \ \ {\rm  as } \ \ M_{C}^{B}
\quad {\rm where} \quad   {B} := (D \times \{0,1\}) \cap B', \
 {C} := (E\times \{0,1\}) \cap C' \, .
\ee
In the particular case  $ D = \{ i' \} $,  $ E := \{ i \} $, $i,i' \in \Z^b$,  we use the simpler notations  
\be\label{line1}
M_i := M_{\{ i \}} \quad {\rm (it \ is \ either \ a \ line \ or \ a \ group \ of \ two \ lines \ of} \   M {\rm )},
\ee
\be\label{line2}
M^{i'} := M^{\{ i' \}} \quad {\rm (it \ is \ either \ a \ column \ or \ a \ group \ of \ two \ columns \ of} \ M {\rm )}, 
\ee
and 
\be\label{line3} 
M_i^{i'} := M_{\{ i \}}^{\{ i' \}} \, , 
\ee
it is a $ m \times m' $-complex matrix, where $ m \in \{1,2 \} $ (resp. $ m' \in \{1,2 \} $) is
the cardinality of $ C $ (resp. of $  B $) defined in \eqref{BCtilde} with    $ E := \{ i \} $ (resp. $ D = \{ i' \} $). 

We endow the vector-space of the $ 2\times 2$ (resp. $ 2 \times 1 $, $1\times 2$, $1\times 1$) 
complex matrices with a norm $ | \ | $ such that 
$$
|U W | \leq |U\| W| \, , 
$$
whenever the dimensions of the matrices make their multiplication possible, and $|U| \leq |V|$ if $U$ is a submatrix of $V$. 

\begin{remark}
The notations in \eqref{BCtilde}, \eqref{line1}, \eqref{line2}, \eqref{line3}, may be not very specific, but it is deliberate:  
it is convenient not to distinguish the index $a \in \{0,1\} $, which is irrelevant
in the definition of the $ s $-norms, in Definition \ref{defnormatr}.
\end{remark}

We also set the $ L^2 $-operatorial norm 
\be\label{L2norm}
\| M^B_C \|_0 := \sup_{h \in \Hb_B, h \neq 0} \frac{\| M^B_C h \|_0}{\| h \|_0}  
\ee
where  $ \| \ \|_0 := \| \ \|_{L^2 } $.

\begin{definition} \label{defnormatr} {\bf ($s$-norm)}
The $ s $-norm of a matrix $ M \in \matr^B_C $ is defined by
\be\label{defM}
\norma M \norma_s^2 := K_0 \sum_{n \in \Z^b} [M(n)]^2 \langle n \rangle^{2s}  
\ee
where $ \langle n \rangle := \max (|n|,1)$ , 
\be\label{Mpri}
[M(n)] := \begin{cases}
\max_{i-i'=n, i \in {\ov C}, i' \in {\ov B}} |M^{i'}_i|  \ \ \, \, \quad   {\rm if}  \ \   n \in \ov{C}-\ov{B}\\
0 \qquad \qquad \qquad \qquad \qquad  {\rm if}  \ \  n \notin \ov{C}-\ov{B} 
\end{cases}
\ee
with $ \ov{B} := {\rm proj}_{\Z^b} B $, $ \ov{C} := {\rm proj}_{\Z^b} C $
(see \eqref{proj}), and the  constant $ K_0 > 0 $ is introduced in  \eqref{def:Hs}.
\end{definition}

It is easy to check that $\nors{\ }$ is a norm on $ \matr_C^B $. It verifies $ \nors{\ } \leq \norma \ \norma_{s'} $, 
$ \forall s \leq s' $, and
$$ 
\forall M \in \matr^B_C \, , 
\quad \forall B' \subseteq B \ , \ C' \subseteq C \ , \quad \nors{M^{B'}_{C'}} \leq \nors{M} \, . 
$$
The $ s $-norm is designed to estimate the off-diagonal decay of matrices like  $ T $ in (\ref{Tmatrix})
with $ p, q, V \in H^s $.

\begin{lemma}\label{lem:multi}
The matrices $ T_1 $, $ T_2 $ in \eqref{T1}, \eqref{Lomega} 
with $ p, q, V \in H^s $, satisfy
\be\label{multiplication}
\nors{T_1} \leq K \| (q,p) \|_s \, ,  \quad \nors{T_2} \leq K \| V \|_s \, .
\ee
\end{lemma}
\begin{pf}
By \eqref{Mpri}, \eqref{Tmatrix12} we get
$$
[T_1 (n)] := \max_{i-i'=n} \Big| \left(
\begin{array}{cc}
p_{i-i'} & q_{i-i'} \\
\ov{q_{i-i'}} & p_{i-i'}
\end{array}
\right)   \Big|  
\leq K (|p_n|+|q_n|) \, .
$$
Hence, the definition in \eqref{defM} implies
$$
\norma T_1 \norma_s^2 = K_0 \sum_{n\in \Z^b} [T_1(n)]^2 \la n \ra^{2s}
\leq K_1 \sum_{n \in \Z^b} (|p_n|+|q_n|)^2 \la n \ra^{2s}
\leq K_2 \|(p,q)\|_s^2 
$$
and \eqref{multiplication} follows. The estimate for $\nors{T_2} $ is similar. 
\end{pf}

In order to prove that the matrices with finite 
$ s $-norm satisfy the interpolation inequalities (\ref{interpm}), and then the algebra property
(\ref{algebra}), the guiding principle
is the analogy between these matrices and the operators of the form \eqref{T1}, 
i.e. the multiplication operators for functions.
We introduce the subset $ \mapp $ of $ \cap_{s \geq 0} H^s  $ 
formed by the trigonometric polynomials with positive Fourier coefficients
$$
\mapp 
:= \Big\{ h = \sum h_{l,j} e^{\ii (l \cdot \varphi + j \cdot x)} \
{\rm with} \ h_{l,j} \neq 0 \ {\rm for \ a \ finite \ number \ of} \ (l,j) \ {\rm only \  and \ }  
h_{l,j} \in \R_+ \Big\}\, . 
$$
Note that the  sum and the product of two functions in $ \mapp $ remain in $ \mapp $. 

\begin{definition}\label{domina} 
Given $ M \in \matr^B_C $, $ h \in \mapp $, we say that $M$ is dominated by $ h $, and we write $ M \prec h $,
if 
\be\label{Msec}
[M(n)] \leq h_n \, , \quad \forall n \in \Z^b \, ,
\ee
in other words  if  $\  |M_i^{i'}| \leq h_{i-i'} $ , $ \forall i' \in {\rm proj}_{\Z^b} B $, $ i \in {\rm proj}_{\Z^b} C $. 
\end{definition}

It is easy to check  ($B$ and $C$ being finite) that
\be \label{normdom}
\norma M \norma_s = \min \Big\{ \|h\|_s \  : \  h\in \mapp \ , \  M \prec h \Big\} \quad {\rm and} \quad 
\exists h \in \mapp \ ,  \  \forall s\geq 0 \, ,\  \nors{M} =  \|h\|_s \,   . 
\ee

\begin{lemma} \label{domprod}
For $ M_1 \in \matr^C_D $, $M_2 \in \matr^B_C$, $M_3 \in \matr^C_D$, we have 
$$
M_1\prec h_1 \, , \
M_2\prec h_2 \, , \ M_3\prec h_3 \quad \Longrightarrow \quad  M_1+M_3 \prec h_1 + h_3  \quad {\rm and} \quad M_1M_2 \prec h_1 h_2 \, .
$$
\end{lemma}

\begin{pf}
Property $ M_1+M_3 \prec h_1 + h_3 $ is straightforward. 
For $i\in {\rm proj}_{\Z^b}  D$, $i'\in {\rm proj}_{\Z^b} {B}$, we have
\begin{eqnarray*}
| (M_1M_2)^{i'}_i |  = \Big| \sum_{q \in \ov{C} := {\rm proj}_{\Z^b} C} (M_1)^q_i (M_2)_q^{i'} \Big| &\leq &
\sum_{q \in \ov{C}} |(M_1)^q_i|  |(M_2)_q^{i'}| \leq  \sum_{q \in \ov{C}} (h_1)_{i-q} (h_2)_{q-i'} \\
&\leq & \sum_{q \in \Z^b} (h_1)_{i-q} (h_2)_{q-i'}=(h_1h_2)_{i-i'}
\end{eqnarray*}
implying $ M_1M_2 \prec h_1 h_2  $ by Definition \ref{domina}. 
\end{pf}

We immediately deduce from   (\ref{lions}) and (\ref{normdom})
the following interpolation estimates. 

\begin{lemma} \label{prodest}
{\bf (Interpolation)} 
$ \forall s \geq s_0 > (d+\nu)/2 $ there is $ C(s) \geq 1 $, with $C(s_0)=1$,  such that, for any finite subset $ B, C, D \subset \Z^b \times \{ 0, 1 \} $,
$ \forall M_1 \in \matr^C_D$, $ M_2 \in \matr^B_C $, 
\be \label{interpm}
\nors{M_1 M_2} \leq  (1/2) \norso{M_1} \nors{M_2} + (C(s)/2) \nors{M_1} \norso{M_2} \, ,   
\ee
in particular, 
\be \label{algebra}
\norma M_1 M_2 \norma_s \leq  C(s) \nors{M_1} \nors{M_2} \, .
\ee
\end{lemma}

Note that the constant $ C(s) $ in Lemma \ref{prodest} is independent of $ B $, $ C $, $ D $. By (\ref{algebra}) with $ s = s_0 $,  
we get (recall that $ C(s_0 ) = 1 $)

\begin{lemma} \label{prodest1}
For any finite subset $ B, C, D \subset \Z^b \times \{ 0, 1 \}  $, for all $  M_1 \in \matr^C_D $, $ M_2 \in \matr^B_C $, we have
\be \label{norsoest}
\norso{M_1M_2} \leq \norso{M_1}\norso{M_2} \, , 
\ee
and, $ \forall M \in \matr^B_B $, $ \forall n \geq 1 $, 
\label{powerest} 
\be\label{Mnab}
\norma M^n \norma_{s_0} \leq \norma M \norma_{s_0}^n \qquad
{\rm and}  \qquad   
\norma M^n \norma_{s} \leq  C(s) \norma M \norma_{s_0}^{n-1}  \norma M \norma_{s} \, , \ \forall s \geq s_0 \, .
\ee
\end{lemma}

\begin{pf}
The second estimate in (\ref{Mnab}) is obtained from (\ref{interpm}),  using $C(s) \geq 1$. 
\end{pf}

The $ s $-norm of a matrix  $ M \in \matr^B_C $ controls also the Sobolev $ H^s $-norm. Indeed, 
we  identify $ \Hb_B $ with the space $ \matr_B^{\{0\}} $ of column matrices and 
the Sobolev norm $ \|  \ \|_s $ is equal to the $ s $-norm $ \norma \ \norma_s $, i.e.
\be\label{vettori=}
\forall w \in \Hb_B \, , \ \ \| w \|_s  =  \nors{w} \, , \quad  \forall s \geq 0 \, .
\ee
Then $ M w \in \Hb_C $ and the next lemma is a particular case of Lemma \ref{prodest}.

\begin{lemma}\label{sobonorm} {\bf (Sobolev norm)}
$ \forall s \geq s_0 $ there is $ C(s) \geq 1 $ such that, for any finite subset $ B, C \subset \Z^b \times \{0, 1\}$, 
\be\label{opernorm}
\| Mw \|_s \leq (1/2) \norma M\norma_{s_0} \|w\|_s + (C(s)/2)  \norma M \norma_s \|w\|_{s_0} \, , \quad 
\forall M \in \matr^B_C \, , \ w \in {\bf H}_B \, .
\ee
\end{lemma}

\noindent
The following lemma is the analogue of the smoothing properties (\ref{S1})-(\ref{S2})
of the projection operators.

\begin{lemma} \label{norcomp}
{\bf (Smoothing)} 
Let $ M  \in \matr^B_C $ and $ N \geq 2 $. Then, $ \forall s' \geq  s \geq 0 $,
\be \label{Sm1}
M_i^{i'}=0 \, , \ \forall |i - i' | < N \quad \Longrightarrow  \quad  \nors{M} \leq N^{-(s'-s)} \norma M \norma_{s'} \, , \ee
\be \label{Sm2} 
M_i^{i'}=0 \, , \ \forall |i - i' | > N 
\quad  \Longrightarrow  \quad 
\begin{cases} 
\norma M \norma_{s'}   \leq N^{s'-s} \nors{M} \, \\
\nors{M} \leq N^{s+b}  \| M \|_{0} \, . 
\end{cases} \ee
\end{lemma} 

\begin{pf}
Estimate (\ref{Sm1}) and the first bound of (\ref{Sm2})
follow from the definition of the norms $\norma \  \norma_s $. The second bound 
of  (\ref{Sm2}) follows by the first bound in (\ref{Sm2}), noting that $ | M^{i'}_i | \leq \| M \|_0 $, $ \forall i, i' $, 
$$
\nors{M}  \leq N^{s}  \norma{M} \norma_0 \leq 
N^s \sqrt{(2N+1)^{b}} \| M \|_{0} \leq   N^{s+b} \| M \|_{0}
$$
for $ N \geq 3 $.
\end{pf}

In the next lemma 
we bound the $ s $-norm of a matrix in terms of the $ (s + b)$-norms of its lines. 

\begin{lemma} \label{norextracted}
{\bf (Decay along lines)}
Let $ M \in \matr^B_C$. Then, $ \forall s\geq 0$, 
\be\label{Mdecay} 
\nors{M} \leq K_1 \max_{i \in {\rm proj}_{\Z^b} C} \norma M_{\{i\}} \norma_{s+b} 
\ee
(we could replace the index $ b $ with any $ \a > b / 2 $).
\end{lemma}

\begin{pf}
For all $ i \in \ov{C} := {\rm proj}_{\Z^b} C $, $ i' \in \ov{B} := {\rm proj}_{\Z^b} B $, $ \forall s \geq 0 $, 
$$
|M_i^{i'}|  \leq \frac{ \norsb{M_{\{i\}}} }{\langle i - i' \rangle^{s+b}} \leq   \frac{m(s+b)}{\langle i - i' \rangle^{s+b}} 
$$
where $ m(s+b) := \max_{i \in \ov{C}} \norsb{M_{\{i\}}} $. 
As a consequence 
$$
\nors{M} = \Big( \sum_{n \in \ov{C} - \ov{B}} (M[n])^2 \langle n \rangle^{2s} \Big)^{1/2} \leq  
m(s+b) \Big( \sum_{n \in \Z^b} \langle n \rangle^{-2b}   \Big)^{1/2} = m(s+b) K(b) 
$$
implying (\ref{Mdecay}).
\end{pf}

The $ L^2 $-norm  and $ s_0 $-norm of a matrix are related.

\begin{lemma} Let $ M \in {\cal M}_B^C $. Then, for $ s_0 > (d+\nu) / 2 $,
\be\label{schur}
\| M \|_0 \leq \norso{M} \, .  
\ee
\end{lemma}

\begin{pf}
Let $ m \in {\cal H}_+ $ be  such that $ M \prec m $ and $ \nors{M} = \| m \|_s $ for 
all $ s \geq 0 $, see (\ref{normdom}). Also for $H\in {\cal M}^{\{ 0\}}_C$, there is $ h \in {\cal H}_+ $
such that $ H \prec h $ and $ \nors{H} = \|h\|_s $, $ \forall s \geq 0 $. Lemma \ref{domprod} implies that 
$ M H \prec m h $ and so
$$
\norma MH \norma_0 \leq \|mh\|_0 
\leq |m|_{L^\infty} \|h\|_0
\stackrel{\eqref{embedding}} \leq \|m\|_{s_0}  \|h\|_0 = \norso{M} \norma H \norma_0 \, , \quad
\forall H\in {\cal M}^{\{ 0\}}_C \, .
$$
Then (\ref{schur}) follows (recall \eqref{vettori=}).
\end{pf}

It will be  convenient to use the notion of left invertible operators.

\begin{definition} {\bf (Left Inverse)}
A matrix $ M \in \matr_C^B $ is left invertible if there exists  
$ N \in \matr^C_B $  such that $ NM = I_B $. Then $ N $ is called a left inverse of  $ M $. 
\end{definition}

Note that $ M $ is left invertible if and only if $ M $ (considered as a linear map) is injective
(then ${\rm dim} \, \Hb_C \geq {\rm dim} \, \Hb_B $). The left inverses  of $ M $  are not unique
if ${\rm dim} \, \Hb_C > {\rm dim} \, \Hb_B $ : 
they are uniquely defined only on the range of $ M $.

\smallskip

We shall often use the following perturbation lemma for left invertible operators. Note that 
the bound  (\ref{NR12}) for the perturbation in $ s_0 $-norm only, allows to estimate
the inverse (\ref{inv2}) also in $ s \geq s_0 $ norm. 

\begin{lemma} \label{leftinv} {\bf (Perturbation of left invertible matrices)}
If $ M \in \matr_C^B $ has a left inverse $ N \in \matr^C_B $ , 
then  
\be\label{NR12}
\forall P \in \matr_C^B \qquad {\rm with} \qquad 
\norma N \norma_{s_0}  \norma P \norma_{s_0} \leq 1/2 \, , 
\ee 
the matrix $ M + P $ has a left inverse $ N_P $ that satisfies
\be \label{inv1}
\norma N_P \norma_{s_0} \leq 2\norma N \norma_{s_0} \,,
\ee
and, $ \forall s \geq s_0 $,
\begin{eqnarray}  \label{inv12}
\norma N_P \norma_{s} & \leq  & \Big(1+ C(s)\norma N \norma_{s_0} \norma P \norma_{s_0} \Big) \norma N \norma_s + 
C(s) \norma N \norma^2_{s_0} \norma P \norma_s \\
& \leq & C(s)\Big( \norma N \norma_s + \norma N \norma^2_{s_0} \norma P \norma_s \Big)\, . \label{inv2}
\end{eqnarray}
Moreover, 
\be\label{NR0}
\forall P \in \matr_C^B \qquad {\rm with} \qquad \| N \|_0  \| P \|_0 \leq 1/2 \, , 
\ee 
the matrix $ M + P $ has a left inverse $ N_P $ that satisfies
\be \label{inv10}
\| N_P \|_0 \leq 2 \| N \|_0 \, .
\ee
\end{lemma}

\begin{pf}
We simplify notations denoting  $ C(s) $ any constant that depends  on $ s $ only. 
\\[1mm]
{\bf Step I.} {\it Proof of (\ref{inv1}).} 
\\[1mm]
The matrix $ N_P = A N $  with $ A  \in \matr^B_B $ is a left inverse of $ M + P $ if and only if
$$
I_B = A N (M+P) = A (I_B + NP) \, , 
$$
i.e. if and only if $ A $ is the inverse of $ I_B + NP \in \matr^B_B $. 
By (\ref{NR12}) $\norso{NP} \leq 1/2$, hence  the matrix $ I_B + NP $ is invertible and
\be \label{devnr}
N_P = A N = (I_B + NP)^{-1} N = \sum_{p=0}^\infty (-1)^p (NP)^p N 
\ee
is a left inverse of $ M + P $.  
Estimate (\ref{inv1})
is an immediate consequence of 
(\ref{devnr}), (\ref{norsoest}) and (\ref{NR12}). 
\\[1mm]
{\bf Step II.} {\it Proof of (\ref{inv12}).} 
\\[1mm]
For all $ s \geq s_0 $ 
\begin{eqnarray}
\forall p \geq 1, \ \nors{(NP)^p N} & \stackrel{(\ref{interpm})} \leq &  
C(s) \norso{N} \nors{(NP)^p}  + C(s)  \nors{N}  \norso{(NP)^p} \nonumber \\
& \stackrel{ (\ref{Mnab}) } \leq & {C(s)} \norso{N} \norso{NP}^{p-1} \nors{NP}  + {C(s)} \nors{N}  \norso{NP}^{p}  \nonumber \\
& \stackrel{(\ref{NR12}), (\ref{interpm})} \leq & 
C(s) 2^{-p} ( \norso{N} \norso{P}  \nors{N} +  \norso{N}^2 \nors{P}) \, . \label{stimahigh}
\end{eqnarray}
We derive (\ref{inv12}) by
$$
\nors{N_P} \stackrel{(\ref{devnr})} \leq  \nors{N} + \sum_{p=1}^\infty \nors{(NP)^pN} 
\stackrel{(\ref{stimahigh})} \leq  \nors{N} +
C(s) (  \norso{N} \norso{P}  \nors{N} +  \norso{N}^2 \nors{P}) \, .
$$
Finally  (\ref{inv10})  
follows from (\ref{NR0})
as in Step I because the operatorial $ L^2 $-norm (see (\ref{L2norm}))
satisfies the algebra property as the $s_0$-norm in (\ref{norsoest}). 
\end{pf}

\section{The multiscale analysis: estimates of  Green functions}\setcounter{equation}{0}\label{multiscale}

The main result of this section is the multiscale Proposition \ref{propinv}.
In the whole section $\d \in (0,1)$  is fixed and $\tau'>0$, $\Theta \geq 1$
are  real  parameters, on which we shall impose some condition in 
Proposition \ref{propinv}. 

Given $ \Omega, \Omega' \subset  E \subset \Z^b \times \{ 0, 1 \} $ we define 
$$
{\rm diam}(E) := \sup_{k,k' \in E} |k-k'| \, , \qquad 
 {\rm d}(\Omega, \Omega') := \inf_{k \in \Omega, k' \in \Omega'} |k-k'|  \, , 
$$
where, for $ k= (i,a) $, $ k' := (i', a')$ we set 
$$ 
|k - k' | := \max\{|i-i'|, |a-a'|\} \, .
$$ 
\begin{definition}\label{goodmatrix}
{\bf ($N$-good/bad matrix)} The matrix $ A \in {\cal M}_E^E $, with $ E \subset \Z^b \times \{0,1\}$, 
$ {\rm diam}(E) \leq 4 N $,  is $ N $-good if $ A $ is invertible and
\be\label{Ngoodmat}
\forall s \in [s_0, s_1] \   , \ \ \nors{A^{-1}} \leq N^{\tau'+\d s}.
\ee
Otherwise $ A $ is $ N $-bad.
\end{definition}

We first define the  regular and singular sites of a matrix. 

\begin{definition}\label{regulars} {\bf (Regular/Singular sites)} 
The index $ k := (i,a) \in \Z^{b} \times \{0,1\} $  is  {\sc regular} for $A $ 
if $ |A_k^k| \geq \Theta $. Otherwise $ k $ is  {\sc singular}.
\end{definition}
 
 Now we need a more precise notion adapted to the induction process. 

\begin{definition}\label{ANreg}
{\bf ($(A,N)$-good/bad site)}
For $ A \in \matr^E_E $, we say that $ k \in E \subset \Z^b \times \{ 0, 1 \}  $ is
\begin{itemize}
\item $(A,N)$-{\sc regular} if there is $ F \subset E$ such that
${\rm diam}(F) \leq 4N$,  ${\rm d}(k, E\backslash F) \geq N$ and
$A_F^F$ is $N$-good. 
\item $(A,N)$-{\sc good}  if it is regular for $A$ or $(A,N)$-regular. Otherwise we say that $ k $ is $(A,N)$-{\sc bad}.
\end{itemize}
\end{definition}

Let us consider  the new larger scale 
\be\label{newscale}
N' = N^\chi 
\ee
with $ \chi > 1 $.  

\smallskip

For a matrix $ A \in \matr_E^E  $ we define $ {\rm Diag}(A) := ( \d_{kk'} A_k^{k'})_{k, k' \in E} $.

\begin{proposition} {\bf (Multiscale step)} \label{propinv}
Assume
\be\label{dtC}
\d \in (0,1/2) \, ,  \ \tau' > 2 \tau + b + 1 \, , \ C_1 \geq 2 \, , 
\ee
and, setting $ \kappa := \tau' + b + s_0 $,
\be\label{chi1}
\chi (\t' - 2 \t - b) >  3 (\kappa + (s_0+ b) C_1 ) \, , \ 
\chi \delta > C_1 \, ,
\ee
\be\label{s1}
S \geq  s_1 > 3 \kappa + \chi (\tau + b) + C_1 s_0 \, . 
\ee
For any given $ \Upsilon > 0 $, there exist $\Theta := \Theta (\Upsilon, s_1) > 0 $ large enough 
(appearing in  Definition \ref{regulars}),  and  $  N_0 (\Upsilon, \Theta ,  S) \in \N $  
 such that:
\\[1mm]
$ \forall N \geq N_0(\Upsilon, \Theta  , S) $, 
$ \forall E \subset \Z^b \times \{0,1\} $ with 
${\rm diam}(E) \leq 4N'=4N^\chi $ (see \eqref{newscale}), if $ A \in \matr_E^E $ satisfies 
\begin{itemize}
\item 
{\bf (H1)} $\norsone{A- {\rm Diag}(A)} \leq \Upsilon $ 
\item 
{\bf (H2)} $ \| A^{-1} \|_0 \leq (N')^{\tau}$
\item
{\bf (H3)} 
There is a partition of the  $(A,N)$-bad sites $ B = \cup_{\alpha} \Omega_\alpha$ with
\be\label{sepabad}
{\rm diam}(\Omega_\alpha) \leq N^{C_1} \, , \quad {\rm d}(\Omega_\alpha , \Omega_\beta) \geq N^2 \ , \ \forall \alpha \neq \beta \, ,
\ee
\end{itemize}
then $ A $ is $ N' $-good. More precisely
\be\label{A-1alta}
\forall s \in [s_0,S]  \ , \  \  \nors{A^{-1}} \leq \frac{1}{4} ({N'})^{\tau' } \Big( ({N'})^{\d s}+ \nors{A- {\rm Diag}(A)} \Big) \,. 
\ee
\end{proposition}

The above proposition says, roughly, the following. If $ A $ has a sufficient off-diagonal decay (assumption (H1) and (\ref{s1})), 
and if the sites that can not be inserted in  good ``small''  submatrices (of size $ O(N) $) along the diagonal of $ A $  are
sufficiently separated (assumption (H3)), then the $ L^2 $-bound $ (H2) $ for $ A^{-1} $ implies that the ``large'' matrix $ A $ 
(of size $ N' = N^\chi $ with $ \chi $ as in (\ref{chi1})) is good, and 
$ A^{-1} $ satisfies also  the bounds (\ref{A-1alta}) in $ s $-norm for $ s > s_1 $.
It is remarkable that the bounds for $ s > s_1 $
follow only by informations  on the $N$-good submatrices  in $ s_1 $-norm 
(see Definition \ref{goodmatrix}) plus, of course, the $s$-decay of $ A $. 

\smallskip

According to (\ref{chi1}) the exponent $ \chi $, which measures the new scale $ N' >> N $, 
is large with respect to the size of the bad clusters $ \O_\a $, i.e. with respect to $ C_1 $. 
The intuitive meaning is that, for $ \chi $ large enough,
the ``resonance effects" due to the bad clusters are ``negligible" at the new larger scale. 

\smallskip

The constant $ \Theta \geq 1 $ which defines the regular sites (see Definition \ref{regulars}) 
must be large enough with respect to $ \Upsilon $, 
i.e. with respect to  the off diagonal part $ {\cal T} := A - {\rm Diag}(A) $, 
see (H1) and Lemma \ref{defmatrMN}.
In the application to matrices like $ A $ in \eqref{Aomega}  the constant 
$ \Upsilon $ is proportional to $  \| V  \|_{s_1} + \e \| (p,q) \|_{s_1} $.

\smallskip

The exponent $ \t \geq \t (b) $ shall be taken large
in order to verify condition (H2), imposing lower bounds on the modulus of the eigenvalues of $ A $.
Note that $ \chi $ in (\ref{chi1}) can be taken  
large independently of $ \t $, choosing, for example, $ \t' := 3 \t + 2b $ (see remark \ref{defchi}).

\smallskip

Finally, the Sobolev  
index $ s_1 $ has to be  large with respect to $ \chi $ and $ \t $, according to (\ref{s1}).
This is also natural:  if the decay is sufficiently strong, 
then the ``interaction" between different clusters of $ N$-bad sites is weak enough.

\begin{remark}
In (\ref{sepabad})  we have fixed the separation $ N^2 $ between the bad clusters 
just for definiteness:
any separation $ N^\mu $, $ \mu > 0 $, would be sufficient. Of course,  
the smaller $ \mu > 0 $ is, the larger the Sobolev exponent $ s_1 $ has to be.
See remark \ref{rem:separ} for other comments on assumption (H3).
\end{remark}


\begin{remark}\label{advanta}
 An advantage of the multiscale Proposition \ref{propinv} with respect to  analogous lemmata 
 in \cite{B5} (see for example Lemma 14.31-\cite{B5}) is to require only an $ L^2 $-bound for the inverse of $ A $,
 and not for submatrices.  For this we use the notion of left inverse matrix in the proof. 
 \end{remark}

The proof of Proposition \ref{propinv} is divided in several lemmas. In each of them we shall assume that 
the hypotheses of Proposition \ref{propinv} are satisfied.  We set 
\be\label{AD+T}
{\cal T}  := A - {\rm Diag}(A) \, , 
\qquad \norsone{{\cal T}} \stackrel{(H1)} \leq \Upsilon \, . 
\ee
Call $G$ (resp. $B$) the set of the $(A,N)$-good (resp. bad) sites.
The partition $ E = B \cup G $ induces the orthogonal decomposition $ \Hb_E = \Hb_B \oplus \Hb_G $ and we write 
$$ 
u = u_B + u_G \qquad {\rm where}  \qquad u_B := \Pi_B u \, ,  \ u_G := \Pi_G u \, .
$$
The next Lemmas \ref{defmatrMN} and \ref{defAprime} 
say that the Cramer system $ A u = h $
can be nicely reduced along the good sites $ G $, giving rise to a (non-square) system $ A' u_B = Z h $, 
with a good control of the $ s $-norms of the matrices $ A' $ and $ Z $.
Moreover $ A^{-1} $ is a left inverse of $ A' $.

\begin{lemma} \label{defmatrMN}
{\bf (Semi-reduction on the good sites)} 
Let $ \Theta^{-1} \Upsilon \leq c_0 (s_1) $ be small enough. 
There exist $ {\cal M} \in \matr^E_G $, $ {\cal N} \in \matr^B_G $ satisfying, if $ N \geq N_1( \Upsilon)$ is large enough, 
\be  \label{Nm}
\norso{{\cal M}} \leq  c N^\kappa \, , \quad \norso{{\cal N}} \leq c \, \Theta^{-1} \Upsilon \, , 
\ee
for some $ c := c(s_1) > 0 $, and, $ \forall s \geq s_0 $, 
\be\label{Nmalta}
\nors{{\cal M}} \leq C(s) N^{2\kappa } (N^{s-s_0}+ N^{-b}\norma {\cal T} \norma_{s+b}) \, , \quad   
\nors{{\cal N}}  \leq C(s) N^{\kappa } (N^{s-s_0}+ N^{-b} \norma {\cal T} \norma_{s+b}) \, , 
\ee
such that 
$$
Au = h \quad \Longrightarrow \quad u_G = {\cal N} u_B + {\cal M} h \, . 
$$
Moreover
\be \label{redex}
u_G = {\cal N} u_B + {\cal M} h \quad \Longrightarrow \quad
\forall k \ {\rm regular} \, , \ (Au)_k = h_k \, .
\ee
\end{lemma}

\begin{pf} It is based on ``resolvent identity" arguments like in \cite{B5}. The use of the  
$s$-norms  introduced in section \ref{sec:off}  makes the proof very neat. 
\\[1mm]
{\bf Step I.}  {\it There exist $ \Gamma , L  \in \matr^E_G $ satisfying  
\be\label{Gammas0} 
\norso{\Gamma} \leq  C_0 (s_1) \Theta^{-1} \Upsilon  \, , \quad \norso{L} \leq  N^{\kappa} \, , 
\ee 
and, $ \forall s \geq s_0 $, 
\be \label{Gammam}
\nors{\Gamma}  \leq C(s) N^\kappa (N^{s-s_0} + N^{-b} \norma {\cal T} \norma_{s+b}) \, , 
\quad \nors{L}  \leq  C(s) N^{\kappa +  s -s_0} \, , 
\ee
such that} 
\be\label{Au=h}
Au = h \quad  \Longrightarrow \quad u_G+ \Gamma u = Lh \, .
\ee
Fix any $ k \in G $ (see Definition \ref{ANreg}). If $ k $  is regular, let $ F := \{k \} $,  
and,  if $  k $ is not regular but $(A,N)$-regular, 
let $ F \subset E  $ such  that  $ {\rm d}(k, E\backslash F) \geq N $, $ {\rm diam}(F) \leq 4N $,
$ A_F^F  $ is  $N$-good.
We have
\be \label{eqF}
A u = h  \quad  \Longrightarrow \quad
A^F_F u_F + A^{E\backslash F}_F u_{E\backslash F} = h_F \quad  \Longrightarrow \quad
u_F + Q u_{E\backslash F} = (A_F^F)^{-1} h_F 
\ee
where 
\be\label{defB}
Q := (A_F^F)^{-1} A^{E\backslash F}_F=(A_F^F)^{-1} {\cal T}^{E\backslash F}_F \in \matr_F^{E\backslash F} \, . 
\ee
The matrix $ Q $ satisfies 
\be \label{estB}
\norsone{Q} \stackrel{(\ref{algebra})} \leq C(s_1) \norsone{(A_F^F)^{-1}} \norsone{{\cal T}} 
\stackrel{(\ref{Ngoodmat}), (\ref{AD+T})}  \leq C(s_1) N^{\tau' + \d s_1} \Upsilon
\ee
(the matrix  $ A_F^F  $ is  $N$-good).
Moreover, $ \forall s \geq s_0 $, using the interpolation Lemma \ref{prodest}, and  $ {\rm diam}(F) \leq 4N $, 
\begin{eqnarray}
\norma Q \norma_{s+b} & \stackrel{(\ref{interpm})}\leq & 
C(s) ( \norma (A_F^F)^{-1} \norma_{s+b} \norso{{\cal T}} + \norso{(A_F^F)^{-1}} \norma {\cal T} \norma_{s+b} ) \nonumber \\
&  \stackrel{(\ref{Sm2})} \leq & 
C(s) ( N^{s+b -s_0} \norma (A_F^F)^{-1} \norma_{s_0} \norso{{\cal T}} + \norso{(A_F^F)^{-1}} \norma {\cal T} \norma_{s+b} ) \nonumber \\ 
&  \stackrel{(\ref{Ngoodmat}), (\ref{AD+T})} \leq & C(s) N^{(\d -1) s_0 }( N^{s+b+ \tau' } \Upsilon + N^{\tau'+   s_0} \norma {\cal T} \norma_{s+b} ) \,. \label{Balta}
\end{eqnarray}
Applying the projector $ \Pi_{\{k\}} $ in (\ref{eqF}), we obtain
\be\label{primos}
Au=h \quad  \Longrightarrow \quad 
u_k + \sum_{k' \in E} \Gamma_k^{k'} u_{k'}= \sum_{k'\in E} L_k^{k'} h_{k'}  
\ee
that is  (\ref{Au=h})  with
\be\label{defGL}
\Gamma_k^{k'}  := 
\begin{cases} 0 \	 \  \ \quad{\rm if} \ \ k' \in F \\
Q_k^{k'}  \ \  \ \, {\rm if} \ \ k' \in E \setminus F
\end{cases}
\quad {\rm and} \qquad 
L_k^{k'} := \begin{cases}
[(A_F^F)^{-1}]_k^{k'} \  \ \ \, {\rm if} \ \ k' \in F \\
0 \quad  \qquad \quad \ \,  \ \ {\rm if} \ \ k' \in E \setminus F. 
\end{cases}  
\ee
If $ k $ is regular then $ F = \{ k \} $, and, by Definition \ref{regulars}, 
\be\label{invAi}
| A_k^k | \geq \Theta \, .
\ee
Therefore, by (\ref{defGL}) and (\ref{defB}), the $ k $-line of $ \Gamma $ satisfies 
\be\label{lin1}
\norsob{\Gamma_k} \leq \norsob{(A_k^k)^{-1} {\cal T}_k} \stackrel{ (\ref{invAi}), (\ref{AD+T})} \leq 
C(s_0) \Theta^{-1} \Upsilon \, .
\ee
If $ k $ is not regular but $(A,N)$-regular, since  ${\rm d}(k, E \backslash F) \geq N $ 
we have, by (\ref{defGL}), that $\Gamma_k^{k'}=0$ for 
$ | k - k' | \leq N $. Hence, by Lemma \ref{norcomp}, 
\begin{eqnarray}
\norsob{\Gamma_k} & \stackrel{(\ref{Sm1})} \leq & N^{-(s_1-s_0-b)} \norsone{\Gamma_k} 
\stackrel{(\ref{defGL})} \leq N^{-(s_1-s_0-b)}  \norsone{Q} 
\stackrel{(\ref{estB})} \leq C(s_1) \Upsilon N^{\tau' + s_0 + b  - (1-\d)s_1} \nonumber \\ 
& \leq &  C(s_1) \Theta^{-1} \Upsilon  \label{lin2}
\end{eqnarray}
for $ N \geq N_0 (\Theta) $ large enough. Indeed
the exponent  $ \tau' + s_0 + b  - (1-\d)s_1 <  0 $ because $ s_1 $ is large enough  
according to (\ref{s1}) and $ \d \in (0, 1/2) $ (recall $ \kappa := \tau' + s_0 + b $).
In both cases (\ref{lin1})-(\ref{lin2}) imply that each line $ \Gamma_k $ decays like
$$
\norsob{\Gamma_k}  \leq  C(s_1) \Theta^{-1} \Upsilon \, ,\quad  \forall k \in G \, .
$$ 
Hence, by Lemma \ref{norextracted}, 
$\norso{\Gamma} \leq C'(s_1) \Theta^{-1} \Upsilon $, which is 
the first inequality in (\ref{Gammas0}). 
Likewise we prove the second estimate in (\ref{Gammas0}).
Moreover, $ \forall s \geq s_0 $, still by Lemma \ref{norextracted}, 
$$
\nors{\Gamma}  \leq K \sup_{k \in G} \norma \Gamma_k \norma_{s+b}  
\stackrel{(\ref{defGL})} \leq  K \norsb{Q} \stackrel{(\ref{Balta})}\leq C(s) N^\kappa 
(N^{s-s_0} + N^{-b} \norsb{{\cal T}})  
$$
where  $ \kappa := \tau' + s_0 + b $ and for $ N \geq N_0(\Upsilon)$.

The second estimate in (\ref{Gammam}) follows by $ \norso{L} \leq  N^{\kappa} $ (see (\ref{Gammas0})) and (\ref{Sm2})
(note that  by (\ref{defGL}), since diam$F \leq 4N$, we have $ L_k^{k'} = 0 $ for all $ |k - k' | > 4N $).
\\[2mm]
{\bf Step II.}  By (\ref{Au=h}) we have
\be\label{Auin}
Au = h \quad \Longrightarrow  \quad (I_G  + \Gamma^G)u_G  = Lh - \Gamma^B u_B \, .
\ee
By (\ref{Gammas0}), if $ \Theta $ is large enough (depending on $ \Upsilon$, namely on the potential $ V_0 $), 
we have $\norso{\Gamma^G} \leq 1/2 $. Hence,
by  Lemma \ref{leftinv}, $ I_G  + \Gamma^G $ is invertible and 
\be\label{Linver0}
\norso{(I_G  + \Gamma^G)^{-1}} \stackrel{(\ref{inv1})} \leq 2 \, , 
\ee  
\be\label{Linver}
\forall s \geq s_0 \ ,\ \ \nors{(I_G  + \Gamma^G)^{-1}} 
\stackrel{(\ref{inv2})} \leq C(s) (1+ \nors{\Gamma^G})
\stackrel{(\ref{Gammam})} \leq  C(s) N^{\kappa }(N^{s-s_0}  + 
N^{-b} \norma {\cal T} \norma_{s+b}) \, .
\ee
By (\ref{Auin}), $\  Au=h \Longrightarrow u_G = {\cal M} h + {\cal N} u_B\  $, with
\be\label{poi}
\qquad {\cal M} := (I_G  + \Gamma^G)^{-1}L \quad {\rm and} \quad
{\cal N} := - (I_G  + \Gamma^G)^{-1} \Gamma^B 
\ee
and  estimates (\ref{Nm})-(\ref{Nmalta}) follow by Lemma \ref{prodest}, (\ref{Linver0})-(\ref{Linver}) and (\ref{Gammas0})-(\ref{Gammam}).

Note that  
\be\label{conclus}
u_G + \Gamma u = Lh  \quad \iff \quad u_G = {\cal M} h + {\cal N} u_B \, . 
\ee
As a consequence, if $u_G = {\cal M} h + {\cal N} u_B$ then, by
(\ref{defGL}), for $ k $ regular,
$$
u_k + (A_k^k)^{-1} \sum_{k'\neq k} A_k^{k'} u_{k'}=(A_k^k)^{-1} h_k \, ,
$$
hence $ (Au)_k = h_k $, proving \eqref{redex}.
\end{pf}

\begin{lemma} \label{defAprime}
{\bf (Reduction on the bad sites)} We have
$$
Au = h \quad \Longrightarrow  \quad A' u_B = Zh  
$$
where
\be\label{A'Z}
A' := A^B + A^G {\cal N} \ \in  {\cal M}^B_E \, ,  \qquad Z := I_E - A^G {\cal M} \  \in  {\cal M}^E_E \, ,
\ee
satisfy
\be\label{A'a}
\norso{A'} \leq c(\Theta)  \, , \qquad \nors{A'} \leq  C(s, \Theta) N^{\kappa } (N^{s-s_0}
+ N^{-b}\norma {\cal T} \norma_{s+b})\, , 
\ee
\be\label{Za}
\norso{Z} \leq c N^\kappa \, ,  \quad  \nors{Z} \leq C(s, \Theta) N^{2\kappa } (N^{s-s_0} 
+ N^{-b} \norma {\cal T} \norma_{s+b}) \, .
\ee
Moreover $(A^{-1})_B $ is a left inverse of $ A' $. 
\end{lemma}

\begin{pf}
By Lemma \ref{defmatrMN},
$$
Au = h \quad \Longrightarrow \quad
\begin{cases}
A^G u_G + A^B u_B = h \\
u_G = {\cal N} u_B + {\cal M} h 
\end{cases} \quad \Longrightarrow \quad 
(A^G {\cal N} + A^B ) u_B = h-A^G {\cal M} h \, , 
$$
{\it i.e.} $ A' u_B = Z h $. Let us prove  estimates (\ref{A'a})-(\ref{Za}) for $ A' $ and $ Z $.
\\[1mm]
{\bf Step I.} {\it $ \forall \, k $ regular we have $ A'_k = 0 $, $ Z_k = 0 $.}
\\[1mm]
By (\ref{redex}), for all $ k $ regular, 
$$ 
\forall h \, , \  \forall u_B \in \Hb_B \ , \quad 
\Big(  A^G ({\cal N}u_B + {\cal M}h)+ A^B u_B \Big)_k = h_k \  , \quad
i.e. \quad  (A'u_B)_k = (Zh)_k \,  ,
$$
which implies $ A'_k = 0 $ and $ Z_k = 0 $. 
\\[1mm]
{\bf Step II.} {\it Proof of (\ref{A'a})-(\ref{Za})}.
\\[1mm] 
Call $ R \subset E$  the regular sites in $ E $. For all $ k \in E \backslash R $, 
we have $ |A_k^k| < \Theta $ (see  Definition \ref{regulars}). Then 
(\ref{AD+T})  implies 
\be\label{passagge}
\norso{A_{E \backslash R}} \leq \Theta + \norso{{\cal T}}  \leq c(\Theta) \, , 
\quad \nors{A_{E \backslash R}} \leq \Theta + \nors{{\cal T}} \, ,
\ \forall s \geq s_0 \, .
\ee
By Step I and the definition of $ A' $ in (\ref{A'Z}) we get
$$
\nors{A'}  = \nors{A'_{E\backslash R}} \leq \nors{A^B_{E\backslash R}} + \nors{A^G_{E\backslash R} {\cal N}} \, .
$$
Therefore, Lemma \ref{prodest}, (\ref{passagge}), (\ref{Nm}), (\ref{Nmalta}), imply  
$$
\nors{A'} \leq C(s, \Theta) N^{\kappa}  (N^{s-s_0} + N^{-b} \norma {\cal T} \norma_{s+b})  
\qquad {\rm and} \qquad \norso{A'} \leq c(\Theta) \, ,
$$
proving (\ref{A'a}). The bound (\ref{Za}) follows similarly.
\\[1mm]
{\bf Step III.} {\it $ (A^{-1})_B $ is a left inverse of $ A' $}.
\\[1mm] 
By $ A^{-1} A' = $ $ A^{-1} ( A^B + A^G {\cal N}) = $ $ I^B_E + I^G_E {\cal N} $  we get
$$
(A^{-1})_B A'=(A^{-1} A')_B=I^B_B - 0= I^B_B
$$ 
proving  that $(A^{-1})_B$ is a left inverse of $ A' $. 
\end{pf}

Now $ A' \in \matr^B_E $, and the set $B$ is partitioned in clusters $\Omega_\a$ of size
$O(N^{C_1})$, far enough one from another, see (H3).  Then, up to a remainder of very small $ s_0 $-norm 
(see (\ref{noroR0})), 
$ A' $ is defined by the submatrices $(A')^{\Omega_\a}_{\Omega'_\a}$ where $\Omega'_\a$ is some neighborhood
of $\Omega_\a$ (the distance between two distinct $\Omega'_\a$ and $\Omega'_\b$ remains large). Since
$ A' $ has a left inverse with $L^2$-norm $O({N'}^\tau)$, 
so have the submatrices $(A')^{\Omega_\a}_{\Omega'_\a}$. Since these
submatrices are of size $O(N^{C_1})$, the $s$-norms of their inverse will be estimated as
$O(N^{C_1s} {N'}^\tau)= O({N'}^{\tau + \chi^{-1} C_1 s})$, see (\ref{decaW0}). 
By Lemma \ref{leftinv},  provided $ \chi $ is chosen large enough,
$ A' $ has a left inverse $ V $ with $s$-norms satisfying (\ref{LeftY}).
The details are given in the following lemma.

\begin{lemma} \label{defY} {\bf (Left inverse with decay)}
The matrix $ A' $ defined in Lemma \ref{defAprime} has a left inverse $ V $ which satisfies 
\be\label{LeftY}
\forall s \geq s_0 \ , \  \
\nors{V}\leq 
C(s) {N}^{2\chi \tau + \kappa + 2 (s_0+b) C_1}(N^{C_1 s} + \norsb{{\cal T}})  \, . 
\ee
\end{lemma}

\begin{pf}
Define $ \dom \in \matr_E^B $ by
\be\label{grossodia}
\dom_{k'}^k :=\begin{cases}
(A')_{k'}^{k} \quad \, \rm{if}  \  \ ({\it k}, {\it k}') \in \cup_\a (\Omega_\alpha \times \Omega'_\alpha) \\
0 \qquad \quad \ \hbox{if}  \  \ ({\it k}, {\it k}') \notin \cup_\a (\Omega_\alpha \times \Omega'_\alpha) 
\end{cases}
\quad {\rm where} \quad \Omega'_\alpha := \{ k \in E \  : \ {\rm d}(k,\Omega_\alpha) \leq N^2/4 \} \, .
\ee
{\bf Step I.} {\it $ \dom $ has a left inverse $ W \in \matr^E_B $ with $ \| W \|_0 \leq 2 ({N'})^{\tau} $.}
\\[1mm]
We define $ \remain := A' - \dom $. By the definition (\ref{grossodia}), 
if  $ {\rm d}(k',k) < N^2 / 4 $ then $ \remain_{k'}^k = 0 $ and so
\begin{eqnarray} 
\norso{\remain} & \stackrel{(\ref{Sm1})} \leq &  
4^{s_1} N^{-2(s_1 - b - s_0)} \norma{\remain} \norma_{s_1 - b}  
\leq 4^{s_1} N^{-2(s_1- b - s_0)} \norma{A'} \norma_{s_1 - b} \nonumber  \\
&  \stackrel{(\ref{A'a}), (\ref{AD+T})}  \leq  & C(s_1) N^{-2(s_1 - b - s_0)} N^{\kappa} ( N^{s_1-b-s_0}+N^{-b} 
\Upsilon) \leq C(s_1) N^{2 \kappa - s_1}  \label{noroR0}
\end{eqnarray}
for $ N \geq N_0 (\Upsilon) $ large enough. 
Therefore 
\begin{eqnarray}\label{rsmall}
\| \remain \|_0 \| (A^{-1})_B \|_0 \stackrel{(\ref{schur})} \leq  \norso{\remain} \| A^{-1} \|_0 
& \stackrel{(\ref{noroR0}),(H2)}  \leq & C(s_1)  N^{2\kappa - s_1} (N')^{\t} \nonumber \\
& \stackrel{(\ref{newscale})} = & C(s_1) N^{2\kappa  - s_1+ \chi \t}   \stackrel{(\ref{s1})} \leq 1/2 
\end{eqnarray}
for $ N \geq N(s_1) $.
Since $ (A^{-1})_B  \in \matr^E_B $ is a left inverse of $ A' $ (see Lemma \ref{defAprime}), Lemma 
\ref{leftinv} and (\ref{rsmall}) imply that  
$ \dom = A' - R $ has a left inverse $ W \in \matr^E_B $, and
\be\label{W00}
\| W \|_0 \stackrel{(\ref{inv10})} \leq 2  \| (A^{-1})_B \|_0 \leq  2  \| A^{-1} \|_0 \stackrel{(H2)}\leq
2 ({N'})^{\tau} \, . 
\ee
{\bf Step II.} {\it $ W_0 \in \matr^E_B $ defined by
\be\label{W0ii}
(W_0)^{k'}_k := \begin{cases}  W^{k'}_k \quad  \quad {\rm if} \ \ \, (k,k') \in \cup_\a
(\Omega_\alpha \times \Omega'_\alpha) \\ 
0 \, \qquad \quad \, \hbox{\rm if } \ \ (k,k') \not\in \cup_\a
(\Omega_\alpha \times \Omega'_\alpha) \end{cases}
\ee
is a left inverse of $ \dom $ and $ \nors{W_0} \leq C(s) {N}^{(s+b)C_1 + \chi \tau  } $, $ \forall s \geq s_0 $.}
\\[1mm]
Since $ W \dom = I_B $, we prove that $ W_0 $ is a left inverse of $ \dom $ showing that 
\be\label{vai}
(W-W_0)\dom = 0 \, .
\ee
Let us prove (\ref{vai}). 
For $ k \in B = \cup_\a \Om_\a $, there is $\a$ such that  $k \in \Omega_\alpha $,  and
\be\label{sum}
\forall k' \in B \ , \ ((W-W_0)\dom )_k^{k'}=\sum_{q \notin \Omega'_\a} (W-W_0)_k^q \dom^{k'}_q
\ee
since $ (W - W_0)_k^q = 0 $ if $ q \in \Omega'_\a $, see the Definition (\ref{W0ii}). 
\\[1mm]
{\sc Case I:} $ k' \in \Om_\a $. Then $ \dom^{k'}_q = 0 $  in (\ref{sum}) and so $((W-W_0)\dom )_k^{k'} = 0 $.
\\[1mm]
{\sc Case II:} $ k' \in \Omega_\beta $ for some $\beta \neq \a $. Then,
since $ \dom^{k'}_q = 0 $ if $ q \notin \Omega'_\beta$, we obtain by (\ref{sum}) that
$$
((W-W_0)\dom )_k^{k'}=\sum_{q \in \Omega'_\beta} (W-W_0)_k^q \dom^{k'}_q
\stackrel{(\ref{W0ii})} = \sum_{q \in \Omega'_\beta} W_k^q \dom^{k'}_q
 \stackrel{(\ref{grossodia})}  =
\sum_{k \in E} W_k^q \dom^{k'}_q
= (W \dom)_k^{k'} = (I_B)_k^{k'} = 0 \, .
$$
Since $ {\rm diam}(\Omega'_\a) \leq 2N^{C_1} $,  definition (\ref{W0ii}) implies 
$ (W_0)_k^{k'} = 0 $ for all $ | k - k' | \geq 2N^{C_1} $.  Hence,  $ \forall s \geq 0 $,
\be\label{decaW0}
\nors{W_0} \stackrel{(\ref{Sm2})} \leq C(s) N^{(s + b) C_1 } \| W_0 \|_0
\stackrel{(\ref{W00})} \leq C(s) {N}^{(s + b)C_1 +\chi \t } \, . 
\ee
{\bf Step III.} {\it $ A' $ has a left inverse $ V $ satisfying (\ref{LeftY}).}
\\[1mm]
Now $ A' = \dom + \remain $, $ W_0 $ is a left inverse of $ \dom $, and 
$$
\norso{W_0} \norso{\remain} 
\stackrel{(\ref{decaW0}), (\ref{noroR0})} 
\leq C(s_1) N^{(s_0 + b) C_1 + \chi \tau + 2 \kappa -  s_1} \stackrel{(\ref{s1})} \leq 1/2
$$
(we use  also that $ \chi > C_1 $ by (\ref{chi1}))
for $ N \geq N(s_1) $  large enough. Hence, by Lemma  \ref{leftinv}, $A'$ has a left inverse $ V $ with
\be\label{Vs0}
\norso{V} \stackrel{(\ref{inv1})} \leq 2 \norso{W_0} \stackrel{(\ref{decaW0})} \leq C N^{ (s_0 + b)C_1 + \chi \tau }
\ee
and, $ \forall s \geq s_0 $, 
\begin{eqnarray*}
\nors{V} & \stackrel{(\ref{inv2})} \leq & C(s) ( \nors{W_0} + \norso{W_0}^2  \nors{\remain} ) 
\leq C(s) ( \nors{W_0} + \norso{W_0}^2  \nors{A'} ) \\
& \stackrel{(\ref{decaW0}), (\ref{A'a})} \leq & C(s) {N}^{2\chi \tau + \kappa + 2 (s_0+b) C_1}(N^{C_1 s} + \norsb{{\cal T}})
\end{eqnarray*}
proving (\ref{LeftY}).
\end{pf}

\noindent
{\sc Proof of Proposition \ref{propinv} completed.} 
Lemmata \ref{defmatrMN}, \ref{defAprime}, \ref{defY} imply
$$
Au = h  \quad \Longrightarrow \quad 
\begin{cases}
u_G = {\cal M} h + {\cal N}u_B \\ 
u_B = V Z h 
\end{cases}  
$$
whence 
\be\label{eccoci}
(A^{-1})_B= V Z  \qquad {\rm and} \qquad (A^{-1})_G= {\cal M} + {\cal N} V Z={\cal M} + {\cal N}(A^{-1})_B \, . 
\ee
Therefore, $ \forall s \geq s_0 $, 
\begin{eqnarray}
\nors{(A^{-1})_B} & \stackrel{(\ref{eccoci}), (\ref{interpm})} \leq & C(s)( \nors{V} \norso{Z} + \norso{V} \nors{Z})  \nonumber \\
& \stackrel{(\ref{LeftY}), (\ref{Za}), (\ref{AD+T}), \eqref{Vs0}} \leq & C(s) N^{2\kappa + 2 \chi \tau + 2 (s_0+b)C_1 } 
(N^{C_1 s} + \norsb{{\cal T}}) \nonumber  \\
& \leq & C(s) {(N')}^{\a_1} ( {(N')}^{\a_2 s} + \nors{{\cal T}}) \nonumber 
\end{eqnarray}
using $ \norsb{{\cal T}} \leq C(s) (N')^b \nors{{\cal T}} $ (by (\ref{Sm2})) and defining
$$ 
\a_1 := \, 2 \tau + b + 2 \chi^{-1}( \kappa +  C_1 (s_0 + b))\, , \quad  \a_2 := \chi^{-1} C_1 \, .  
$$
We obtain the same bound for $ \nors{(A^{-1})_G} $. Hence, for $s\in [s_0,S]$, 
\begin{eqnarray*}
\nors{A^{-1} } & \leq & \nors{(A^{-1})_B} + \nors{(A^{-1})_G} \leq C(s) {(N')}^{\a_1} ({(N')}^{\a_2 s}+\nors{{\cal T}}) \\
& \stackrel{(\ref{chi1})} \leq & \frac14 {(N')}^{\t'} ({(N')}^{\d s}+\nors{{\cal T}}) 
\end{eqnarray*}
for $ N \geq N(S) $ large enough, proving (\ref{A-1alta}). 

\section{Separation properties of the bad sites}\label{sec:sepa}\setcounter{equation}{0}

The aim of this section is to verify the separation properties of the bad sites 
required in the multiscale Proposition \ref{propinv}. 

Let  $ A := A(\e,\l,\teta ) $ be the infinite dimensional matrix defined in (\ref{matrpar}). Given $ N \in \N $
and $ i = (l_0,j_0 ) $, recall that the submatrix  $ A_{N,i} $ is defined in (\ref{ANl0}). 

\begin{definition} \label{GBsite} {\bf ($N$-good/bad site)}
A site $ k := (i,a) \in \Z^b \times \{0,1\} $ is: 
\begin{itemize}
\item $ N$-{\sc regular } if $ A_{N,i} $ is $ N $-good (Definition \ref{goodmatrix}). Otherwise we say that $ k $ is $ N$-{\sc singular}.
\item $N$-{\sc good} if  
\be\label{regNreg} 
k  \ {is \ regular} \ ({Definition} \, \ref{regulars})
\quad  {\rm or} \quad  
{ \it all \  the \ sites} \  k' \  { with} \  {\rm d}(k',k) \leq N \ { are} \ N-{ regular} \, .
\ee
Otherwise, we say  that  $ k $ is $ N $-{\sc bad}. 
\end{itemize}
\end{definition}

\begin{remark}\label{good2}
It is easy to see that a site $ k $ which is $ N $-good  according to Definition \ref{GBsite}, is $(A_E^E,N)$-good according to Definition \ref{ANreg}, 
for any set $ E = E_0 \times \{0,1\} $ containing $ k $ where 
$ E_0  \subset \Z^b $ is a product of intervals of length $\geq N$.
We introduce these different definitions 
for merely technical reasons:
it is more convenient to prove separation properties of $ N $-bad sites for  infinite dimensional matrices.
On the other hand, for a finite matrix $ A_E^E $, 
we need the notion of $(A_E^E,N)$-good sites in order to perform the ``resolvent identity" 
also near the boundary $ \partial E $, see Step I of Lemma \ref{defmatrMN}. 
\end{remark}

We define
\be\label{tetabad}
B_{N}(j_0; \e, \l )  := \Big\{ \teta \in \R \, : \,   A_{N, j_0}(\e, \l,\theta) \ {\rm is} \ N-bad  \Big\} \, . 
\ee

\begin{definition} {\bf ($ N$-good/bad parameters)} \label{def:freqgood}
A couple $ (\e, \l) \in \R^2 $ is $ N$-good for $ A $ if
\be\label{BNcomponents}
\forall \, j_0 \in \Z^d  \, , \quad 
B_{N}(j_0;\e, \l) \subset \bigcup_{q = 1, \ldots, N^{2d+\nu+4}} I_q 
\ee
where $ I_q $ are  intervals with measure $ | I_q| \leq N^{-\t} $. 
Otherwise, we say $ (\e, \l)  $ is $ N$-bad.
We define 
\be\label{good}
{\cal G}_{N} := {\cal G}_{N}({\bf u}) := \Big\{ (\e ,\l) 
\in [ 0 ,\e_0] \times  \Lambda \, : \,  (\e, \l)  \ \  {\rm is \ } \ N-{\rm good \ for \ } A  \Big\} \, . 
\ee
\end{definition}

The main result of this section is the following proposition.  It will
enable to verify the assumption (H3) of Proposition \ref{propinv} 
for the submatrices $ A_{N',j_0}(\e, \l,\theta) $, 
see Lemmata \ref{H1H3} and  \ref{S3n+1}.

\begin{proposition}\label{prop:separation}
{\bf (Separation properties of $ N $-bad sites)}
There exist $ C_1 := C_1(d, \nu) \geq 2 $  and  $ N_1 := N_1(\nu,d, \g_0, \t_0, m,  \Theta)$   such that if $N \geq N_1$ and 
\begin{itemize}
\item {\bf (i)} $ (\e, \l) $ is $ N$-good for $ A $
\item {\bf (ii)} 
$ \tau > \chi \t_0 $ ($ \, \t_0 $ is the diophantine exponent of $ \bar \om $ in (\ref{diophan0})),
\end{itemize}
then
$\  \forall \teta \in \R $,   
the $ N $-bad sites $k:= (l,j,a) \in \Z^\nu \times \Z^d \times \{0,1\} $  of $ A(\e, \l, \teta) $ with $|l| \leq N' $
admit a partition $ \cup_\a \Omega_\a $ in disjoint clusters satisfying
\be\label{separ}
{\rm diam}(\Omega_\a) 
\leq N^{C_1(d,\nu)} \, ,   \quad {\rm d}(\Omega_\a, \Omega_\b) > N^2 \, , \ 
\forall \a \neq	 \b \, .
\ee
\end{proposition}

We underline that  the estimates (\ref{separ}) are {\it uniform} in $ \teta $.

\begin{remark}\label{rem:separ}
The $ N $-bad sites  
appear  necessarily 
in clusters with increasing size $ O(N^{C_1})$,
due to the multiplicity of the eigenvalues of the Laplacian; this happens already for the singular sites  
of periodic solutions, i.e. for $ \nu = 1 $, see \cite{BB}.
It is also natural that the separation between  clusters of $ N $-bad sites increases with $ N $, because,
roughly speaking,  the $ N $-bad sites correspond small divisors of size $ O(N^{-\a}) $.
\end{remark}

\begin{remark}\label{regus}
The geometric structure of the bad and singular sites, determines 
the regularity of the solutions of Theorem \ref{thm:main}. Actually,  
the  solutions of Theorem \ref{thm:main}  have the {\it same} Sobolev regularity in time and space
because the $ N$-bad clusters are separated in the space-time Fourier indices, see (\ref{separ}).
\end{remark}

We first estimate the time Fourier components of the $ N $-singular sites. 
We use that, by (\ref{diophan0}), the frequency vectors $ \om = \l \bar \om $, $ \forall \l \in [1/2, 3/2] $, are diophantine, namely
\be\label{diophan}
|\om \cdot l | \geq \frac{\g_0}{|l|^{\tau_0}} \, , \quad \forall l \in \Z^{\nu} \setminus \{ 0 \} \, ,
\ee 
and we use  the  ``complexity" information (\ref{BNcomponents}) on the set
$ B_{N}(j_0;\e, \l) $. This kind of argument was used in 
\cite{E1} and \cite{B5}. 

\begin{lemma}\label{Ntime} Assume (i)-(ii) of Proposition \ref{prop:separation}. Then, 
$ \forall j_1 \in \Z^d $, the number of $ N $-singular sites $ (l_1,j_1,a_1) \in \Z^\nu \times \Z^d \times \{0,1\} $ 
with $ |l_1| \leq N'  $  does not exceed $ 2 N^{2d+\nu+4} $.
\end{lemma}

\begin{pf}
If $(l_1,j_1,a_1) $ is $ N$-singular then 
$ A_{N,l_1,j_1}(\e, \l, \teta) $ is $ N$-bad (see Definitions \ref{GBsite} and \ref{goodmatrix}).
By (\ref{shifted}), we get that $ A_{N,j_1}(\e, \l, \teta + \l \bar{\om} \cdot l_1) $ is $ N$-bad, namely
$ \teta + \l \bar{\om} \cdot l_1 \in B_{N}(j_1;\e, \l ) $ (see (\ref{tetabad})).
By assumption, $ (\e, \l)$ is $ N$-good, and, therefore, (\ref{BNcomponents}) holds.

We claim that in each interval $ I_q $ there is at most one element $ \theta + \om \cdot l_1  $ with 
$\o = \l \bar{\o}$, $ |l_1| \leq N' $. 
Then, since there are at most $ N^{2d+\nu+4} $  intervals $ I_q $ (see (\ref{BNcomponents}))
and $ a \in \{0,1\} $, the lemma follows. 

We prove the previous claim by contradiction. 
Suppose that there exist $ l_1 \neq l_1' $ with $ |l_1|, |l_1'| \leq N'$, 
such that $ \om \cdot l_1 + \theta $, $ \om \cdot l_1' + \theta \in I_q $.  Then 
\be\label{up}
|\om \cdot (l_1 - l_1')| = |(\om \cdot l_1 + \theta) - (\om \cdot l_1' + \theta) | \leq |I_q| \leq N^{- \tau} \, .
\ee
By (\ref{diophan}) we also have
\be\label{lo}
|\om \cdot (l_1  - l_1')| \geq \frac{\g_0}{|l_1-l_1'|^{\tau_0}}  \geq \frac{\g_0}{ (2 N' )^{\tau_0}} = 2^{-\t_0} \g_0 N^{- \chi \tau_0} \, . 
\ee
By assumption (ii) of Proposition \ref{prop:separation} 
the inequalities (\ref{up}) and (\ref{lo}) are in contradiction, for $ N \geq N_0(\g_0, \t_0) $ large enough. 
\end{pf}

We now estimate  also the spatial components of the sites 
\be\label{sn}
{\cal S}_N := \Big\{ k = (l,j,a) \in \Z^{\nu + d}  \times \{0,1\} \,  :   \,  |l| \leq N', \, k  \ {\rm is} \, 
{\rm \, singular \,  and}  \, N-{\rm singular} \, {\rm for } \, A(\e, \l, \teta) 
\Big\} \, . 
\ee
In order to achieve a partition in clusters of $ {\cal S}_N $ we use the  notion of ``chain" of singular sites already used
for the search of periodic solutions of NLS and NLW in higher dimension in \cite{B4}, \cite{BB}. 

\begin{definition} {\bf ($ M $-chain)} 
A sequence $ k_0, \ldots , k_L 
\in \Z^{d+\nu} \times \{0,1\}  $ of distinct integer vectors
satisfying, for some $ M  \geq 2 $, 
$  | k_{q+1} - k_q | \leq M $, $ \forall q = 0, \ldots, L - 1 $, is called a $M$-chain of length $ L $.
\end{definition}

Proposition \ref{prop:separation} will be a  consequence of the following lemma. 
Here we exploit that     the sites $ k= (i,a) $ in $ {\cal S}_N $ are singular, see Definition \ref{regulars}. 

\begin{lemma}\label{thm:separation}
There is  $ C(d, \nu) > 0 $ such that, 
$ \forall \theta \in \R $, $ \forall N $, any $ M $-chain of sites in ${\cal S}_N $ 
has length  
\be\label{lenghtM} 
L \leq (M N)^{C(d, \nu)} \, . 
\ee 
\end{lemma}

\begin{pf}
Let $ k_q = (l_q,j_q,a_q) $, $ q = 0, \ldots , L  $, be a $ M$-chain of sites in $ {\cal S}_N $. Then  
\be\label{chain}
\max\{|l_{q+1} - l_q |, | j_{q+1} - j_q |\} \leq M \, , \quad  \forall q \in [0,L] \, , 
\ee
and, in particular, by Definition \ref{regulars} and (\ref{diao}), 
$$
|-\om \cdot l_q + \|j_q \|^2 + m - \theta | < \Theta \ \ {\rm (if} \ a_q = 1 {)} \quad {\rm or} \quad
|\om \cdot l_q + \|j_q \|^2 + m + \theta | < \Theta \ \ {\rm (if} \ a_q = 0 {)} \, .
$$
We deduce one of the following $\teta$-independent inequalities  
$$
| \pm \om \cdot (l_{q+1} - l_q) + (\|j_{q+1}\|^2 \pm \|j_q\|^2 ) | \leq 2( \Theta +m) \, . 
$$
By (\ref{chain}) we get $ | \|j_{q+1}\|^2 \pm \|j_q\|^2 | \leq  2(\Theta +m) + |\om| M \leq K_1 M $ for some
 $  K_1 := K_1 (\Theta, m) $. Since $| \|j_{q+1}\|^2 - \|j_q\|^2 | \leq 
 \|j_{q+1}\|^2 + \|j_q\|^2 $, in any case 
$ | \|j_{q+1}\|^2 - \|j_q\|^2 | \leq K_1M $.
Therefore
\be\label{jl}
\forall q, q_0 \in [0,L] \   , \   | \|j_{q}\|^2 - \|j_{q_0}\|^2 | \leq |q - q_0 | K_1 M  
\ee
and, using also  (\ref{chain}), 
\be
|j_{q_0} \cdot (j_{q}- j_{q_0})| =  \frac12 \Big| \| j_q \|^2 - \| j_{q_0} \|^2 -  \| j_q - j_{q_0}\|^2 \Big|
\leq K_2 | q - q_0|^2 M^2 \,  .  \label{jll0}
\ee
Let us introduce the subspace of $\R^d$
$$
G={\rm Span}_\R \{ j_{q}-j_{q'} \ : \ 0\leq q,q' \leq L \ \}=
{\rm Span}_\R \{ j_{q}-j_{0} \ : \ 0\leq q \leq L \ \}
$$
and let us call $g$ ($1\leq g \leq d$) the dimension of $ G $. Define  $ \d :=  (2d+1)^{-2}$. 
The constants  $C $ below (may) depend on  $ \Theta, m, d, \nu $.
\\[1mm]
{\bf Case I}. $ \forall q_0 \in [0, L] $, 
$ {\rm Span}_\R \{ j_{q}- j_{q_0} \, : \,  | q - q_0 | \leq L^\d \, ,  \ q \in [0,L] \, \} = G \, $.
\\[1mm]
We select a basis of $ G $ from  $ j_{q}- j_{q_0} $ ($|q-q_0| \leq L^\d$), say  $ f_1 \, , f_2 \, , \ldots \, , f_g  \in G $.
By (\ref{chain})  we have
\be\label{boundfi}
| f_i | \leq M L^\d \,  , \qquad \forall i =1, \ldots, g \, .  
\ee
Decomposing in this basis the orthogonal projection of $j_{q_0}$ on $G$, 
\be\label{compo}
P_G j_{q_0} = \sum_{i=1}^g x_i f_i 
\ee 
and taking the scalar products with $ f_p $, $ p =1, \ldots, g $, we get the linear system
$$ 
F x = b  \quad {\rm with} \quad \ F^{i}_p := f_i \cdot f_p  \, , \  b_p :=  P_G j_{q_0} \cdot f_p= j_{q_0} \cdot f_p  \, .
$$
Since $ \{  f_i \}_{i =1, \ldots , g} $ is a basis of $ G $ the matrix $ F $ is invertible. Since the coefficients 
of $ F $ are integers, $ |{\rm det}(F) | \geq 1 $. By Cramer rule, using that (\ref{boundfi}) implies 
$ |F^{i}_p| \leq C |f_i| |f_p| \leq ( M L^\d )^2 $, we deduce that 
\be\label{A-1}
|(F^{-1})_i^{i'}| \leq C ( M L^\d )^{2(g-1)} \, , \quad \forall i,i'  = 1, \ldots, g \, .
\ee
By (\ref{jll0}), we have $ |b_i | \leq K_2 (M L^{\d})^2 $, $ \forall i =1, \ldots, g $, and (\ref{A-1}) implies
\be\label{xb}
|x_{i'}| \leq C (ML^\d)^{2g} \, , \quad \forall i' =1, \ldots,  g \, . 
\ee
From (\ref{compo}), (\ref{boundfi}), (\ref{xb}),  we deduce 
$ |P_G j_{q_0} | \leq C (ML^\d)^{2g+1}$, $ \forall q_0 \in [0,L] $, 
and
$$
|j_{q_1}-j_{q_2}| = |P_G j_{q_1} - P_G j_{q_2} | \leq C (ML^\d)^{2g+1} \leq C (ML^\d)^{2d+1}, 
\quad \forall (q_1,q_2)  \in [0,L]^2 \, . 
$$
Since all the $ j_{q}$ are in  $ \Z^d $, 
their number (counted without multiplicity) does not exceed $ C (ML^\d)^{(2d+1)d} $.
Thus we have obtained the bound
\be\label{Nspace}
\sharp \{ j_{q} \ ; \ 0\leq q \leq L \} 
\leq C (ML^\d)^{(2d+1)d}  \, . 
\ee
Now by Lemma  \ref{Ntime}, for each $q_0 \in [0,L]$, the number of $q\in [0,L]$ such that
$j_q=j_{q_0}$  is  at most $ 2 N^{2d + \nu + 4} $, 
and so
$$
L \leq C (ML^\d)^{(2d+1)d} 2 N^{2d+\nu+4} \, . 
$$
Since $ \d (2d+1)d < 1 \slash 2 $,  we get 
\be\label{lunghezza}
L \leq M^{2d(d+1)} N^{2(2d+\nu+4)} 
\ee
for $N$ large enough, proving \eqref{lenghtM}. 
\\[1mm]
{\bf Case II.} There is $ q_0 \in [0,L] $ such that  
$$
 \mu := {\rm dim} \, {\rm Span} \{ j_{q} - j_{q_0} \, : \,  | q - q_0| \leq L^\d \, , \ q \in [0,L] \}  \leq g -1 \, ,
$$
namely all the vectors $ j_q $ stay in a affine subspace of dimension $ \mu \leq g - 1 $. 
Then we repeat on the sub-chain $ j_q $, $ | q - q_0 | \leq L^\d $, the argument of case I, 
to obtain a bound for $L^\d$ (and hence for $L$).

Applying  at most $ d $-times the above procedure, we obtain a bound for $L$ of the form 
$ L \leq  (M N)^{C(d,\nu)} $, proving the lemma.
\end{pf}

We introduce the following equivalence relation in ${\cal S}_N$.

\begin{definition}\label{equivalence} 
We say that $ x \equiv y $ if  there is a $M$-chain $ \{ k_q \}_{q = 0, \ldots, L} $
in ${\cal S}_N$ connecting $ x $ to $ y $, namely $ k_0 = x $, $ k_L = y $.
\end{definition}

\noindent 
{\sc Proof of Proposition \ref{prop:separation} completed.} 
Set $ M := 2 N^2 $. 
By the previous equivalence relation we get a partition 
$$
{\cal S}_N = \bigcup_\a \Omega_\a' 
$$
in disjoint equivalent  classes, satisfying, by Lemma \ref{thm:separation}, 
\be\label{sepaquasi}
{\rm d}(\Omega_\a', \Omega_\b') > 2 N^2 \, , \quad {\rm diam}(\Omega_\a') 
\stackrel{(\ref{lenghtM})} \leq 2 N^2 (2 N^3)^{C(d,\nu)} \, .
\ee
All the sites outsides $ {\cal S}_N $ are regular or $ N $-regular, see (\ref{sn}). 
As a consequence all the sites outside 
$$
\bigcup_\a	 \Omega_\a'' \qquad {\rm where} \qquad  \Om_\a'' := \Big\{ k \in \Z^b \times \{0,1\} 
\, : \, {\rm d}(k, \Om_\a') \leq N \Big\} 
$$
are $ N $-good, see (\ref{regNreg}). 
Hence the $ N $-bad sites (see Definition \ref{GBsite}) of $ A(\e, \l,\teta ) $ with $|l| \leq N' $ are included in 
$$
\bigcup_\a	 \Omega_\a 
\qquad {\rm where} \qquad  \Om_\a :=  \Om_\a'' \cap \{(l,j,a) \, : \, |l| \leq N' \}  \, .
$$
Then (\ref{separ}) follows by  (\ref{sepaquasi}) with
$ C_1 := 3 C(d,\nu) + 3 $, for $ N \geq N_0 (d,\nu,m,\Theta , \g_0 , \tau_0) $ large enough.

\section{Measure and ``complexity" estimates}\label{sec:measure}\setcounter{equation}{0}

We define
\begin{eqnarray}\label{tetabadweak}
B_{N}^0 (j_0;\e, \l) & := & \Big\{ \teta \in \R \, : \,   
\| A_{N, j_0}^{-1} (\e, \l,\theta)\|_0 > N^{\t}  \Big\} \\
& = &
\Big\{ \teta \in \R \, : \,  \exists {\rm \ an \ eigenvalue \ of \ }  A_{N, j_0} (\e, \l,\theta) \ {\rm with \ modulus \ less \ than} \ N^{-\t}  \Big\}
\label{la2}
\end{eqnarray}
where $ \|  \ \|_0 $ is the operatorial $ L^2 $-norm defined in (\ref{L2norm}).
The equivalence between (\ref{tetabadweak}) and (\ref{la2}) is a consequence of the self-adjointness of $A_{N, j_0} (\e, \l,\theta)$. 
We also define

\begin{eqnarray}\label{weakgood}
{\cal G}_{N}^0 := {\cal G}_{N}^0 ({\bf u}) & := & \Big\{ (\e , \l) \in [0,\e_0] \times  \Lambda \, : \,   
\forall \, j_0 \in \Z^d \, , \quad
B_{N}^0(j_0; \e, \l) \subset \bigcup_{q = 1, \ldots, N^{2d+\nu+4}} I_q \label{BNcomponent2} \\
& & \ \  {\rm where} \ I_q \ {\rm are \  disjoint \ intervals \ with \ measure} \ | I_q| \leq N^{-\t}   \Big\}  \, . \nonumber
\end{eqnarray}

\begin{remark}
The difference between the sets $ {\cal G}_N^0 $ defined in (\ref{weakgood}) 
and ${\cal G}_N $ 
defined in (\ref{good})
relies in the different definition of $ B_N^0 (j_0;\e, \l) $ in (\ref{tetabadweak}) and $ B_N (j_0;\e, \l) $ in (\ref{tetabad}).
For all $ \teta \notin B_N(j_0;\e, \l) $ the matrices $ A_{N,j_0}(\e, \l,\teta) $ are $ N $-good, i.e. satisfy bounds on 
$ \nors{A_{N,j_0}^{-1}(\e, \l,\teta)} \leq N^{\d s+ \t'} $ 
for  $ s \in [s_0,s_1] $, while for all $ \teta \notin B_N^0(j_0;\e, \l) $ we only have the $ L^2$- bound 
$ \| A_{N,j_0}^{-1}(\e, \l,\teta)\|_0 \leq N^\tau $. Using the multiscale Proposition \ref{propinv} and the 
separation Proposition \ref{prop:separation} (which holds for any  $ \teta $) we shall prove 
inductively that the parameters that stay in ${\cal G}_{N_k}^0 (u_k)$ along the Nash-Moser scheme are in fact also 
in ${\cal G}_{N_k} (u_k)$.
\end{remark}

The aim of this section is to prove the following proposition. 

\begin{proposition}\label{PNmeas}
There is a constant $C>0$ such that,
for $ N \geq N_0 (V,d,\nu) $ large enough and 
\be\label{ipopicco}
\e_0 \b_0^{-1} (\| T_1 \|_0 + \| \partial_\l T_1 \|_0) \leq c 
\ee 
small enough ($ \b_0 $ is defined in \eqref{eq:posi} and $ T_1 $ in \eqref{T1}$)$, 
the set $ {\cal B}_{N}^0 := ({\cal G}_{N}^0 )^c \cap ([0,\e_0] \times \Lambda) $ 
has measure 
\be\label{measBN0}
|{\cal B}_{N}^0 | \leq C \, \e_0 N^{-1} \, .
\ee
\end{proposition}

Proposition \ref{PNmeas} is derived from several lemmas based on  
basic properties of eigenvalues of self-adjoint matrices,
which are a consequence of their variational characterization.  

\begin{lemma}\label{variatione}
i) Let $ A(\xi) $ be a family of  square matrices in ${\cal M}^E_E $, $ C^1 $ in the real parameter $ \xi  \in \R $. 
Assume that there  is an invertible matrix $ U $ such that the matrices $\wtilde{A}(\xi) := A(\xi) U $ are self-adjoint
and $ \partial_\xi \wtilde{A}(\xi) \geq \b I $, $ \b >  0  $. Then, for any $ \alpha >0 $, the measure 
\be\label{me1}
\Big| \Big\{ \xi \in \R \, : \, \|A^{-1}(\xi)\|_0 \geq \alpha^{-1} \Big\}\Big| \leq  2 |E| \alpha \beta^{-1} \|U\|_0 
\ee
where $ | E | $ denotes the cardinality of the set $ E $.
\\[1mm]
ii) In particular, if $ A = Z + \xi W $ with $ Z, W $ selfadjoint, $ W $ invertible and 
$ \beta_1 I \leq Z \leq \beta_2 I $, $ \beta_1 > 0 $, then
\be\label{me2}
\Big| \Big\{ \xi \in \R \, : \, \|A^{-1}(\xi)\|_0 \geq \alpha^{-1}  \Big\}\Big| \leq  2 |E| \alpha \beta_2 \beta_1^{-1} \|W^{-1}\|_0 \, . 
\ee
\end{lemma}

\begin{pf}
i)  The eigenvalues of the self-adjoint matrices 
$\wtilde{A}(\xi)$ can be listed as $ C^1$ functions $ \mu_k (\xi ) $, $ 1 \leq k \leq |E| $. Now
\begin{eqnarray*}
\Big\{ \xi \in \R \, : \, \|A^{-1}(\xi)\|_0 \geq \alpha^{-1}  \Big\} &\subset&
\Big\{ \xi \in \R \, : \, \|\wtilde{A}^{-1}(\xi)\|_0 \geq (\alpha \|U\|_0)^{-1}  \Big\} \\ 
& = &
\Big\{ \xi \in \R \,  :  \, \exists k \in [1, |E|] \,  , \, |\mu_k(\xi)| \leq \alpha \|U\|_0 \Big\}
\end{eqnarray*}
because $\wtilde{A}(\xi) $ is selfadjoint. Since
$ \partial_\xi \wtilde{A}(\xi) \geq \b I $, we have
$ \partial_\xi \mu_k(\xi) \geq \b >  0 $ and the measure estimate \eqref{me1} follows readily. 

ii) Applying i) with $ U = W^{-1}Z $ and self-adjoint matrices $ \wtilde{A}(\xi) = ZW^{-1}Z + \xi Z $, 
we get 
$$
\Big| \Big\{ \xi \in \R \, : \, \|A^{-1}(\xi)\|_0 \geq \alpha^{-1}  \Big\}\Big| \leq 2 |E| \alpha \beta_1^{-1} 
\|W^{-1}\|_0 \|Z\|_0 \leq  2 |E| \alpha \beta_2 \beta_1^{-1} \|W^{-1}\|_0 ,
$$
which is \eqref{me2}.
\end{pf}

From the variational characterization of the eigenvalues
of selfadjoint matrices we can derive :

\begin{lemma}\label{Lips}
Let $ A $, $ A_1 $ be self adjoint matrices. Then their eigenvalues (ranked in nondecreasing order)
satisfy the Lipschitz property
\be\label{continua}
|\mu_k (A) - \mu_k (A_1)| \leq \| A - A_1 \|_0  \, .
\ee
\end{lemma}

The  continuity property (\ref{continua}) 
of the eigenvalues allows to derive  a ``complexity estimate" for $ B_N^0 (j_0;\e, \l) $
knowing its measure, more precisely the measure of  
\be\label{B2N0}
B_{2,N}^0(j_0; \e, \l) :=
\Big\{   \teta \in \R \, : \,  \| A_{N, j_0}^{-1} (\e, \l,\theta)\|_0 >  N^{\t} / 2  \Big\} \, .
\ee

\begin{lemma}\label{complessita}
$ \forall j_0 \in \Z^d $, $ \forall (\e,\l) \in [0,\e_0] \times \Lambda $, we have 
$ B^0_N (j_0; \e, \l) \subset \cup_{q=1,..., 2 \, {\mathtt M} N^\tau } I_q $
where $ I_q $ are  intervals with $ |I_q| \leq N^{-\tau } $ and  
$ {\mathtt M} := | B_{2,N}^0(j_0; \e, \l) |$.
\end{lemma}

\begin{pf}
If $ \theta \in  B_{N}^0(j_0;\e, \l)  $, by   (\ref{continua}) and since $ \| Y \|_0 = 1 $ (see (\ref{def:Y})), 
we deduce that 
\begin{eqnarray}
\Big[ \theta -  N^{-\t} , \theta +  N^{-\t} \Big]  & \subset &  B_{2,N}^0(j_0;\e, \l) \nonumber   \\ 
& = & 
\Big\{ \teta \in \R \, : \,  \exists {\rm \ an \ eigenvalue \ of \ }  
A_{N, j_0} (\e, \l,\theta) \ {\rm with \ modulus \ less \ than} \ 2N^{-\t}  \Big\} .  \nonumber
\end{eqnarray} 
Hence $ B_{N}^0(j_0;\e, \l)  $ is included in an union of  intervals $J_m$ of disjoint interiors,
\be\label{inequa}
B_{N}^0 (j_0;\e, \l)  \subset \bigcup_{m} J_m \subset B_{2,N}^0(j_0;\e, \l), \quad {\rm with \ length} \quad  
|J_m|  \geq  2 N^{-\t}
\ee
(if some of the intervals $ [\theta - N^{-\t} , \theta + N^{-\t} ]$ overlap, then we glue them together).  
We decompose each 
$ J_m $  as an union of (non overlapping) intervals $ I_q $ of length between 
$ N^{-\tau}/2 $ and  $ N^{-\tau} $. Then, 
by (\ref{inequa}), we  get a new covering  
$$
B_{N}^0 (j_0;\e, \l)  \subset \bigcup_{q=1, \ldots, Q} I_q  \subset B_{2,N}^0(j_0;\e, \l) \quad  {\rm with} \ \   
N^{-\t} / 2 \leq |I_q | \leq N^{-\t} 
$$
and, since the intervals $ I_q $ do not overlap, 
$$
Q N^{- \t}  \slash 2 \leq \sum_{q = 1}^Q | I_q |   
\leq |  B_{2,N}^0(j_0;\e, \l) | =:   {\mathtt M} \, . 
$$
As a consequence
$ Q \leq 2 \,  {\mathtt M} \, N^{\t} $, which proves the lemma.
\end{pf}

\noindent

We estimate the measure $ |B_{2,N}^0(j_0; \e, \l)| $  differently for  $ |j_0| \geq 2N $ or $ |j_0| < 2 N $.
In the next lemmas we assume 
\be\label{assum1}
N \geq N_0 (V,\nu,d) > 0 \ {\rm large \ enough} \quad {\rm and} \quad 
 \e \| T_1 \|_0 \leq 1 \, .
\ee

\begin{lemma}\label{sectio1} 
$ \forall |j_0 | \geq 2 N $, $ \forall (\e,\l) \in [0,\e_0] \times \Lambda $,  we have
$ |B_{2,N}^0(j_0; \e, \l)| \leq CN^{- \t +d +\nu} $.
\end{lemma}

\begin{pf}
Recalling  (\ref{ANj0})  and (\ref{matrpar}), we have
\be\label{svipi}
A_{N,j_0}(\e, \l, \theta) = A_{N, j_0}(\e, \l) + \theta Y_{N, j_0} = 
D_{N,j_0} (\l) + T_{N,j_0}(\e, \l) + \teta Y_{N, j_0}. 
\ee 
We claim that, if $ |j_0| \geq 2 N  $ and $  N \geq N_0(V,d,\nu)$, see \eqref{assum1}, then  
\be \label{positiv} 
4 d|j_0|^2 I \geq A_{N, j_0}(\e, \l) \geq \frac{|j_0|^2}{8} I \, . 
\ee
Indeed by (\ref{svipi}) and \eqref{continua}, the eigenvalues  $ \lambda_{l,j}  $ of $A_{N, j_0}(\e, \l)$ satisfy
\be \label{vpA}
\lambda_{l,j} = \d_{l,j}^\pm  + O(\e\|T_1\|_0 + \|V \|_0)  \quad {\rm where} \quad
\d_{l,j}^\pm :=   \|j\|^2  \pm \om \cdot l \, .
\ee
Since $ |\om | = |\l| | \bar \o | \leq 3/2 $ (see (\ref{baromega})), $\|j\| \geq |j|$ (see \eqref{supeuc}),  
$|j-j_0| \leq N $, $ |l| \leq N $, 
we have 
\be\label{lo1}
\d_{l,j}^\pm    \geq (|j_0| - |j-j_0|)^2   - \nu |\o | |l| \geq 
(|j_0| - N)^2   -  \frac32  \nu N \geq \frac{|j_0|^2}{6}
\ee
for $ |j_0| \geq 2 N  $ and $ N \geq N_0(\nu) $ large enough. 
Moreover, since $\|j\|^2 \leq d |j|^2$, 
\be\label{up1}
\d_{l,j}^\pm  \leq  d(|j_0| + |j-j_0|)^2 + \nu|\o | |l| \leq d(|j_0| + N)^2 +  2 \nu N \leq 3d |j_0|^2  
\ee
for $ N \geq N_0(\nu) $ large enough. 
Hence \eqref{vpA}, \eqref{lo1}, \eqref{up1}, \eqref{assum1} imply (\ref{positiv}). 
As a consequence, by Lemma \ref{variatione}-ii) with $ W = Y_{N,j_0} $, $ \|W^{-1}\|_0 = 1 $, 
we deduce  $|B_{2,N}^0(j_0; \e, \l)| \leq C N^{- \t +d +\nu} $.
\end{pf}

Lemmas \ref{complessita} and \ref{sectio1} imply that:

\begin{corollary}\label{lem:complexity} 
$ \forall |j_0| \geq 2N  $, $ \forall (\e,\l) \in [0,\e_0] \times \Lambda $,  we have
$$ 
B_{N}^0 (j_0; \e, \l) \subset \bigcup_{q=1, \ldots, N^{d+\nu+2}} I_q 
$$
where $ I_q $ are intervals  satisfying $  |I_q | \leq N^{- \t} $. 
\end{corollary}

We now consider the cases $ |j_0| < 2 N $.

\begin{lemma}\label{corol1}
$ \forall | j_0 | < 2 N $, $ \forall (\e,\l) \in [0,\e_0] \times \Lambda $, we have
$$ 
B_{2,N}^0(j_0; \e, \l) \subset I_N := ( - 11d N^2, 11d N^2) \, .
$$
\end{lemma}

\begin{pf}
The eigenvalues of $ \teta Y $ are $ \pm \theta $ and  \eqref{supeuc}  implies
$\|j\| ^2 \leq d (|j_0  | + |j-j_0| )^2 \leq 9 d N^2 $. Hence, by  (\ref{svipi}), (\ref{vpA}), 
$ |l| \leq N $, (\ref{baromega}), \eqref{assum1},   
$$
\| A_{N, j_0}(\e, \l) \|_0 \leq \| D_{N, j_0}(\l) \|_0 +  \| T_{N, j_0}(\e, \l) \|_0 \leq 
2 \nu  N + 9d N^2 + C(1+\| V \|_0) \leq 10 dN^2   
$$
for $ N \geq N(V,d,	\nu) $ large enough.  By Lemma \ref{Lips}, 
if $  \teta \notin I_N $ all the eigenvalues of 
$ A_{N,j_0}(\e, \l, \teta) =  A_{N,j_0}(\e, \l) + \theta Y_{N,j_0}  $ are greater than $ 1 $ (actually $ d N^2 $).
\end{pf}

\begin{lemma}\label{cor2} $ \forall |j_0| < 2 N $, the set 
\be\label{B22N}
{\bf B}^{0}_{2,N}(j_0) := \Big\{  (\e, \l, \teta)  \in [0,\e_0] \times \Lambda \times {\R} \, : \, 
\Big\| A_{N,  j_0}^{-1}(\e, \l,\teta) \Big\|_0 > N^{\t}/2  \Big\}
\ee
has measure 
\be\label{B2N} 
|{\bf B}^{0}_{2,N}(j_0)| \leq \e_0 N^{-\t+d+\nu+3}  \, . 
\ee 
\end{lemma}

\begin{pf}
By Lemma \ref{corol1},   $ {\bf B}^{0}_{2,N}(j_0)  \subset [0,\e_0] \times \Lambda \times I_N $.
In order to estimate the  ``bad" $ (\e, \l, \teta) $ where at least one eigenvalue of 
$ A_{N,j_0}(\e, \l, \teta )$ is less than $ N^{-\t} $, we introduce the variables
\be\label{changevaria}
\xi := \frac{1}{\l} \, , \quad \eta := \frac{\teta}{\l} 
\quad {\rm where} \quad (\xi, \eta) \in [2 \slash 3, 2] \times 2 I_N  
\ee
and we consider the self adjoint matrix
\be\label{etaxi}
\frac{1}{\l} A_{N,j_0} (\e, \l,\teta) 	\stackrel{(\ref{svipi})} =
{\rm diag}_{|l| \leq N, |j-j_0|\leq N} 
\left(\begin{array}{cc} 
 - \, \bar \o \cdot  l  & 0 \\ 
  0 &  \bar \o \cdot l  \end{array}\right) \,
+ \xi  P_{N,j_0}
\, - \e \xi T_1(\e , 1 \slash \xi )  + \eta Y 
\ee
where 
$$
P := \left(\begin{array}{cc} 
 - \Delta + V(x)  & 0\\ 0 &
 - \Delta + V(x)  \end{array}\right) \quad
{\rm satisfies} \quad P \stackrel{\eqref{eq:posi}} \geq \b_0 I  \, .
$$
The derivative with respect to $ \xi $ of the matrix in (\ref{etaxi}) is 
$$
P_{N,j_0} - \e T_1(\e, 1/ \xi ) + \frac{\e}{\xi} \partial_\l T_1( \e, 1/ \xi ) 
\stackrel{\eqref{ipopicco}} \geq \frac{\b_0}{2} I \, , 
$$
i.e. positive definite (for $\e_0$ small enough).
By Lemma \ref{variatione}, for each fixed 
$ \eta $, the set of $ \xi \in [2/3,2] $ such that at least one eigenvalue is $ \leq  N^{-\t} $ has measure
at most $ O(  N^{- \t + d + \nu })$. Then, integrating on $ \eta \in I_N  $, whose length is $ |I_N| = O(N^2)$, on 
$\e \in [0,\e_0 ] $,  
and since the change of variables (\ref{changevaria}) has a Jacobian of modulus
 $\geq  1 / 8 $, we deduce (\ref{B2N}).
\end{pf}

By the same arguments (see also the proof of Lemma \ref{NRfre})
we also get the following measure estimate that will be used in section \ref{sec:NM}, see $ (S4)_n $.

\begin{lemma}\label{measure0} 
The complementary of the set 
\be\label{Binver}
{\mathtt G}_{N} := {\mathtt G}_N({\bf u}) := \Big\{  (\e,\l)   \in [0,\e_0] \times \Lambda \, : \, \| A_{N}^{-1}( \e, \l) \|_0 \leq  N^{\t}   \Big\}
\ee
has measure  
\be\label{measGN0} 
| {\mathtt G}^c_{N} \cap ([0,\e_0] \times \Lambda) | \leq  \e_0 N^{- \t +d+\nu+1} \, . 
\ee 
\end{lemma}

\begin{remark}
For periodic solutions (i.e. $ \nu = 1 $), a similar eigenvalue variation argument which exploits $ - \Delta \geq 0 $ was used 
in the Appendix of \cite{B3} and in \cite{BP}.
\end{remark}

As a consequence of Lemma \ref{cor2}, for ``most" $ (\e, \l)$ the measure  of 
$ B_{2,N}^0 (j_0;  \e, \l ) $ is ``small".

\begin{lemma}\label{intermed}
$ \forall |j_0| < 2 N $, the set
$$
{\cal F}_{N}(j_0) :=  \Big\{ (\e,\l) \in [0,\e_0] \times \Lambda \, : \, |B_{2,N}^0 (j_0; \e, \l)| \geq 
\frac12 N^{- \t + 2d + \nu + 4} \Big\}
$$
has measure 
\be\label{sectio}
| {\cal F}_{N}(j_0)| \leq 2 \e_0 N^{- d -1 } \, .
\ee
\end{lemma}

\begin{pf}
By Fubini theorem (see \eqref{B22N} and \eqref{B2N0})
\be\label{Fubini}
|{\bf B}^{0}_{2,N}(j_0)| = \int_{[0,\e_0] \times \Lambda} | B_{2,N}^0(j_0; \e, \l) | d \e \, d \l  \, .
\ee
Let  $ \mu := \t - 2 d - \nu - 4  $.  By (\ref{Fubini}) and (\ref{B2N}), 
\begin{eqnarray}
\e_0 N^{-\t+d+ \nu + 3} & \geq & 
\int_{[0,\e_0] \times \Lambda} |B_{2,N}^0 (j_0; \e, \l) | d \e \, d \l \nonumber \\
& \geq & \frac12 N^{-\mu} \Big| \Big\{ (\e,\l) \in 
[0,\e_0] \times \Lambda  \, : \, |B_{2,N}^0 (j_0; \e, \l)| \geq \frac12 N^{-\mu} \Big\} \Big| \nonumber 
:=  \frac12 N^{-\mu} |{\cal F}_{N}(j_0)| \nonumber
\end{eqnarray}
whence (\ref{sectio}).
\end{pf}

By Lemma \ref{intermed}, for all $ (\e, \l) \notin {\cal F}_{N}(j_0) $ we have  the measure 
estimate $ |B_{2,N}^0(j_0; \e, \l)| <  N^{- \t + 2d + \nu + 4} / 2 $. 
Then, Lemma \ref{complessita}  implies

\begin{corollary}\label{lem:complexity1} 
$ \forall |j_0| < 2N  $, $ \forall (\e,\l) \notin {\cal F}_N(j_0) $, we have  
$ B_{N}^0 (j_0; \e, \l) \subset \bigcup_{q=1, \ldots, N^{2d+\nu+4}} I_q $ with $ I_q $ intervals 
satisfying $ |I_q | \leq N^{- \t} $. 
\end{corollary}

Proposition \ref{PNmeas}  is a direct consequence of the following lemma. 

\begin{lemma}\label{inclusion}
$ {\cal B}_{N}^0 \subseteq \bigcup_{|j_0| < 2 N} {\cal F}_{N}(j_0)  $.
\end{lemma}

\begin{pf}
Corollaries \ref{lem:complexity}  and \ref{lem:complexity1} imply that
$$ 
(\e,\l) \notin \bigcup_{|j_0| < 2 N} {\cal F}_{N}(j_0)   \quad \Longrightarrow \quad 
(\e,\l) \in  {\cal G}_N^0 
$$
(see the definition in (\ref{weakgood})).
The lemma follows. 
\end{pf}

\noindent
{\sc Proof of Proposition \ref{PNmeas} completed}.  
By Lemma \ref{inclusion} and (\ref{sectio}) we get
\be\label{meafin}
| {\cal B}_{N}^0 | \leq \sum_{|j_0| < 2 N} |{\cal F}_{N}(j_0)|  <
(2 N+1)^d  |{\cal F}_{N}(j_0)|  \leq (2 N+1)^d 2 \e_0  N^{-d-1} \leq  C \e_0 N^{-1} \, .
\ee

\section{Nash Moser iterative scheme}\label{sec:NM}\setcounter{equation}{0}

Consider the orthogonal splitting 
$$ 
{\bf H}^s = H_n \oplus H_n^\bot
$$
where $ {\bf H}^s $ is defined in \eqref{bHs} and 
\begin{eqnarray}
H_n & := & \Big\{ u := {\bf u} = (u^+,u^-) \in {\bf H}^s \, : \, u = \sum_{|(l,j)| \leq N_n} \, u_{l,j}  \, e^{\ii (l \cdot \vphi + j \cdot x)}  
\Big\} \label{Hn} \\
H_n^\bot & := & \Big\{ u  := {\bf u}  = (u^+,u^-) \in {\bf H}^s \,  : \, u = \sum_{|(l,j)|> N_n} \,   u_{l,j} \, e^{\ii (l \cdot \vphi + j \cdot x)} \Big\} \, , \nonumber
\end{eqnarray}
with $ u_{l,j} := ( u^+_{l,j}, u^-_{l,j}) \in \C^2 $, and 
\be\label{defNn}
N_n := N_0^{2^n} \, ,  \ \qquad   {\rm namely} \ \ \ N_{n+1} = N_n^2  \, ,  \ \forall n \geq 0 \, . 
\ee
In the proof we shall take  $ N_0 \in \N $ large enough depending on $ \e_0 $ and $ V $, $ d $, $ \nu $, see \eqref{N0ge}. 
We denote by 
\be\label{Pn}
P_n : {\bf H}^s \to H_n \qquad \qquad {\rm and} \qquad \qquad   
P^{\bot}_n: {\bf H}^s \to H_n^\bot
\ee
the orthogonal projectors onto $ H_n $ and $ H_n^{\bot} $. 
The following ``smoothing" properties hold, $ \forall n \in \N $, $ s \geq 0 $, $ r \geq 0 $,  
\begin{eqnarray} 
& & \| P_n u \|_{s+r}  \leq  N_n^r \| u \|_s \, , \qquad  \forall u \in {\bf H}^s \label{S1} \\
&& \|  P_n^\bot u \|_{s}  \leq   N_n^{- r } \| u \|_{s+r} \, , \quad  \,  \forall u \in {\bf H}^{s+r} \, . \label{S2}
\end{eqnarray}
More  generally, for $ j_0 \in \Z^d $, we denote $ P_{N,j_0} $
the orthogonal projector   from $ {\bf H}^{s} $ onto  the subspace
\be\label{HNj0}
H_{N,j_0} := \Big\{ u \in {\bf H}^s \, : \, u = \sum_{|(l,j-j_0)| \leq N} \, u_{l,j}  \, e^{\ii (l \cdot \vphi + j \cdot x)}  
 \Big\} \, .
\ee
With the above notation
$ H_n = H_{N_n,0} $, see \eqref{Hn}, and $ P_n := P_{N_n,0} $, see \eqref{Pn}. 
Moreover we also denote
$ \Pi_{N,j_0} $ the orthogonal projector from $H^{s_0}(\T^d)$  (functions only of the $ x $-variable)
onto  the space
\be\label{E0}
E_{N,j_0} := \Big\{ u(x) := \sum_{|j-j_0| \leq N}  u_j e^{\ii j \cdot x} \, , \ u_j \in \C  \Big\} \, .
\ee
The  composition operator on  Sobolev spaces
$$
f: {\bf H}^s \to {\bf H}^s \, , \qquad 
f(u) (t,x) := \left(\begin{array}{c}  f(\vphi,x,u^- u^+ ) u^+  \\  f(\vphi,x,u^- u^+ ) u^-    \end{array}\right)  \, , 
$$
where $ f \in C^q (\T^\nu \times \T^d \times \R; \R)$ with  
\be\label{defk} 
q \geq S + 2 
\ee
satisfies the following standard properties (see e.g.  \cite{LM}): $ \forall s \in [s_1, S] $,
$ s_1 > (d+\nu)/2 $, 
\begin{itemize}
\item {\bf (F1)} ({\bf Regularity}) 
$ f \in C^2 ({\bf H}^{s}; {\bf H}^{s} )$. 
\end{itemize}
\begin{itemize}
\item {\bf (F2)} ({\bf Tame estimates}) 
$ \forall u, h \in {\bf H}^{s} $ with $ \| u \|_{s_1} \leq 1 $, 
\be\label{fDftame} 
\| f(u) \|_{s}  \leq C(s) \| u \|_{s} \, , \ \  \| (Df)(u) h \|_s 
\leq C(s) ( \| h \|_s + \| u \|_{s} \| h \|_{s_1}) \, .
\ee
\be\label{secondorder}
\| D^2 f(u)[h,v] \|_{s} \leq C(s) \Big(\|u\|_{s} \|h\|_{s_1} \|v\|_{s_1} + \|v\|_{s} \|h\|_{s_1} + \|v\|_{s_1} \|h\|_{s} \Big) 	\, . 
\ee
\end{itemize}
As a consequence we get
\begin{itemize} 
\item {\bf (F3)} 
({\bf Taylor Tame estimate}) 
$ \forall u \in {\bf H}^{s} $ 
with $ \| u \|_{s_1} \leq 1 $,  $ \forall h  \in {\bf H}^{s} $ with $\| h \|_{s_1} \leq 1$, 
\be\label{quadratic}
\| f( u + h) - f (u) -  (Df)(u) \, h \|_{s}  \leq  C(s) ( \| u\|_{s}  \| h \|_{s_1}^2 + \| h \|_{s_1}   \| h \|_{s}) \, .
\ee
In particular, for $ s = s_1 $, 
\be\label{P3s}
\| f (u+ h) - f (u) - (Df)(u) \, h \|_{s_1} \leq C(s_1) \| h \|_{s_1}^2 \, . 
\ee
\end{itemize}
The values of the constants $ s_1 $ and $ S $ are fixed in \eqref{Sgr} below. 

\begin{remark}
The differential $ (Df)(u) $ is the operator  $ T_1 $
defined in \eqref{T1} with $ (p,q) $ as in \eqref{pq}. 
\end{remark}

By Lemma \ref{lem:multi} and the first inequality in \eqref{fDftame}
applied to the composition operators in (\ref{pq}),
the T\"oplitz matrix $ T_1 $ which represents $ Df(u) $ satisfies, $ \forall s \in [s_1,S] $, 
\be\label{decayTu}
\norma T_1 \norma_s = \norma (Df)(u) \norma_s 
\leq C(s) (1+ \| u  \|_{s}) \, .
\ee
For simplicity of notation we  denote $(g, \bar g)$ simply by $ g $.  We shall use that 
$ g $ and the potential $ V $ satisfy
\be\label{gk}
\| g \|_{C^q} \leq C \,  , \quad 
\| V \|_{C^q} \leq C \, , 
\ee
for some fixed constant $ C $. 

With the above more concise notations, the vector NLS-equation \eqref{vnls} becomes 
\be\label{riscritta}
L_\om u = \e (f(u)+g) \, . 
\ee
For definiteness we fix the Sobolev indices
$ s_0 < s_1 < S  $
as
\be\label{Sgr}
s_0 := b = d + \nu \, , \qquad  
s_1 := 10 (\t + b)C_2 \, , \qquad 
S := 12 \t' +  8(s_1 + 1)   \, ,
\ee
where  
\be\label{tautau0}
C_2 :=  6(C_1 + 2)  \,  , \ \t :=   \max\{ d + \nu + 2, 2 \, C_2 \,  \t_0 + 1 \} \, ,
   \  \t' := 3 \t+ 2 b \, , \  \tau_0 := \nu 
\ee
(the constant $ \tau_0   $ is introduced in (\ref{diophan0}))
and $ C_1 := C_1(d,\nu) \geq 2 $ is defined in Proposition \ref{prop:separation}.
Note that $ s_0, s_1, S $ defined in \eqref{Sgr} 
depend only on $ d $ and $ \nu $.

We  also fix the constant $ \d $  in Definition \ref{goodmatrix} as
\be\label{delta14}
\d := 1/ 4  \, .
\ee

\begin{remark}\label{defchi}
By (\ref{Sgr})-(\ref{delta14})  the hypotheses (\ref{dtC})-(\ref{s1}) of Proposition  \ref{propinv} 
are satisfied for any $ \chi \in [C_2, 2 C_2) $, 
as well as assumption (ii) of Proposition \ref{prop:separation}. 
We assume $ \t \geq d + \nu + 2 $ in view of  (\ref{measGN0}). 
The strongest condition for $ S $ appears in the proof of  Lemma  \ref{deriva}.
\end{remark}

Setting 
$$ 
\t_1 :=  d + \nu 
$$ 
and $ \gamma > 0 $, we shall implement the first steps 
of the Nash-Moser iteration restricting $ \l $ to the set
\begin{eqnarray}\label{diofs}
\bar {\cal G} & := &  \Big\{ \l  \in \Lambda \, : \, 
\Big\| \Big(\pm \l \bar \om \cdot l + \Pi_{0} (-\D+ V(x))_{| E_{0}}\Big)^{-1} \Big\|_{L^2_x} 
\leq \frac{N_0^{\t_1}}{\g}, 
 \,  \forall  \, |l| \leq N_0   \Big\} \nonumber \\
& = & 
\Big\{ \l  \in \Lambda \, : \, 
|\pm \lambda \bar \om \cdot l +  \mu_j | \geq \g N_0^{-\t_1}, 
\, \forall  \, |j| \leq N_0, \, |l| \leq N_0    \Big\} 
\end{eqnarray}
where $ \mu_j $ are the eigenvalues of $ \Pi_{0} (- \D + V(x) )_{| E_{0}} $ where
$ \Pi_0 := \Pi_{N_0 ,0} $, $ E_0 := E_{N_0 ,0} $ are defined in  \eqref{E0}.
We shall prove in Lemma \ref{NRfre}  the measure bound  $ |\bar {\cal G}| = 1 - O(\g) $ (since $ \t_1 \geq d + \nu  $).
The constant $ \g $ will be fixed in \eqref{N0ge}.

We also define
\be\label{def:sigma}
\s := \t' + \d s_1 + 2 \, . 
\ee
Given a set $ A $ we denote $ {\cal N}(A, \eta) $ the open neighborhood of $ A $ of width $ \eta $ (which is empty if $ A $
is empty).

\begin{theorem}\label{cor1} {\bf (Nash-Moser)} 
There exist $ \bar c , \bar \gamma > 0 $ $($depending on $ d,\nu,V$,$\g_0,\beta_0)$ such that, if 
\be\label{smallsto}
N_0 \geq 2 \g^{-1}  \, , \  \gamma \in (0, \bar \g)  \, , 
\qquad { and } \qquad  
\e_0 N_0^{S} \leq \bar c  \, , 
\ee
then there is a sequence $(u_n)_{n \geq 0} $ of $ C^1$ maps 
$ u_n:$ $ [0,\e_0) \times \Lambda \to  {\bf H}^{s_1} \cap {\cal U} $
 (see \eqref{subvs}) satisfying 
\begin{itemize}
\item[$ {\bf (S1)}_n $] \ $ u_n (\e, \l ) \in H_n \cap {\cal U} $,  $ u_n ( 0, \l ) = 0 $, 
$ \| u_n \|_{s_1} \leq 1$,  $ \| \partial_{(\e,\l)} u_n \|_{s_1} \leq C(s_1) N_0^{\t_1+ s_1+ 1} \g^{-1} $.
\item[$ {\bf (S2)}_n $]  \ $ (n \geq 1) $ \  For all  $ 1 \leq k \leq n $, 
$\|u_k-u_{k-1}\|_{s_1} \leq N_k^{-\s -1 } $, 
$ \|\partial_{(\e,\l)}(u_k-u_{k-1})\|_{s_1} \leq N_k^{-1/2}$. 
\item[$ {\bf (S3)}_n $] \  $ (n \geq 1) $ 
\be\label{GN0N}
\| u - u_{n-1} \|_{s_1} \leq N_n^{-\sigma} \quad  \Longrightarrow \quad
\bigcap_{k= 1}^n {\cal G}_{N_k}^0 (u_{k-1})  \subseteq  {\cal G}_{N_{n}}(u) 
\ee
where  
$ {\cal G}^0_{N } (u ) $ (resp. ${\cal G}_{N } (u ) $) is defined in  (\ref{weakgood}) (resp. in (\ref{good}))  . 
\item[$ {\bf (S4)}_n $]
Define the set
\be\label{Gscavo}
{\cal C}_n := 
\bigcap_{k= 1}^n {\mathtt G}_{N_k}(u_{k-1})  \bigcap_{k= 1}^n 
{\cal G}_{N_k}^0(u_{k-1})   \bigcap \Big([0,\e_0] \times
{\bar {\cal G}}\Big) \, , 
\ee
where  ${\mathtt G}_{N_{k}}(u_{k-1}) $ is defined in (\ref{Binver}), $ \bar {\cal G} $ in (\ref{diofs}),
$ {\cal G}_{N_k}^0(u_{k-1}) $ in (\ref{weakgood}). 
 
If $ (\e, \l ) \in {\cal N} ( {\cal C}_n, N_n^{-\s}) $ then  
$  u_n(\e,\l) $ solves the equation
$$
P_{n} \Big(L_{\om} u - \e (f( u ) + g) \Big)= 0  \, . \leqno{(P_n)}   
$$
\item[$ {\bf (S5)}_n $] \   $ U_n :=  \|  u_n \|_{S} $,  
$ U_n' :=  \|  \partial_{(\e,\l)}{ u}_n \|_S $  (where $ S $ is defined in (\ref{Sgr}))
satisfy
$$
(i) \;\; U_n \leq  \displaystyle  N_n^{ 2( \t' + \d s_1 + 1)} \, , 
\qquad  (ii) \;\; U_n' \leq  N_{n}^{4\t'+ 2 s_1 + 4}  \, . 
$$
\end{itemize}
The sequence $ (u_n)_{n \geq 0} $  converges 
in $ C^1 $ norm to a map
\be\label{uC1} 
u \in C^1 ([0,\e_0) \times \Lambda ,  {\bf H}^{s_1}) \quad {\rm with} \quad u(0, \l) = 0 
\ee
and,  if $ (\e,\l) $ belongs to  the Cantor like set
\be\label{Cinfty}
{\cal C}_\infty := \bigcap_{n \geq 0} {\cal C}_n 
\ee
then $ u(\e,\l) $  is a solution of (\ref{vnls}), i.e.  \eqref{riscritta}, with $ \om = \l \bar \om $.
\end{theorem}

The sets of parameters $ {\cal C}_n $  in $ (S4)_n $ are decreasing, i.e.
$$
\ldots \subseteq {\cal C}_{n} \subseteq {\cal C}_{n-1} \subseteq \ldots \subseteq  {\cal C}_{0} \subset 
[0,\e_0] \times \bar {\cal G} \subset [0,\e_0] \times \Lambda \, , 
$$
and it could happen that $ {\cal C}_{n_0} = \emptyset $ for some $ n_0 \geq 1 $. 
In such a case $ u_n = u_{n_0} $, $ \forall n \geq n_0 $
(however the map $ u $ in \eqref{uC1} is always defined), and $ {\cal C}_\infty = \emptyset  $.
Later, in (\ref{N0ge}), we shall specify the values of $ \g, \e_0, N_0 $, in order to
verify that $ {\cal C}_\infty   $ has asymptotically full measure, i.e. (\ref{Cmeas}) holds.

The proof of  Theorem \ref{cor1} 
is based on an improvement of the  Nash-Moser theorems  in \cite{BB07}, \cite{BB}, \cite{BBP}.
The main difference 
is that  the ``tame exponent" $\t' + \d s $ in (\ref{os2}) depends on the Sobolev index $ s $.
We have chosen
$ \d = 1 /4 $ in (\ref{delta14}) for definiteness. The Nash-Moser  iteration would converge for any $ \d < 1 $,
see section \ref{ideas}. 

Another difference with respect to the scheme in \cite{BB07}, \cite{BB}, \cite{BBP}, is that we perform, 
at the same time, the Nash-Moser iteration and the multiscale argument 
for proving the invertibility of the linearized operators, see Lemma \ref{invLn+1}. This is more convenient
for proving measure estimates.

\subsection{Initialization of the Nash-Moser scheme}\label{sec:ini}

We perform the first step of the Nash-Moser iteration restricting
$ \l \in {\cal N}(\bar {\cal G}, 2 N_0^{-\s} )$ (the set $ \bar {\cal G} $ is defined in \eqref{diofs}).

\begin{lemma}\label{lem:ini} 
For all $  \l \in {\cal N}(\bar {\cal G}, 2N_0^{-\s}) $, 
the operator 
\be\label{calL0} 
{\cal L}_0 :=  P_0 (L_{\l \bar \om})_{| H_0}  
\ee 
(where $ L_{\om} $ is defined in (\ref{LomegaV})) is invertible and
\be\label{Lom}
\| {\cal L}_0^{-1} \|_{s_1} \leq  2N_0^{ \tau_1 + s_1} \g^{-1} \, . 
\ee
\end{lemma}

\begin{pf}
With the notations of \eqref{diofs}, for all $\l \in {\cal N} (\bar {\cal G},2N_0^{-\s})$, 
\be  \label{eigenv}
\forall |(l,j)|\leq N_0 \, , \ \   |\pm \l \bar{\om} \cdot l  + \mu_j| \geq \g N_0^{-\tau_1} - 
2 |\bar{\om}| N_0^{1-\s}  \geq \frac{\gamma}{2} N_0^{-\tau_1},
\ee
provided $ N_0 \geq 4\g^{-1}|\bar{\om}| $ (recall \eqref{def:sigma}, \eqref{tautau0} and $ \t_1 := d + \nu $). Then 
$ \| {\cal L}_0^{-1} \|_0 \leq 2 \g^{-1} N_0^{\t_1} $ and 
(\ref{Lom}) follows
by the smoothing property \eqref{S1}. 
\end{pf}

A fixed point of
\be\label{eq:F0}
F_0 : H_0 \to H_0 \, , \quad F_0 (u) :=  \e {\cal L}_0^{-1} P_0 ( f(u)+ g) \, ,  
\ee
is a solution of equation ($P_0$).  

\begin{lemma}\label{steppino} For $\e \g^{-1} N_0^{\tau_1 +s_1+\s} \leq c(s_1) $ small, 
 $ \forall \l \in {\cal N}(\bar {\cal G}, 2N_0^{-\s}) $,
the map $ F_0 $ is a contraction in 
$ {\mathtt B}_0(s_1) := \{ u \in H_0 \, : \, \| u \|_{s_1} \leq \rho_0 := N_0^{-\s} \} $.
\end{lemma}

\begin{pf}
The map $ F_0 $ maps $ {\mathtt B}_0(s_1) $ into itself, because, $ \forall \| u \|_{s_1} \leq \rho_0 $,  
$$
\| F_0(u) \|_{s_1} \stackrel{(\ref{Lom})} \leq 2 \e \g^{-1} N_0^{\tau_1+s_1} (\| f(u) \|_{s_1} + \|g\|_{s_1})
\stackrel{(F2), (\ref{gk})} 
\leq \e \g^{-1} N_0^{\tau_1 +s_1} C(s_1) 
\leq \rho_0 
$$
for $\e \g^{-1} N_0^{\tau_1 +s_1+\s}  $ is small enough. 
Moreover, $ \forall \| u \|_{s_1} \leq \rho_0 $,  
\be\label{DF0}
\| (D F_0)(u) \|_{s_1} = \e \| {\cal L}_0^{-1} P_0 (Df)(u)_{| H_0} \|_{s_1}  \stackrel{(\ref{Lom}), (F2)} 
\leq  \e  N_0^{\tau_1+s_1} \g^{-1}  C(s_1) 
\leq 1/2\, ,
\ee
implying that the map $ F_0 $ is a contraction in  $ {\mathtt B}_0(s_1) $. 
\end{pf}

Let $ {\wtilde u}_0 (\e, \l ) $ denote the unique solution of ($P_0$) in  $ {\mathtt B}_0 (s_1) $
defined for all $ (\e,\l ) \in  [0,\e_0] \times {\cal N}(\bar {\cal G}, 2N_0^{-\s}) $. For $ \e =  0 $ 
the map $ F_0 $ in (\ref{eq:F0}) has $ u = 0 $ as a fixed point. 
By uniqueness we deduce $ {\wtilde u}_0 (0,\l) = 0 $. 
Since the contracting map $F_0$ leaves $B_0(s_1) \cap {\cal U} $ invariant (see \eqref{subvs}), 
we deduce that $ {\wtilde u}_0 (\e, \l ) \in {\cal U} $. 
Moreover, by (\ref{DF0}),  the operator
\be\label{L0ep}
{\cal L}_0(\e) := P_0 \Big( L_\om - \e (Df)({\wtilde u}_0) \Big)_{|H_0} = {\cal L}_0 - \e P_0 (Df)({\wtilde u}_0)_{|H_0}
=  {\cal L}_0 \Big( I -  (DF_0)( {\wtilde u}_0 ) \Big) 
\ee 
is invertible and 
\be\label{pertu0}
\| {\cal L}_0^{-1}(\e) \|_{s_1} \leq 
2 \| {\cal L}_0^{-1} \|_{s_1} \stackrel{(\ref{Lom})} \leq  4 N_0^{\tau_1+s_1} \g^{-1}  \, .
\ee
The implicit function theorem implies that
$ {\wtilde u}_0 \in C^1(  [0,\e_0] \times {\cal N}(\bar {\cal G}, 2N_0^{-\s}); H_0) $ and
\be\label{derivate0}
\partial_{\e} {\wtilde u}_0 = {\cal L}_0^{-1}(\e) P_0(f({\wtilde u}_0) + g) \, , \quad
\partial_{\l} {\wtilde u}_0 = -  {\cal L}_0^{-1}(\e)  (\partial_\l {\cal L}_0) {\wtilde u}_0  \, . 
\ee
Then, by \eqref{derivate0}, \eqref{pertu0} and $ \partial_\l L_\om = {\rm diag}(\pm \ii \bar \om \cdot \partial_\vphi)  $,
 we get
\be\label{bound0}
\| \partial_{\e} {\wtilde u}_0\|_{s_1} \leq    N_0^{\tau_1+s_1}  \g^{-1} C(s_1) \, , \ \ 
\| \partial_{\l} {\wtilde u}_0\|_{s_1} \leq 4 |\bar{\om}|   N_0^{\tau_1+s_1}  \g^{-1}  \|{\wtilde u}_0 \|_{s_1+1} \leq  
C N_0^{\tau_1+ s_1+ 1- \s}\g^{-1} 
\ee
using that $\|{\wtilde u}_0 \|_{s_1+1} \leq N_0 \|{\wtilde u}_0 \|_{s_1} \leq N_0 N_0^{-\s} $.

Finally we define the $ C^1 $ map $ u_0 := \psi_0 {\wtilde u}_0 : [0,\e_0] \times \Lambda \to H_0 $ with
cut-off function $ \psi_0 :  \Lambda \to [0,1] $,  
\be\label{psi0}
\psi_0 :=
\begin{cases} 
1 \quad   {\rm if}   \  \l   \in  {\cal N}(\bar {\cal G}, N_0^{- \s }) \\
0  \quad  {\rm if} \  \l   \notin   {\cal N}(\bar {\cal G}, 2N_0^{- \s }) 
\end{cases} 
\quad  {\rm and} \qquad 
|D_\l \psi_0 | \leq N_0^\s C \, . 
\ee
Then \eqref{psi0}, $ \| {\wtilde u}_0 \|_{s_1} \leq N_0^{-\s} $ and  \eqref{bound0} imply 
(we have $ \partial_\e \psi_0 \equiv 0 $)
\be\label{u0s1}
\| u_0 \|_{s_1} \leq N_0^{-\s} \, , \quad \| \partial_{(\e,\l)} u_0 \|_{s_1} \leq C(s_1) N_0^{\t_1+ s_1 + 1} \g^{-1} \, . 
\ee
The statement $ (S1)_0 $ is proved.  
Note that  $ (S2)_0 $, $ (S3)_0 $ are  empty. 
Finally, also property  $ (S4)_0 $ is proved because, by \eqref{psi0} the function $ u_0(\e,\l) $ 
solves the equation ($P_0$) for all 
$ (\e,\l ) \in {\cal N}({\cal C}_0, N_0^{-\s}) $, since $ {\cal C}_0 = [0, \e_0] \times {\bar {\cal G}} $. 
 
\smallskip 

For the next steps of the induction we need the following lemma which establishes a property which replaces $ (S3)_n $
for the first steps of the induction. 

\begin{lemma}\label{Inizioind}
There exists $ N_0 := N_0 (S, V) \in \N $ and  $ c(s_1) > 0 $ such that, if  
\be\label{e0N0small} 
\e_0 N_0^{\t' + \d s_1}  \leq  c(s_1) \, , 
\ee
then $ \forall N_0^{1/C_2} \leq N \leq N_0 $, $ \forall \| u \|_{s_1} \leq 1 $,  $ {\cal G}_{N}(u) = [0,\e_0] \times \Lambda $. 
\end{lemma}

In order to prove Lemma \ref{Inizioind}  we prefix the following Lemma.

\begin{lemma}\label{alta}
For $ N \geq {\tilde N} (S,V) $ large enough, if 
\be\label{simplisigma}
\Big\| \Big( \vartheta \, {\rm I} + \Pi_{N,j_0} (-\D+ V(x))_{| E_{N,j_0}} \Big)^{-1} \Big\|_{L^2_x}  \leq N^{\t} \, , \ \  \vartheta \in \R \, ,
\ee
(see the definition of  $ E_{N,j_0} $ in \eqref{E0}) then, $ \forall s \in [s_0,S] $,  
\be\label{Goodstatic}
\Big|\!\!\Big| \Big( \vartheta {\rm I} +  \Pi_{N,j_0} (-\D+ V(x))_{|E_{N,j_0}} \Big)^{-1} 
\Big|\!\!\Big|_{s}  \leq \frac12 N^{\t' + \d s} \, .  
\ee
\end{lemma}

\begin{pf}
We apply  a simplified version of Proposition \ref{propinv} to
$  \vartheta {\rm I}+  \Pi_{N,j_0} (-\D+ V(x))_{| E_{N,j_0}} $. 
We sketch the main modifications only.  The scale $ N' $ in Proposition \ref{propinv} is here 
replaced by $ N $. 
Assumption (H1) follows from the regularity of the potential $ V(x) $ (see Lemma \ref{lem:multi}) and 
(H2) is \eqref{simplisigma}. With respect to 
 Proposition \ref{propinv}, we use a stronger version of assumption (H3), calling ``good sites'' the regular sites only, namely 
 the $ j \in \Z^d $, $ |j - j_0| \leq N $, such that 
$$
 |d_j| \geq \Theta \qquad  {\rm where}  \qquad d_j :=  \vartheta + \| j \|^2 + m 
$$
and $ m $ denotes the average of the potential $ V(x) $, see \eqref{average}.  
This is enough because here 
the  singular sites satisfy  separation properties.
For $ \Theta^{-1} \| V  \|_{s_1}  $ small enough 
we have the analogue of Lemma \ref{defmatrMN} (the proof is simpler because all 
the good   sites satisfy $ |d_j| \geq \Theta $). The separation properties of the singular sites 
$ j \in \Z^d $, $ |j - j_0| \leq N $, such that  $ | d_j | < \Theta $,
is proved as in section \ref{sec:sepa}: a $ M $-chain of  
singular sites has length at most $ L \leq M^{C_3(d)} $, see Lemma \ref{thm:separation}
and \eqref{Nspace}. 
Then, taking $ M := N^{\d / 2(1+ C_3(d))}$ we get a partition of the singular sites  in clusters 
$ \Om_\a $ satisfying 
$$ 
{\rm d}(\Om_\a, \Om_\b) > N^{\d / 2(1+ C_3(d))}  \quad  {\rm and} \quad 
{\rm diam}({\Omega}_\a) \leq M L \leq M^{1+ C_3(d)} =  N^{\d / 2}  \,.
$$  
Estimate \eqref{Goodstatic} follows by the arguments of Lemmas
\ref{defAprime}, \ref{defY} in section \ref{multiscale}.
\end{pf}

\begin{pfn}{\sc of Lemma}   \ref{Inizioind}. 
We claim that, $ \forall (\e,\l) \in [0, \e_0] \times \Lambda $, $ \forall j_0 \in \Z^d $,
\be\label{smallBN0}
B_{N}(j_0; \e,\l ) \subset \bigcup_{|(l,j-j_0)| \leq N} 
\Big\{  \theta \in \R \, : \, | \d_{l,j}^\pm (\teta)| \leq  N^{-\t} \Big\} 
\ee
where 
$$ 
\d_{l,j}^\pm (\teta) := \pm (\om \cdot l + \teta) + {\tilde \mu}_j \, , \ \om = \l \bar \om \, , \
 {\tilde \mu}_j := \, {\rm eigenvalues \ of }  \ \Pi_{N,j_0} (- \Delta + V(x) )_{|E_{N,j_0}} 
$$ 
(which depend on $ N $)  
and the subspace $ E_{N,j_0}$ is defined in \eqref{E0}.
Actually (\ref{smallBN0}) is equivalent to
\be\label{N0good}
|\d_{l,j}^\pm (\teta)|  >  N^{-\t} \, , \ \forall \, |(l,j-j_0)| \leq N \quad
\Longrightarrow \quad  {A}_{N,j_0}(\e, \l, \teta) \ {\rm is} \ N-{\rm good} 
\ee
with $ A = {\cal L}(u) = L_\om + \theta Y - \e (Df)(u) $.
We first prove that 
the left hand side condition in \eqref{N0good} implies
\be\label{stiLom}
Q_{N,j_0} := P_{N,j_0} (L_\om + \theta Y)_{| H_{N,j_0}} \quad {\rm satisfies} \quad
\nors{  Q_{N,j_0}^{-1} } \leq \frac12 N^{\t' + \d s} \, , \ \forall s \in [s_0,S] \, ,
\ee
(the subspace $ H_{N,j_0} $ is defined in \eqref{HNj0}).
Indeed, the operator $ L_\om $ is diagonal in time Fourier basis.
The left hand side condition in \eqref{N0good} is equivalent to
$$
\Big\| \Big( \pm (\l \bar \om \cdot l + \theta) {\rm I} + \Pi_{N,j_0} (-\D+ V(x))_{| E_{N,j_0}} \Big)^{-1} \Big\|_{L^2_x}  < N^{\t} \,   , \ \, \forall |l|  \leq N \, .
$$
Lemma \ref{alta} implies, for $ N \geq N_0^{1/C_2} \geq {\tilde N}(V, S) $,  that
$$
\Big|\!\!\Big|  \Big( \pm (\l \bar \om \cdot l + \theta){\rm I} + 
\Pi_{N,j_0} (- \Delta + V(x))_{|E_{N,j_0}}\Big)^{-1} 
\Big|\!\!\Big|_s  \leq \frac12 N^{\t' + \d s}   \, , \ \, \forall |l|  \leq N \, ,
$$
and \eqref{stiLom} follows because  $ Q_{N,j_0} $ is diagonal in time Fourier basis.

We now prove \eqref{N0good} by a perturbative argument. 
By (\ref{decayTu}) and $ \| u \|_{s_1} \leq 1 $ we have $ 
\norma (Df)(u) \norma_{s_1} \leq C(s_1) $. Hence
\be \label{QDf}
\e \norsone{ Q_{N,j_0}} \norsone{ (Df)(u)} \stackrel{\eqref{stiLom}} \leq \e N^{\t' + \d s_1 } C(s_1) 
 \leq \e_0 N_0^{\t' + \d s_1}  C(s_1) \stackrel{\eqref{e0N0small}} \leq 1/2 \, .
\ee
Then, by Lemma \ref{leftinv},  the matrix $  A_{N,j_0}(\e,\l,\theta ) =
P_{N,j_0} (L_\om +\theta Y - \e (Df)(u))_{| H_{N,j_0}} $ is invertible and
\be\label{ANini}
\forall s \in [s_0,s_1] \,   ,  \  \  \norma A_{N,j_0}^{-1}(\e,\l,\theta )  \norma_{s} \stackrel{\eqref{inv1}}\leq 2 \norma Q_{N,j_0}^{-1} \norma_{s} 
\stackrel{\eqref{stiLom}} \leq N^{\t'+ \d s} \, , 
\ee
namely it is $ N $-good. 

Finally, by (\ref{smallBN0}),   $ B_{N}(j_0; \e, \l )$ 
is included in an union of $ 2(2N+1)^b $ intervals of measure $\leq 2N^{-\tau}$, hence of 
$ 4 (2N+1)^b \leq N^{2d+\nu+4}$ intervals $ I_q $ of measure $ |I_q| \leq N^{-\tau}$. This proves that
any $(\e, \l) \in [0,\e_0]  \times \Lambda $ is $ N$-good  (see 
Definition \ref{def:freqgood}) for $ A = {\cal L}(u)$, 
namely that $ (\e, \l) $ is in $ {\cal G}_N (u) $, see \eqref{good}.
\end{pfn}

Finally we prove  $ (S5)_0 $. 
With estimates similar to the proof of $ (S1)_0 $ using the smallness condition  on $ \e_0 $ in 
\eqref{smallsto}, we deduce  $ (S5)_0 $-(i).
In order to estimate $ \partial_{(\e, \l)}  u_0 $, we use that 
the inverse of the operator $ {\cal L}_0(\e) = {\cal L}_0 - \e P_0 Df({\wtilde u}_0)_{| H_0} $ 
defined  in \eqref{L0ep} ($ {\cal L}_0 $ is defined in \eqref{calL0}) satisfies,
for $\l \in {\cal N} (\ov{\cal G} , 2 N_0^{-\s})$, 
\be\label{vaipi}
\norma {\cal L}_0^{-1}(\e) \norma_{s} \leq N_0^{\t'+ \d s} \, , \quad \forall s \in [s_1, S] \, .
\ee
Indeed,  note that by (\ref{eigenv}), for $N=N_0$ and $\theta=0$, the real numbers 
$|\d^{\pm}_{l,j} (0)|$  defined after (\ref{smallBN0})  are bounded from below by $\g N_0^{-\tau_1}/2 \geq N_0^{-\tau}$. Hence
$ {\cal L}_0 = Q_{N_0,0} $ satisfies \eqref{stiLom}, and
 Lemma \ref{leftinv} implies, $ \forall s \in [s_1,S] $,
\begin{eqnarray*}
\norma {\cal L}_0^{-1}(\e) \norma_{s} & \stackrel{\eqref{inv12}, \eqref{stiLom}} \leq & 
\Big(  1+C(s) \e \norso{Q_{N_0,0}^{-1}} \norso{(Df)(\wtilde{u}_0)}\Big) \frac{N_0^{\tau'+\d s}}{2} 
+ C(s) \e (N_0^{\tau'+\d s_0})^2 \, \nors{(Df)(\wtilde{u}_0)} \\
&  \stackrel{\eqref{stiLom}, \eqref{decayTu}, (S5)_0} \leq &   \Big(1+C(s)\e N_0^{\tau'+\d s_0} \Big)  \frac12 N_0^{\tau'+\d s}
+ C(s) \e N_0^{2(\tau'+\d s_0) + 2(\tau'+\d s_1 +1)} \\
&\stackrel{\eqref{smallsto}, \eqref{Sgr}} \leq &  N_0^{\tau'+\d s} 
\end{eqnarray*} 
since $ 4 \tau'+ 4\d s_1 +2 <S $.  
The bound $ (S5)_0 $-(ii) follows easily from (\ref{vaipi}). 
Let us give the details for  $ \partial_\e u_0 $ (which is not small with $ \e $).
We have
\begin{eqnarray}
\| \partial_{\e} {\wtilde u}_0 \|_S & \stackrel{\eqref{derivate0}} = &  
\| {\cal L}_0^{-1}(\e) P_0(f({\wtilde u}_0) + g) \|_S \nonumber \\
& \stackrel{\eqref{opernorm}} \leq & \norma {\cal L}_0^{-1}(\e)  \norma_{s_1} \| f({\wtilde u}_0) + g \|_S + 
C(S)  \norma {\cal L}_0^{-1}(\e)  \norma_S \| f({\wtilde u}_0) + g\|_{s_1} \nonumber \\
& \stackrel{\eqref{vaipi}, (F2), \eqref{gk}} \leq & 
C(S) N_0^{\t'+ \d s_1} (\| {\wtilde u}_0\|_S + 1) + C'(S) N_0^{\t'+ \d S} \nonumber \\
& \stackrel{(S5)_0-(i)} \leq &
C'(S) N_0^{3(\t'+ \d s_1)+2}  + C'(S) N_0^{\t'+ \d S} 
 \leq  N_0^{4 \t' + 2 s_1 + 4} \nonumber
\end{eqnarray}
by  \eqref{Sgr} and $ \d = 1 /4 $. Then  $(S5)_0$-(ii) is proved.

\subsection{Iteration of the Nash-Moser scheme}

Suppose, by induction, that we have already defined 
$ u_n \in C^1([0,\e_0] \times \Lambda; H_n \cap {\cal U}) $ 
and that properties $(S1)_k$-$(S5)_k $ hold for all $k\leq n$. 
We are going to define $u_{n+1}$ and prove the statements $(S1)_{n+1}$-$(S5)_{n+1}$.
Consider the operators $ {\cal L}(u) $ (introduced in (\ref{Linve})), 
\be\label{calnin}
{\cal L}(u) := {\cal L}(\om, \e , u)  := L_\om - \e (Df)(u) \, . 
\ee
In order to carry  out a modified Nash-Moser scheme, we shall study  the invertibility
of 
\be\label{caln+1}
{\cal L}_{n+1}(u_{n}) := P_{n+1} {\cal L}(u_{n})_{| H_{n+1}} 
\ee
and the tame estimates of its inverse,  
applying  Proposition \ref{propinv}. We  distinguish two cases. 
\\
If $ 2^{n+1} > C_2 $ (the constant $ C_2 $ is fixed in (\ref{tautau0})), 
then there exists a unique $ p \in [0,n] $ such that
\be\label{n+1np}
N_{n+1}  = N_p^{\chi}  \, , \quad  \chi = 2^{n+1-p} \in  [C_2, 2 C_2) \, .  
\ee
If $ 2^{n+1} \leq C_2 $ then there exists  $ \chi  \in [C_2, 2 C_2] $ such that 
\be\label{Nbasso}
N_{n+1} = {\bar N}^\chi \, , \ \  {\bar N} := [ N_{n+1}^{1/C_2}] \in (N_0^{1/\chi}, N_0) \, .  
\ee
If \eqref{n+1np} holds 
we  consider  in Proposition \ref{propinv} the two scales $ N' = N_{n+1} $, 
$ N = N_p $, see (\ref{newscale}). If \eqref{Nbasso}  holds, we set $N' = N_{n+1} $, $ N = \bar N $.

\smallskip

A key point of the whole induction process 
is that the separation properties of the bad sites of $ {\cal L} (u_n) + \teta Y $
hold uniformly for {\it all} $ \teta \in \R $ and  $ j_0 \in \Z^d $.

\begin{lemma}\label{H1H3} 
For all 
$$ 
(\e,\l) \in \bigcap_{k =  1}^{n+1} {\cal G}_{N_k}^0 (u_{k-1}) \, , \
 \teta \in \R \, , \ j_0 \in \Z^d \, ,
$$ 
the hypothesis  (H3) of Proposition \ref{propinv} apply
to $ A_{N_{n+1},j_0}(\e,\l,\teta ) $ where $ A(\e,\l,\teta) := {\cal L} (u_n) + \teta Y $. \end{lemma}

\begin{pf}
We give the proof when \eqref{n+1np} holds.  By remark \ref{good2}, a site 
\be\label{defE}
k \in E := \Big((0,j_0) + [-N_{n+1}, N_{n+1}]^b\Big) \times \{0,1\} \, ,  
\ee
which is $ N_p $-good for 
$ A(\e, \l, \teta ) := {\cal L}(u_n) + \teta Y  $ 
(see Definition \ref{GBsite} with $ A = A(\e,\l,\teta) $) is  also 
$$
(A_{N_{n+1},j_0}(\e, \l, \teta), N_p)-{\rm good}
$$ 
(see Definition \ref{ANreg} with 
$ A = A_{N_{n+1},j_0}(\e, \l ,\teta) $). As a consequence the
\be\label{badincl}
\Big\{ \ (A_{N_{n+1},j_0}(\e, \l, \teta), N_p){\rm -bad \ sites} \ \Big\}  \
\subset \  \Big\{ N_p{\rm -bad \ sites \ of} \  A(\e, \l ,\teta) \ {\rm with \ } |l| \leq N_{n+1} \Big\}.
\ee
and (H3) is proved if the latter $ N_p $-bad sites (in the right hand side of (\ref{badincl})) 
are contained in a disjoint union $ \cup_\a \Om_\a $ of clusters satisfying (\ref{sepabad}) (with $ N = N_p $). 
This is a consequence of  Proposition \ref{prop:separation} applied to 
the infinite dimensional matrix  $ A(\e, \l ,\teta) $. 
We claim that 
\be\label{inclup}
\bigcap_{k =  1}^{n+1} {\cal G}_{N_k}^0 (u_{k-1}) \subset {\cal G}_{N_p}(u_n) \, , \ {\rm i.e.} \
{\rm any} \ \, (\e,\l) \, \in \bigcap_{k =  1}^{n+1} {\cal G}_{N_k}^0 (u_{k-1}) \ \, {\rm is} \ \, N_p-{\rm good \ for} \  \, A(\e,\l, \teta) \, ,  
\ee
and then assumption (i) of Proposition \ref{prop:separation} holds. Indeed, if $p=0$ then \eqref{inclup} 
is trivially true  because $ {\cal G}_{N_0}(u_n) = [0,\e_0] \times \Lambda $, 
by Lemma \ref{Inizioind} and $(S1)_n $. If $ p \geq 1 $, we have
\be\label{claim1}
\| u_n - u_{p-1} \|_{s_1} \leq \sum_{k=p}^{n} \| u_k - u_{k-1}  \|_{s_1} 
\stackrel{(S2)_k} \leq  \sum_{k=p}^{n} N_k^{-\s-1} \leq
N_p^{-\s}  \sum_{k\geq p} N_k^{-1} \leq N_p^{-\s} 
\ee
and so $ (S3)_p $ implies 
\be\label{claim2}
\bigcap_{k =  1}^p {\cal G}_{N_k}^0 (u_{k-1}) \subset {\cal G}_{N_p}(u_n) \, .
\ee
Assumption (ii) of Proposition \ref{prop:separation} holds by (\ref{tautau0}), since $ \chi \in [C_2, 2 C_2) $.

When \eqref{Nbasso} holds the proof is analogous using Lemma \ref{Inizioind} with $ N = \bar N $ and  $ (S1)_n $. 
\end{pf}

\begin{lemma}\label{S3n+1}
Property $ (S3)_{n+1} $ holds.
\end{lemma}

\begin{pf}
We want to prove that 
$$ 
\| u - u_n \|_{s_1} \leq N_{n+1}^{-\s} \ \ {\rm and} \ \ (\e,\l) \in \bigcap_{k =  1}^{n+1} {\cal G}_{N_k}^0 (u_{k-1})
\quad \Longrightarrow \quad 
(\e,\l) \in {\cal G}_{N_{n+1}} (u) \, .
$$ 
Since $(\e , \l) \in {\cal G}^0_{N_{n+1}} (u_n)$, 
by (\ref{BNcomponent2}) and Definition \ref{def:freqgood} 
it is sufficient to prove that 
$ \forall j_0 \in \Z^d $,  
$$
B_{N_{n+1}} (j_0;\e,\l)(u)  \subset B_{N_{n+1}}^0 (j_0;\e,\l)(u_n)  \, , 
$$
(we highlight the dependence of these sets on $ u $, $ u_n $) 
or, equivalently, by (\ref{tetabadweak}), (\ref{tetabad}), that 
\be\label{inclusion+1}
\| A_{N_{n+1},j_0}^{-1} ( \e, \l , \teta )(u_n) \|_0 \leq N_{n+1}^\t  \quad \Longrightarrow \quad 
A_{N_{n+1},j_0} ( \e, \l , \teta )(u) \ {\rm is} \ N_{n+1} - {\rm good} \, , 
\ee
where $ A(\e, \l , \teta)(u ) = {\cal L}(u) + \teta Y = L_\om  + \teta Y - \e (Df)(u) $.

We prove (\ref{inclusion+1}) applying Proposition \ref{propinv} to $ A := A_{N_{n+1},j_0} ( \e, \l , \teta )(u)$
with $ E $ defined in \eqref{defE},  $  N' = N_{n+1}$, $ N =  N_p $ (resp. $N=\bar N$) if  \eqref{n+1np}
(resp. \eqref{Nbasso}) is satisfied. 
Assumption (H1)  holds with 
\be\label{decadimentoA}
\Upsilon \stackrel{ \eqref{Lu}, \eqref{decayTu}} = C (1+ \| u_n \|_{s_1}+ \norma V \norma_{s_1}) 
\stackrel{ (S1)_n, \eqref{gk}} \leq C' (V) \, .  
\ee  
By Lemma \ref{H1H3}, for all $  \teta \in \R $, $ j_0 \in \Z^d $, 
the hypothesis  (H3) of Proposition \ref{propinv} holds for $ A_{N_{n+1},j_0}(\e, \l ,\teta )(u_n) $. 
Hence, by Proposition \ref{propinv}, for $s\in [s_0,s_1]$, if
$$
\| A_{N_{n+1},j_0}^{-1} ( \e, \l , \teta )(u_n) \|_0 \leq N_{n+1}^\t  
$$
(which is assumption (H2)) then
\be\label{risultAN}
\nors{A^{-1}_{N_{n+1},j_0}(\e, \l ,\teta )(u_n)} \leq 
\frac{1}{4} N_{n+1}^{\tau' } \Big( N_{n+1}^{\d s}+  \norma V \norma_s + \e \nors{ (Df)(u_n) } \Big) \, .
\ee
Finally, since $ \| u - u_n \|_{s_1} \leq N_{n+1}^{-\sigma} $ we have
$$
\norsone{A_{N_{n+1},j_0} ( \e, \l , \teta )(u_n)-A_{N_{n+1},j_0} ( \e, \l , \teta )(u)} \leq C\e \| u - u_n \|_{s_1} \leq N_{n+1}^{-\s} 
$$
and  \eqref{inclusion+1} follows by \eqref{risultAN} and 
a standard perturbative argument (see for instance \eqref{inv1} in  Lemma \ref{leftinv}
with any $s\in [s_0,s_1]$ instead of $s_0$). 
\end{pf}

In order to define $ u_{n+1} $, we write, for $ h \in H_{n+1} $, 
\begin{eqnarray} 
P_{n+1}  \Big(L_\om ( u_n + h ) - \e (f( u_n + h ) + g) \Big) &= & 
P_{n+1} \Big(L_\om  u_n - \e  (f( u_n ) + g)\Big) \nonumber \\
& + & 
P_{n+1} \Big(L_\om h - \e  (Df)( u_n ) h\Big) + R_n ( h ) \nonumber \\
& = & r_n + {\cal L}_{n+1} (u_n) h + R_n ( h )  \label{scritt1}
\end{eqnarray}
where 
$ {\cal L}_{n+1}(u_{n}) $
is defined in  \eqref{caln+1} and 
\be\label{Rnh}
r_n :=  P_{n+1} \Big(L_\om u_n - \e  (f( u_n ) + g)\Big) \, ,  \quad 
R_n ( h ) := - \e  P_{n+1} \Big( f( u_n + h ) - f( u_n ) - (Df)( u_n ) h \Big) \, .
\ee
By $ (S4)_n $, if $ (\e,\l) \in {\cal N}({\cal C}_n, N_{n}^{-\s}) $ then $ u_n $ solves the equation $(P_n) $ and so
\be\label{rnhigh}
r_n =  P_{n+1}P_n^\bot  \Big( L_\om  u_n - \e   (f( u_n ) + g)\Big)  = 
P_{n+1}P_n^\bot  \Big(V_0 \,  u_n - \e   (f( u_n ) + g)\Big) \, , 
\ee
using also that $ P_{n+1} P_n^\bot ( D_\om u_n) = 0 $, see \eqref{Lomega}.
Note that, by (\ref{defNn}) and $ \s \geq 2 $ (see \eqref{def:sigma}), for $ N_0  \geq 2 $, we have the inclusion
\be\label{inclu}
{\cal N}({\cal C}_{n+1}, 2 N_{n+1}^{-\s}) \subset {\cal N}({\cal C}_n, N_{n}^{-\s}) \, . 
\ee

\begin{lemma}\label{invLn+1}
{\bf (Invertibility of $ {\cal L}_{n+1} $)}  
For all $  (\e, \l ) \in  {\cal N}({\cal C}_{n+1}, 2 N_{n+1}^{-\s}) $
the operator  $ {\cal L}_{n+1}(u_n) $ is invertible and,  for $ s  = s_1, S $, 
\be\label{normabassa}
\norma {\cal L}_{n+1}^{-1}(u_n) \norma_s \leq  N_{n+1}^{\t' + \d s } \,  . 
\ee
As a consequence, by (\ref{opernorm}),  $ \forall h \in H_{n+1} $, 
\be\label{os1}
\| {\cal L}_{n+1}^{-1}(u_n) h \|_{s_1} \leq C(s_1) N_{n+1}^{\t' + \d s_1 } \| h \|_{s_1} \, , 
\ee
\be\label{os2}
\| {\cal L}_{n+1}^{-1}(u_n) h \|_S \leq  N_{n+1}^{\t' + \d s_1 } \|h\|_S + C(S) N_{n+1}^{\t' + \d S }  \|h \|_{s_1} \, .
\ee
\end{lemma}

\begin{pf} 
We give the proof when  \eqref{n+1np} holds. The other case is analogous. 
First assume $(\e,\l) \in {\cal C}_{n+1} $, see  \eqref{Gscavo}. Then since 
 $ (\e,\l) \in {\mathtt G}_{N_{n+1}}(u_n) $ (see (\ref{Binver}) with $ A_N(\e,\l) = {\cal L}_{n+1}(u_n) $),  
the operator  $ {\cal L}_{n+1}(u_n) $ is invertible and  
\be\label{Ln+1} \| {\cal L}_{n+1}^{-1}(u_n) \|_0 \leq N_{n+1}^\tau \, . 
\ee 
We  now apply the multiscale Proposition \ref{propinv} to $ A := {\cal L}_{n+1} (u_n) $
with
$$
E := [-N_{n+1}, N_{n+1}]^b  \times \{0,1\} \, , \quad N' = N_{n+1} \, , \quad N =  N_p, \ {\rm see} \  \eqref{n+1np} \, . 
$$
By remark \ref{defchi} and since $ \chi \in [C_2, 2C_2) $ (see (\ref{n+1np})) the assumptions (\ref{dtC})-(\ref{s1}) hold. 
Assumption (H1) holds with \eqref{decadimentoA}.
Assumption (H2) holds by (\ref{Ln+1}). 
Moreover, by the definition of ${\cal C}_{n+1}$, as a particular case of Lemma \ref{H1H3}  -for $ \teta = 0 $,  $ j_0 = 0 $-,
the hypothesis (H3) of Proposition \ref{propinv} holds for $ {\cal L}_{n+1} (u_n) $. 
Then Proposition \ref{propinv} applies and we get that, $\forall (\e,\l) \in {\cal C}_{n+1} $, $ \forall s \in \{ s_1,S \} $, 
$$
\nors{{\cal L}_{n+1}^{-1}(u_n)}  \stackrel{(\ref{A-1alta})} \leq 
\frac{1}{4} N_{n+1}^{\tau' } \Big( N_{n+1}^{\d s}+  \norma V \norma_s + \e \nors{ (Df)(u_n) } \Big) \, ,
$$
whence, for $ s = s_1 $, 
\be\label{fo1}
\norsone{{\cal L}_{n+1}^{-1}(u_n)} 
\stackrel{(\ref{decayTu}),  (S1)_n, \eqref{gk}}\leq \frac{1}{4} N_{n+1}^{\tau' } \Big( N_{n+1}^{\d s_1}+ \norsone{V}  + \e C(s_1) \Big)
\leq \frac12 N_{n+1}^{\t' + \d s_1 } 
\ee
and, for $ s = S $,  recalling that $ U_n := \| u_n \|_S $, 
\begin{eqnarray}\label{fo2}
\norS{{\cal L}_{n+1}^{-1}(u_n)} 
&\stackrel{(\ref{decayTu}), \eqref{gk}} \leq  & \frac14 N_{n+1}^{\tau' } \Big( N_{n+1}^{\d S}+ \norS{V} + \e C(S)(1+U_n) \Big) \nonumber \\
& \stackrel{(S5)_n} \leq &
\frac14 N_{n+1}^{\tau' } \Big( N_{n+1}^{\d S}+ C'(S) N_n^{2(\t' + \d s_1 + 1)}  \Big)
\leq \frac 12 N_{n+1}^{\t' + \d S }  
\end{eqnarray}
by   \eqref{Sgr} and $ \d = 1/4 $.
Assume next  $ (\e',\l') \in {\cal N}({\cal C}_{n+1}, 2 N_{n+1}^{-\s}) $ and let $ (\e,\l) \in {\cal C}_{n+1} $
be such that $ | (\e', \l') - (\e,\l)| < 2 N_{n+1}^{-\s} $. 
We write
$$
{\cal L}_{n+1}(u_n(\e', \l')) =  {\cal L}_{n+1}(u_n(\e, \l)) + {\mathtt R}_{n+1}  
$$
where  $ {\cal L}_{n+1}(u_n(\e, \l)) $  satisfies \eqref{fo1}-\eqref{fo2} and 
$$
{\mathtt R}_{n+1}  := {\cal L}_{n+1}(u_n(\e', \l')) -  {\cal L}_{n+1}(u_n(\e, \l))  \, . 
 $$
By \eqref{caln+1}, \eqref{decayTu}, (F2), \eqref{lions},  \eqref{smallsto}, $(S1)_n $, $(S5)_n $,  
\be\label{Rdiff}
\norsone{{\mathtt R}_{n+1}} \leq C(s_1) N_{n+1}^{-\s+1} \, , \quad 
\norS{{\mathtt R}_{n+1}} \leq C(S)  N_n^{4\tau'+2 s_1 +4} N_{n+1}^{-\s} \, . 
\ee
We apply Lemma \ref{leftinv} with 
$$ 
M = {\cal L}_{n+1}(u_n(\e, \l)) \, , \quad N = {\cal L}_{n+1}^{-1}(u_n(\e, \l)) \, , \quad
P = {\mathtt R}_{n+1} \, .
$$
By \eqref{fo1}, \eqref{Rdiff} and \eqref{def:sigma} the perturbative 
assumption \eqref{NR12} holds 
 with index $ s_1 $ instead of $s_0$.  
Then \eqref{inv1}, \eqref{inv12} (with indices $ s_1 , S $ instead of  $ s_0, s $) imply \eqref{normabassa}
for all $ (\e',\l') \in {\cal N}({\cal C}_{n+1}, 2 N_{n+1}^{-\s})  $,  
by  \eqref{fo1}, \eqref{fo2},  \eqref{Rdiff}, \eqref{def:sigma}.
\end{pf}

By \eqref{scritt1}, setting
\be\label{Fn+1}
F_{n+1}: H_{n+1} \to H_{n+1} \, ,  \qquad
F_{n+1}(h) := - {\cal L}_{n+1}^{-1} (u_n) ( r_n +  R_n ( h )) \, , 
\ee
the equation ($ P_{n+1} $) is equivalent to the fixed point problem $ h = F_{n+1}( h ) $.

\begin{lemma} \label{lemcon}
{\bf (Contraction in $ \| \ \|_{s_1} $-norm)}  
$ \forall  (\e,\l) \in {\cal N}({\cal C}_{n+1}, 2 N_{n+1}^{-\s}) $, 
$ F_{n+1} $  is a contraction in
\be\label{defrhon+1}
{\mathtt B}_{n+1}(s_1) := \Big\{ h \in  H_{n+1} \  : \  \|h\|_{s_1} \leq  \rho_{n+1} :=    N_{n+1}^{-\s-1} \Big\} \, . 
\ee
The unique fixed point $ {\wtilde h}_{n+1}(\e,\l) $ of $ F_{n+1} $ in  $ {\mathtt B}_{n+1}(s_1) $ 
belongs to ${\cal U}$  (see \eqref{subvs}) and
satisfies  
\be\label{hn+1Un}
\| {\wtilde h}_{n+1} \|_{s_1} \leq K(S) N_{n+1}^{\t'+ \d s_1} N_n^{-(S - s_1)} U_n \, . 
\ee
\end{lemma}

\begin{pf} 
For all $ (\e,\l) \in {\cal N}({\cal C}_{n+1}, 2 N_{n+1}^{-\s} ) $, 
by (\ref{Fn+1}) and (\ref{os1}), we have 
\be\label{30}
\| F_{n+1} ( h ) \|_{s_1} 
\leq  C(s_1) N_{n+1}^{\t^{'} + \d s_1} ( \| r_n \|_{s_1} + \| R_n ( h ) \|_{s_1})
\ee
and $ r_n $ has the form (\ref{rnhigh}) because of (\ref{inclu}). 
 Moreover
(recall that $ U_n := \| u_n \|_S $)
\begin{eqnarray}
\| r_n \|_{s_1} + \| R_n ( h ) \|_{s_1} &  \stackrel{(\ref{rnhigh}),(\ref{S2}), \eqref{Rnh}, (\ref{P3s})} \leq & 
  N_n^{- (S-s_1)} ( \| V_0 \,  u_n\|_S +  \e \|  f ( u_n ) \|_{S} + \e \| g \|_{S}) 
+  \e C(s_1) \|h \|_{s_1}^2  \nonumber \\ 
& \stackrel{(\ref{fDftame}), (\ref{gk})} \leq & C(S) N_{n}^{-(S - s_1)} (U_n+1) +  
\e \, C(s_1)  \|h \|_{s_1}^2   \label{rnprima}  \\
& \stackrel{(S5)_n} \leq &  C(S) N_n^{-(S - s_1)} N_n^{2(\t'+\d s_1 + 1)} +  \e \, C(s_1)  \|h \|_{s_1}^2  \, . \label{rnRn}
\end{eqnarray}
(\ref{30}) and (\ref{rnRn}) imply (using also (\ref{defNn})), for some $ K(S), K(s_1) > 0 $, 
\begin{eqnarray}\label{piccocontra}
\| h \|_{s_1} \leq \rho_{n+1}  \quad  \Longrightarrow \quad 
\| F_{n+1} ( h ) \|_{s_1} & \leq &   K(S)  N_{n+1}^{2(\t'+\d s_1) + 1} N_n^{-(S-s_1)} + 
\e K(s_1) N_{n+1}^{\t' + \d s_1}\rho_{n+1}^2 \nonumber \\
& \leq & \rho_{n+1} := N_{n+1}^{-\s -1} \,  ,  \nonumber 
\end{eqnarray}
because the choice of $ S $ in (\ref{Sgr}) 
and of $ \s $ in (\ref{def:sigma})  
imply (for $ N \geq N_0(S) $)
\be\label{piccorho}
 K(S)  N_{n+1}^{2(\t'+\d s_1) + 1} N_n^{-(S-s_1)} \leq \frac{\rho_{n+1}}{2} \, , 
\quad \e K(s_1) N_{n+1}^{\t' + \d s_1} \rho_{n+1} \leq \frac{1}{2} \, .
\ee
Next, differentiating (\ref{Fn+1}) with respect to $  h $ and using (\ref{Rnh}) we get
$$
D_h F_{n+1}(h)[v] = {\cal L}_{n+1}^{-1} (u_n) \e  P_{n+1} \Big(  (Df) ( u_n + h)[v] - (Df)(u_n)[v] \Big)
$$
and, for all $ \| h \|_{s_1} \leq \rho_{n+1} $,  
using \eqref{secondorder}  with $ s = s_1 $,   
$$
\| D_h F_{n+1}(h)[v] \|_{s_1} 
\stackrel{(\ref{os1})} \leq \e K(s_1) N_{n+1}^{\t' +\d s_1} \rho_{n+1} \| v \|_{s_1} 
\stackrel{(\ref{piccorho})} \leq \frac12 \| v \|_{s_1} \, .
$$
Hence $ F_{n+1} $ is a contraction in $ {\mathtt B}_{n+1}(s_1) $. Since $u_n \in {\cal U}$, it is easy to check that 
$F_{n+1}$ leaves $B_{n+1}(s_1) \cap {\cal U} $ invariant, hence $ \wtilde{h}_{n+1} \in {\cal U}$. 
Finally, (\ref{Fn+1}), (\ref{30}), (\ref{rnprima}) and (\ref{piccorho})  imply \eqref{hn+1Un}.
\end{pf}

Since $ {\wtilde h}_{n+1} (\e,\l)$ solves, for all $ (\e,\l) \in {\cal N}({\cal C}_{n+1}, 2N_{n+1}^{-\s}) $, the equation
\be\label{Un+1}
Q_{n+1} (\e , \l , h ) := P_{n+1} \Big(L_\om ( u_n + h )  - \e  (f( u_n + h ) + g) \Big) = 0   \, , \ \  h \in H_{n+1} \, , 
\ee
and $ u_n (0,\l) \stackrel{(S1)_n} = 0 $, we deduce, by the uniqueness of the fixed point, that
$$
{\wtilde h}_{n+1}(0,\l) = 0 \, , \quad \forall (0,\l) \in {\cal N}({\cal C}_{n+1}, 2N_{n+1}^{-\s}) \, .
$$

\begin{lemma} \label{htn+1} {\bf (Estimate in high norm)}
$ \forall  (\e,\l) \in {\cal N}({\cal C}_{n+1}, 2 N_{n+1}^{-\s}) $ we have 
\be\label{halta1} 
\|  {\wtilde h}_{n+1}  \|_{S} \leq K(S) N_{n+1}^{\t' + \d s_1} U_n \, . 
\ee
\end{lemma}

\begin{pf}
We have 
\begin{eqnarray}
\| {\wtilde h}_{n+1} \|_S & \stackrel{\eqref{Fn+1}}= &  
\Big\| {\cal L}_{n+1}^{-1}(u_n) ( r_n + R_n( {\wtilde h}_{n+1})) \Big\|_{S}  \label{disu} \\
& \stackrel{(\ref{os2})} \leq &  
N_{n+1}^{\t' + \d s_1 } \Big(\|  r_n \|_S + \| R_n( {\wtilde h}_{n+1}) \|_S \Big) + C(S) N_{n+1}^{\t' + \d S }  
\Big(\|  r_n \|_{s_1} + \| R_n ({\wtilde h}_{n+1}) \|_{s_1} \Big) \, . \nonumber
\end{eqnarray} 
Now, by \eqref{rnhigh}, $(S1)_n $, (F2), (F3),  
(\ref{gk}), \eqref{defk}, \eqref{Rnh}, and setting $ U_n := \| u_n \|_S $ (we can 
suppose $ U_n \geq 1 $) we get
\be\label{pezz1}
\| r_n \|_S +  \| R_n({\wtilde h}_{n+1}) \|_S  \leq C(S) (U_n  +  \e \rho_{n+1}  \| {\wtilde h}_{n+1} \|_S) 
\ee
and, using also \eqref{rnprima}, \eqref{hn+1Un} and the second inequality in \eqref{piccorho}, 
\be\label{pezz2}
\|  r_n \|_{s_1} + \| R_n ({\wtilde h}_{n+1}) \|_{s_1} \leq C(S)  N_n^{-(S-s_1)} U_n \, .
\ee
Then \eqref{disu}, \eqref{pezz1}, \eqref{pezz2} imply that
\begin{eqnarray}
\| {\wtilde h}_{n+1} \|_S & \leq &  
 C(S) \Big(N_{n+1}^{\t' +\d s_1}  + N_{n+1}^{\t'+ \d S} N_n^{-(S-s_1)}\Big) U_n + 
C(S) \e N_{n+1}^{\t'+\d s_1}  \rho_{n+1} \| {\wtilde h}_{n+1} \|_{S}  \label{quellaS'} \\
& \stackrel{\eqref{Sgr}, \eqref{defrhon+1}} \leq &   C'(S) N_{n+1}^{\t' +\d s_1} U_n + 
\e C(S)  N_{n+1}^{\t'+\d s_1 - \s -1 }  \| {\wtilde h}_{n+1} \|_{S} \nonumber  \\
& \stackrel{\eqref{def:sigma}} \leq & C'(S) N_{n+1}^{\t' +\d s_1} U_n + \frac12  \| {\wtilde h}_{n+1} \|_{S} 
\nonumber 
\end{eqnarray}
for $ \e_0 \leq \e_0(S) $ small.
As a consequence we get $ \| {\wtilde h}_{n+1}  \|_S \leq 2  C'(S) N_{n+1}^{\t'+\d s_1}  U_n $
and (\ref{halta1}) follows. 
\end{pf}

\begin{lemma}\label{deriva}
{\bf (Estimate of the derivatives)} 
The map $ {\wtilde h}_{n+1} \in  C^1({\cal N}({\cal C}_{n+1}, 2 N_{n+1}^{-\s}), H_{n+1}) $ and
\be\label{derivas1}
\| \partial_{(\e,\l)} {\wtilde h}_{n+1} \|_{s_1} \leq  N_{n+1}^{ -1}  \, , \quad 
 \| \partial_{(\e, \l)} {\wtilde h}_{n+1} \|_{S} \leq  N_{n+1}^{\t'+\d s_1 +1} 
 \Big(N_{n+1}^{\t'+\d s_1 +1} U_n + U_n' \Big) \, . 
\ee
\end{lemma}

\begin{pf}
For all $ (\e,\l ) \in {\cal N}({\cal C}_{n+1}, 2 N_{n+1}^{- \s}) $,
$ {\wtilde h}_{n+1}(\e, \l)  $  is a solution of $ Q_{n+1}( \e , \l , {\wtilde h}_{n+1}(\e , \l)) = 0  $, see (\ref{Un+1}).
We have, see \eqref{caln+1},
\be \label{DLQ} 
D_h Q_{n+1} (\e , \l , {\wtilde h}_{n+1}) =  {\cal L}_{n+1} ( u_n + {\wtilde h}_{n+1})
=  {\cal L}_{n+1} ( u_n) - \e P_{n+1} \Big((Df)(u_n+ {\wtilde h}_{n+1}) - (Df)(u_n) \Big)
\ee
which is invertible by Lemma \ref{leftinv} applied with 
$$ 
M \to  {\cal L}_{n+1} ( u_n) \, , 
P \to  - \e P_{n+1}((Df)(u_n+ {\wtilde h}_{n+1}) - (Df)(u_n)) \, , \ s_0 \to s_1 \, . 
$$ 
Indeed the hypothesis \eqref{NR12} follows  from  
(\ref{normabassa}) with $ s = s_1 $, (F1), 
$ (S1)_n $,  Lemma \ref{lem:multi},
$  \| {\wtilde h}_{n+1} \|_{s_1}\leq \rho_{n+1} $ and 
\eqref{piccorho}. 
Therefore Lemma \ref{leftinv} with $ s = s_1 $ implies  
\be\label{invnear}
\Big|\!\!\Big| {\cal L}_{n+1}^{-1}( u_n + {\wtilde h}_{n+1})  \Big|\!\!\Big|_{s_1} \stackrel{\eqref{inv1}} \leq
2 \norma {\cal L}_{n+1}^{-1} (u_n) \norma_{s_1} \stackrel{\eqref{normabassa}} \leq 2  N_{n+1}^{\t' + \d s_1 } \,
\ee
and, by \eqref{inv2}, (\ref{normabassa}) with $ s = S $,  \eqref{halta1},   $ (S5)_n $, \eqref{secondorder},  
$ \d = 1 / 4 $, \eqref{Sgr}, 
\be\label{invnear2}
\Big|\!\!\Big| {\cal L}_{n+1}^{-1}( u_n + {\wtilde h}_{n+1})  \Big|\!\!\Big|_{S} \leq C(S) N_{n+1}^{\t' + \d S}
\, . 
\ee
Hence, the Implicit function theorem implies $ {\wtilde h}_{n+1}  \in C^1 ({\cal N}({\cal C}_{n+1}, 2N_{n+1}^{-\s}), H_{n+1} )$
and  
\be\label{deom}
\partial_{(\e,\l)} {\wtilde h}_{n+1} 
\stackrel{\eqref{DLQ}} = - {\cal L}_{n+1}^{-1}(u_n + {\wtilde h}_{n+1})  
\Big(\partial_{(\e,\l)} Q_{n+1}\Big)(\e,\l, {\wtilde h}_{n+1}) \, . 
\ee
By $ (S4)_n $,  $ u_n (\e,\l) $ solves $(P_n)$ for $(\e,\l) \in {\cal N}({\cal C}_{n+1},  2N_{n+1}^{-\s}) 
\stackrel{\eqref{inclu}}\subset {\cal N}({\cal C}_{n},  N_n^{-\s}) $. Then 
\begin{eqnarray} \label{der1}
(\partial_{\e} Q_{n+1})(\e,\l, {\wtilde h}_{n+1}) & = &
P_{n+1} P_n^\bot ( V_0 \, \partial_\e u_n) + 
 P_{n} (f(u_n)+ g) - P_{n+1} (f( u_n + {\wtilde h}_{n+1}) + g)  \nonumber \\ 
& + &  \e P_n (Df) ( u_n) \partial_\e u_n  -\e  P_{n+1} (Df) (u_n + {\wtilde h}_{n+1}) \partial_\e u_n 
\end{eqnarray}
(we use also that $ P_{n+1} P_n^\bot ( D_\om u_n) = 0 $ since $ u_n \in H_n $, see \eqref{Lomega})
and
\begin{eqnarray} \label{der2}
(\partial_{\l} Q_{n+1})(\e,\l, {\wtilde h}_{n+1}) & = & P_{n+1} P_n^\bot ( V_0 \, \partial_\l u_n) + 
(\partial_\l L_\om) {\wtilde h}_{n+1} \\
& + &  \e P_n (Df) (u_n) \partial_\l u_n - \e P_{n+1} (Df)(u_{n} +{\wtilde h}_{n+1}) \partial_\l u_n \, . \nonumber
\end{eqnarray} 
We deduce from (\ref{invnear})-(\ref{der2}) the estimates (\ref{derivas1}) using also \eqref{opernorm}, 
(F1), (F2), (F3), $ (S1)_n $, (\ref{S2}), 
$ (S5)_n $, \eqref{gk}, (\ref{Sgr}), (\ref{hn+1Un}), \eqref{halta1}. We omit the details. 
\end{pf}

We now define
a $ C^1 $-extension of $ ({\wtilde h}_{n+1})_{|{\cal C}_{n+1}}  $  
onto the whole $ [0,\e_0] \times \Lambda $. 

\begin{lemma} {\bf (Extension)}\label{whit}
There is $ h_{n+1} \in C^1([0,\e_0)\times \Lambda , H_{n+1} \cap {\cal U}) $  satisfying $ h_{n+1}(0,\l) = 0 $, 
\be\label{stimader2} 
\| h_{n+1}\|_{s_1} 	\leq N_{n+1}^{-\s-1} \, , \quad 
\| \partial_{(\e,\l)} h_{n+1}\|_{s_1}	 \leq   N_{n+1}^{-1/2}   
\ee
and $ h_{n+1} $ is equal to $ {\wtilde h}_{n+1} $ on ${\cal N}({\cal C}_{n+1}, N_{n+1}^{-\s}  ) $. 
\end{lemma}

\begin{pf} 
Let
\be\label{defih}
h_{n+1} (\e,\l) := \left\{ \begin{array}{lll} 
\psi_{n+1}(\e,\l) {\wtilde h}_{n+1}(\e,\l) & {\rm if} &  (\e,\l) \in {\cal N}({\cal C}_{n+1}, 2 N_{n+1}^{- \s } )   \\
0 & {\rm if} & (\e,\l) \notin  {\cal N}({\cal C}_{n+1},2 N_{n+1}^{-\s} )  \end{array}
\right.
\ee
where $ \psi_{n+1} $ is a $C^\infty $ cut-off function satisfying
$$ 
0 \leq \psi_{n+1} \leq 1 \, , \quad
\psi_{n+1} \equiv 
\begin{cases} 
1 \quad {\rm if} \  (\e,\l) \in {\cal N}({\cal C}_{n+1},  N_{n+1}^{-\s} ) \\
0  \quad  {\rm if} \  (\e,\l) \notin  {\cal N} ({\cal C}_{n+1}, 2{N_{n+1}^{- \s})} 
\end{cases}  {\rm and} \quad 
|\partial_{(\e,\l)} \psi_{n+1} | \leq N_{n+1}^\s C \, . 
$$
Then $ \| h_{n+1}\|_{s_1} \leq \| {\wtilde h}_{n+1} \|_{s_1} \leq  N_{n+1}^{-\s-1} $ by 
Lemma \ref{lemcon}, and, 
$$
\| \partial_{(\e,\l)} h_{n+1}\|_{s_1} \leq | \partial_{(\e,\l)} \psi_{n+1}| \, \| {\wtilde h}_{n+1} \|_{s_1} + 
\|  \partial_{(\e,\l)} {\wtilde h}_{n+1} \|_{s_1} 
\leq  N_{n+1}^{-1/2}  
$$
thanks to the first estimate  in (\ref{derivas1}), and for $ N_0 $ large.
\end{pf}

Finally we define $ u_{n+1} \in C^1([0,\e_0) \times \Lambda,H_{n+1}\cap {\cal U}) $ as 
\be\label{un+1tilde}
u_{n+1} :=  u_n + h_{n+1} \, . 
\ee
By Lemma \ref{whit}, on $ {\cal N}({\cal C}_{n+1}, N_{n+1}^{-\s } ) $ we have $  h_{n+1} = {\wtilde h}_{n+1} $ 
that solves equation (\ref{Un+1}) and so $ u_{n+1} $ solves equation $ (P_{n+1}) $. Hence
$ (S4)_{n+1} $ holds. By Lemma \ref{whit}, property $ (S2)_{n+1} $ holds. Property $ (S1)_{n+1} $
follows as well because
$$
\| u_{n+1} \|_{s_1} \leq \| u_0 \|_{s_1} + \sum_{k=1}^{n+1} \| h_{k} \|_{s_1} 
\stackrel{(\ref{u0s1}), (S2)_{n+1}}\leq \frac12 + \sum_{k=1}^{n+1} N_{k}^{-\s-1} \leq \frac12 + N_1^{-1} \leq 1 
$$ 
and the estimate 
$ \| \partial_{(\e,\l)} u_{n+1} \|_{s_1} \leq C(s_1) N_0^{\t_1+s_1+1} \g^{-1} $ follows in the same way. 

\begin{lemma}\label{S5n+1}
Property $ (S5)_{n+1} $ holds.
\end{lemma}

\begin{pf}
By the definition of $ U_n $, and since $ \| h_{n+1} \|_S \leq \| {\wtilde h}_{n+1}\|_S $, 
we get 
$$
U_{n+1} \leq U_n + \| {\wtilde h}_{n+1} \|_{S} \stackrel{(\ref{halta1})} \leq K'(S) N_{n+1}^{\t' +   \d s_1}  U_n
\stackrel{(S5)_n} \leq 
K'(S) N_{n+1}^{\t' +   \d s_1} N_n^{2(\t'+\d s_1+1)} \stackrel{(\ref{defNn})} \leq N_{n+1}^{2(\t' +   \d s_1 + 1)} \, . 
$$
The estimate for $ U_{n+1}' $ follows similarly by (\ref{halta1}), (\ref{derivas1}), $(S5)_n$.
\end{pf}

\subsection{Proof of Theorem \ref{thm:main}}\label{sec:measfina}

By Theorem \ref{cor1} it remains to prove that the measure estimate  (\ref{Cmeas}) holds.

\begin{lemma}\label{NRfre}
The set  $ {\cal G} $ defined in \eqref{diofs} satisfies
\be\label{calGg}
| \bar {\cal G} | = 1 - O(\g) \, . 
\ee
\end{lemma}

\begin{pf}
The  $ \l $ such that  (\ref{diofs}) is violated are 
\be\label{eq:1}
{\bar {\cal G}}^c \cap [1/2,3/2] \subseteq 
\bigcup_{|l| \leq N_0, |j| \leq N_0} {\cal R}_{l,j} \quad 
{\rm where} \quad 
{\cal R}_{l,j}^\pm := \Big\{ \l \in [1/2, 3/2 ] \, : \, |\pm \lambda \bar \om \cdot l + \mu_j | < \frac{\g}{N_0^{\t_1}}  \Big\} \, . 
\ee
Dividing by $ \l $,  we have to estimate the $ \xi := 1 / \l  \in [2/3,2] $ such that 
$$
| \pm \bar \om \cdot l + \xi \mu_j | < C \frac{\g}{N_0^{\t_1}} \, . 
$$
The derivative of
the functions $ g_{lj}^\pm(\xi) := \pm \bar \om \cdot l + \xi \mu_j $ satisfies
$ |\partial_\xi g_{lj}^\pm(\xi)|  =  |\mu_j | \geq \b_0 > 0 $, because $ \Pi_0 (- \D + V(x) )_{| E_0} \geq \b_0 I $
 by \eqref{eq:posi}. 
As a consequence,  
we estimate
\be\label{eq:3}
|{\cal R}_{l,j}^\pm| \leq  \frac{C}{ \b_0} \frac{\g}{N_0^{\tau_1} } \, . 
\ee
Then (\ref{eq:1}), (\ref{eq:3}), imply
$$
| {\bar {\cal G}}^c \cap [1/2,3/2] | \leq \sum_{|l| \leq N_0, |j| \leq  N_0, \pm} |{\cal R}^\pm_{l,j}| \leq
C \frac{\g}{\b_0} \frac{N_0^{d+\nu}}{N_0^{\tau_1}} = O(\g)
$$
since $ \t_1  \geq  d+ \nu  $.
\end{pf}

Finally we choose
\be\label{N0ge}
\g := \e_0^{\a} \quad {\rm with} \quad \a := 1 / (S + 1) \, , \quad N_0 := 4 \g^{-1} \, ,
\ee
so that (\ref{smallsto}) is fulfilled for $ \e_0 $ small enough.
The complementary set of $ {\cal C}_\infty $ in $ [0,\e_0] \times \Lambda $ has measure
\begin{eqnarray}
|{\cal C}_\infty^c| & \stackrel{ \eqref{Cinfty}, \eqref{Gscavo}}= & 
\Big| \bigcup_{k \geq 1} {\mathtt G}_{N_k}^c (u_{k-1}) 
\bigcup_{k \geq  1 }  ({\cal G}_{N_k}^0(u_{k-1}))^c   \bigcup 
\Big([0,\e_0] \times {\bar {\cal G}}^c \Big)\Big| \nonumber \\
& \leq & \sum_{k \geq  1 } | {\mathtt G}_{N_k}^c(u_{k-1}) | +
\sum_{k \geq  1 }  |({\cal G}_{N_k}^0(u_{k-1}))^c| +   \e_0 |{\bar {\cal G}}^c | \nonumber \\
& \stackrel{(\ref{measGN0}), (\ref{measBN0}), (\ref{tautau0}), (\ref{calGg})} \leq & 
C \e_0 \sum_{k \geq  1 } N_k^{-1}   
+  C \e_0 \g  \leq C \e_0 (N_0^{-1}  +  \g) 
\stackrel{(\ref{N0ge})} \leq C \e_0^{1 + \a } \nonumber 
\end{eqnarray}
implying (\ref{Cmeas}). 

Theorem \eqref{thm:main}  is proved with $ s(d,\nu) := s_1 $
defined in \eqref{Sgr} and $ q(d,\nu) := S + 3 $, see \eqref{defk}. 

\subsubsection*{Regularity}

Finally, we prove that, if $ V,f,g, $ are  $ C^\infty $ then the solution   $ u(\e,\l) $ 
is in $ C^\infty (\T^d \times \T^\nu ) $.
The argument is the one of Theorem 3 in \cite{BBP}.
The main point is the proof of the following lemma which gives an 
a-priori bound for the divergence of the Sobolev high norms of the approximate solutions $ u_n $,
extending property $(S5)_n $.
Its proof requires only small modifications in Lemmata \ref{invLn+1},  \ref{htn+1}, \ref{S5n+1}.

\begin{lemma} $ \forall S' \geq S $,
\be\label{hS'}
\| u_n \|_{S'} \leq C(S') N_n^{ 2( \t' + \d s_1 + 1)} \, . 
\ee
\end{lemma}

\begin{pf}
First of all,  by the arguments of Lemma \ref{invLn+1}, we get,
the estimate 
\be\label{S'inte}
\norma {\cal L}_{n+1}^{-1}(u_n) \norma_{S'} \leq C(S') \Big(N_{n+1}^{\t'+\d S'} +  N_{n+1}^{\t'} \| u_n\|_{S'} \Big) \, .
\ee
Note that the multiscale Proposition \ref{propinv} is valid for any $ S' > s_1 $, see \eqref{s1}. 
It requires also the condition $ N \geq N_0(\Upsilon, S') $ which is verified for
$ N = N_n $ with $ n \geq n_0(S') $ large enough.

Then, 
following the proof of Lemma \ref{htn+1} we obtain
\begin{eqnarray}
\| {\wtilde h}_{n+1} \|_{S'} 
& \leq &  
N_{n+1}^{\t' + \d s_1 } \Big(\|  r_n \|_{S'} + \| R_n( {\wtilde h}_{n+1}) \|_{S'} \Big) \nonumber \\
& + &
C(S')  \Big(N_{n+1}^{\t'+\d S'} +  N_{n+1}^{\t'} \| u_n\|_{S'} \Big) 
\Big(\|  r_n \|_{s_1} + \| R_n ({\wtilde h}_{n+1}) \|_{s_1} \Big) \, . \label{hn+1S'}
\end{eqnarray} 
We also have the analogue of \eqref{pezz1}-\eqref{pezz2}, namely
$$
\| r_n \|_{S'} +  \| R_n({\wtilde h}_{n+1}) \|_{S'}  \leq C(S') (\| u_n\|_{S'}  +  \e \rho_{n+1}  \| {\wtilde h}_{n+1} \|_{S'}) \, , 
$$
$$
\|  r_n \|_{s_1} + \| R_n ({\wtilde h}_{n+1}) \|_{s_1} \leq C(S')  N_n^{-(S' - s_1)} \| u_n\|_{S'} \, ,
$$
and, by \eqref{hn+1S'}, we deduce 
the analogue of \eqref{quellaS'}, namely 
\be\label{ultiS'}
\| {\wtilde h}_{n+1} \|_{S'} \leq 
C(S') N_{n+1}^{\t'+ \d s_1}\| u_n \|_{S'} + 
C(S') N_{n+1}^{\t'}  N_n^{-(S'-s_1)} \| u_n \|_{S'}^2  + 
\e C(S')  N_{n+1}^{\t'+\d s_1} \rho_{n+1} \| {\wtilde h}_{n+1} \|_{S'} \, .
\ee 
For $ n \geq n_0(S') $ large enough, 
$$
\e C(S') N_{n+1}^{\t'+\d s_1}  \rho_{n+1} \stackrel{\eqref{defrhon+1} }
 = \e C(S') N_{n+1}^{\t'+\d s_1 - \s - 1} \stackrel{\eqref{def:sigma} } \leq \frac12  
$$
and \eqref{ultiS'},  \eqref{Sgr} imply the analogue of \eqref{halta1}, namely
\be\label{penul1}
\| {\wtilde h}_{n+1}  \|_{S'} \leq  
K(S') N_{n+1}^{\t'+ \d s_1}\| u_n \|_{S'} + 
K(S') N_{n+1}^{\t'}  N_n^{-(S'-s_1)} \| u_n \|_{S'}^2  \, .
\ee
Of course, $ h_{n+1} $ defined in \eqref{defih} satisfies \eqref{penul1} as well.
Therefore, as in Lemma \ref{S5n+1},
$$
\| u_{n+1} \|_{S'} \leq \| u_n \|_{S'} +  \| {h}_{n+1} \|_{S'} \leq
2 K(S') N_{n+1}^{\t'+ \d s_1}\| u_n \|_{S'} + 
 K(S') N_{n+1}^{\t'}  N_n^{-(S'-s_1)} \| u_n \|_{S'}^2 
$$
and we deduce that the sequence $ \| u_{n+1} \|_{S'} N_{n+1}^{- 2(\t'+ \d s_1 + 1 )} $ is bounded,
i.e. \eqref{hS'}. 
\end{pf}

By \eqref{hS'} 
we deduce
\be\label{fine3}
 \| h_n \|_{S'} \leq K(S') N_n^{2(\t_1 + \d s_1 + 1)} \, . 
\ee
Now, consider {\it any} $ s  > s_1 $ and write $ s := (1-t)s_1 + t S' $ where $ S' > s $, $ t \in (0,1)$. 
By interpolation 
\be\label{hns'}
\| h_n \|_{s} \leq  K(s_1,S') \| h_n \|_{s_1}^{1-t}  \| h_n \|_{S'}^{t}  
\stackrel{\eqref{defrhon+1}, \eqref{fine3}}\leq  K(S') N_n^{-(\s+1)(1-t) }  N_{n}^{\a t}  = 
K(S') N_n^{-1} 
\ee
having set $ \a := 2(\t_1 + \d s_1 + 1) $, and choosing  $ S' $ (large) such that 
$$ 
t = \frac{s - s_1}{S' - s_1} = \frac{\sigma + 2}{\s + 1 + \a } \, . 
$$
In conclusion, \eqref{hns'} implies that $ \sum_n \| h_n \|_s < + \infty $ and so 
$ u(\e , \l) \in {\bf H}^s $, for {\it any} $ s $.

\noindent
Massimiliano Berti, Dipartimento di Matematica e Applicazioni ``R. Caccioppoli",
Universit\`a degli Studi Napoli Federico II,  Via Cintia, Monte S. Angelo, 
I-80126, Napoli, Italy,  {\tt m.berti@unina.it}.
\\[2mm]
Philippe  Bolle, Universit\'e
d'Avignon et des Pays de Vaucluse, Laboratoire de math\'ematiques d'Avignon (EA 2151), 
F-84018 Avignon, France, {\tt philippe.bolle@univ-avignon.fr}.
\\[2mm]
\indent
This research was supported by the European Research Council under FP7 ``New Connections between dynamical systems
and Hamiltonian PDEs with small divisors phenomena".
\end{document}